%% file: wstart.tex
\documentclass[11pt,letterpaper]{article}
\input{sections/wstart_preamble}
\usepackage[margin=1in]{geometry}

\makeatletter
\def\blfootnote{\gdef\@thefnmark{}\@footnotetext}
\makeatother

\begin{document}

    \title{Faster high-accuracy log-concave sampling \\ via algorithmic warm starts}
    
        \author{
		Jason M.\ Altschuler \\
		NYU\\
            \texttt{ja4775@nyu.edu}
		\and
		Sinho Chewi \\
		MIT \\
		\texttt{schewi@mit.edu}
	}
	\date{\today}
	\maketitle

\input{sections/abstract}

   \newpage
	\setcounter{tocdepth}{2}
	\tableofcontents	
	\normalsize
	\newpage

\input{sections/intro}

\input{sections/prelim}

\input{sections/pabi}

\input{sections/ulmc}

\input{sections/hiacc}

\input{sections/discussion}

   \paragraph*{Acknowledgments.} We thank Mufan (Bill) Li and Matthew Zhang for many insightful conversations regarding hypocoercivity. JMA acknowledges funding from an NYU Faculty Fellowship. SC acknowledges funding from NSF TRIPODS program (award DMS-2022448).

    \newpage{}
	\appendix
	\input{sections/app_background}
	\input{sections/app_pabi}
	\input{sections/app_loacc}
	\input{sections/app_hiacc}

	\small
	\addcontentsline{toc}{section}{References}
	%\bibliographystyle{plainnat}
	%\bibliography{wstart}{}
	\printbibliography{}

\end{document}

%% file: sections/wstart_preamble.tex
\usepackage[utf8]{inputenc}

\usepackage{amsmath, amsthm}
\usepackage{MnSymbol} 
\usepackage{graphicx,color}
\usepackage[hyphens]{url}
\usepackage{dsfont}
\usepackage{booktabs}
\usepackage[normalem]{ulem}
\usepackage{mathtools}
\usepackage{hyperref}
\usepackage[nameinlink,capitalize]{cleveref}
\usepackage{multirow}
\usepackage{algorithm}
\usepackage[noend]{algpseudocode}
\usepackage{tikz} %

\usepackage{soul} %

\usepackage[giveninits=true, maxnames=10, style=alphabetic, sorting=nyt]{biblatex}

\usepackage{subcaption}

\usepackage{tabularx}
\newcolumntype{L}{>{\raggedright\arraybackslash}X}

\DeclareMathAlphabet\mathbb{U}{msb}{m}{n}

\numberwithin{equation}{section}
\newtheorem{theorem}{Theorem}[section]
\newtheorem{cor}[theorem]{Corollary}
\newtheorem{lemma}[theorem]{Lemma}
\newtheorem{remark}[theorem]{Remark}
\newtheorem{prop}[theorem]{Proposition}

\newtheorem{conj}[theorem]{Conjecture}
\newtheorem{defin}[theorem]{Definition}

\newcommand{\cC}{\mathscr{C}}

\newcommand{\cM}{\mathcal{M}}
\newcommand{\cN}{\mathcal{N}}

\newcommand{\cP}{\mathcal{P}}

\newcommand{\cR}{\mathcal{R}}

\newcommand{\R}{\mathbb{R}}

\newcommand{\N}{\mathbb{N}}

\newcommand*{\E}{\mathbb{E}}

\DeclareMathOperator*{\Tr}{tr}

\newcommand*{\Otilde}{\widetilde{O}}

\newcommand*{\polylog}{\mathrm{polylog}}

\newcommand*{\eps}{\varepsilon}

\newcommand*{\defeq}{\coloneqq}

\DeclareMathOperator*{\argmin}{arg\,min}

\renewcommand{\leq}{\leqslant}
\renewcommand{\geq}{\geqslant}
\DeclareMathOperator{\half}{\frac{1}{2}}

\providecommand{\abs}[1]{\left\lvert#1\right\rvert}

\providecommand{\norm}[1]{\lVert{#1}\rVert}

\newcommand{\sig}{\sigma}

\newcommand{\TV}{\mathsf{TV}}
\newcommand{\KL}{\mathsf{KL}}

\newcommand{\xopt}{x^*}

\usepackage{bbm}
\usepackage{blkarray}
\usepackage{cancel}
\usepackage{enumitem}
\usepackage{float}
\usepackage{mathrsfs}
\usepackage{microtype}
\usepackage{ragged2e}
\usepackage{sectsty}
\usepackage{thmtools}
\usepackage[Symbolsmallscale]{upgreek}
\usepackage{array}

\usetikzlibrary{cd}

\hypersetup{citecolor=red, colorlinks=true, linkcolor=blue}

\makeatletter
\patchcmd{\@counteralias}
{\@ifdefinable{c@#1}}
{\expandafter\@ifdefinable\csname c@#1\endcsname}
{}{}
\makeatother

\makeatletter
\newcommand{\opnorm}{\@ifstar\@opnorms\@opnorm}
\newcommand{\@opnorms}[1]{%
	\left|\mkern-1.5mu\left|\mkern-1.5mu\left|
	#1
	\right|\mkern-1.5mu\right|\mkern-1.5mu\right|
}
\newcommand{\@opnorm}[2][]{%
	\mathopen{#1|\mkern-1.5mu#1|\mkern-1.5mu#1|}
	#2
	\mathclose{#1|\mkern-1.5mu#1|\mkern-1.5mu#1|}
}
\makeatother

\newcommand{\PreserveBackslash}[1]{\let\temp=\\#1\let\\=\temp}
\newcolumntype{C}[1]{>{\PreserveBackslash\centering}p{#1}}

\newcommand\bs[1]{\boldsymbol{#1}}
\newcommand\mb[1]{\mathbf{#1}}

\newcommand\mc[1]{\mathcal{#1}}
\newcommand\mf[1]{\mathfrak{#1}}

\newcommand\mss[1]{\mathsf{#1}}
\newcommand\msf[1]{\mathsf{#1}}
\newcommand\on[1]{\operatorname{#1}}

\definecolor{MITBrown}{RGB}{164, 31, 50}

\DeclareMathOperator\esssup{ess\,sup}

\DeclareMathOperator\law{law}

\DeclareMathOperator\one{\mathbbm{1}}

\DeclareMathOperator\prox{prox}

\DeclareMathOperator\var{var}

\newcommand{\D}{\mathrm{d}}

\newcommand\deq{\coloneqq}

\newcommand\mmid{\mathbin{\|}}

\newcommand\T{\mathsf{T}}

\newcommand*\circled[1]{\tikz[baseline=(char.base)]{
		\node[shape=circle,draw,inner sep=2pt] (char) {#1};}}

\newcommand{\hyper}{hyperequilibration}
\newcommand{\epsrgo}{\eps_{\msf{RGO}}}
\newcommand{\Nprox}{N_{\msf{prox}}}
\newcommand{\Nrgo}{N_{\msf{RGO}}}
\newcommand{\patha}{\widehat{\bs \Pi}}
\newcommand{\pathb}{\bs \Pi}

%% file: sections/abstract.tex
\begin{abstract}
    Understanding the complexity of sampling from a strongly log-concave and log-smooth distribution $\pi$ on $\R^d$ to high accuracy is a fundamental problem, both from a practical and theoretical standpoint. In practice, high-accuracy samplers such as the classical Metropolis-adjusted Langevin algorithm (MALA) remain the de facto gold standard; and in theory, via the proximal sampler reduction, it is understood that such samplers are key for sampling even beyond log-concavity (in particular, for distributions satisfying isoperimetric assumptions). 
    
    \par In this work, we improve the dimension dependence of this sampling problem to $\widetilde O(d^{1/2})$, whereas the previous best result for MALA was $\widetilde O(d)$. This closes the long line of work on the complexity of MALA, and moreover leads to state-of-the-art guarantees for high-accuracy sampling under strong log-concavity and beyond (thanks to the aforementioned reduction). 
    \par Our starting point is that the complexity of MALA improves to $\widetilde O(d^{1/2})$, but only under a \emph{warm start} (an initialization with constant R\'enyi divergence w.r.t.\ $\pi$). Previous algorithms took much longer to find a warm start than to use it, and closing this gap has remained an important open problem in the field.
	Our main technical contribution settles this problem by establishing the first $\Otilde(d^{1/2})$ R\'enyi mixing rates for the discretized underdamped Langevin diffusion. For this, we develop new differential-privacy-inspired techniques based on R\'enyi divergences with Orlicz--Wasserstein shifts, which allow us to
    sidestep longstanding challenges for proving fast convergence of hypocoercive differential equations. 
\end{abstract}

%% file: sections/intro.tex
\section{Introduction}\label{sec:intro}

We consider the algorithmic problem of efficiently sampling from a high-dimensional probability distribution $\pi$ on $\R^d$. Due to the many important applications of sampling throughout applied mathematics, engineering, and statistics, significant research effort has been devoted to designing fast sampling algorithms and analyzing their convergence rates. We refer to the book draft~\cite{chewibook} for a recent exposition of the extensive literature and its history.

Yet, despite several decades of progress, many fundamental theoretical questions remain open about the complexity of sampling. Arguably one of the foremost questions in this field is:
\[
	\text{\emph{What is the first-order query complexity for sampling from }} \pi\text{\emph{?}}
\]
Recall that a first-order query refers to accessing $f(x)$ and $\nabla f(x)$ at a query point $x \in \R^d$, where $f$ denotes the negative log-density of $\pi \propto \exp(-f)$ (up to an additive normalization constant). This well-studied notion of a first-order query is inspired on one hand by the fact that such queries do not require knowledge of the normalization constant $\int \exp(-f)$ and thus are readily available in many practical applications, and on the other hand also inspired by the analogous and influential theory of complexity for convex optimization~\cite{nemirovskij1983problem}.

This problem of determining the query complexity for sampling has remained open even for the canonical and seemingly simple class of strongly log-concave and log-smooth (in brief, ``well-conditioned'') distributions $\pi$, let alone in more complicated settings. It is worth emphasizing that this state of affairs for \emph{sampling} is in sharp contrast to that for \emph{optimization}---indeed, the analogous query complexity questions for convex optimization were solved long ago in celebrated results from the 1980s~\cite{nemirovskij1983problem, nesterov2018cvxopt}.

Within the literature, it is of central interest to understand this complexity question in the \emph{high-accuracy regime}\footnote{Throughout, we use the standard terminology \emph{low accuracy} to refer to complexity results which scale polynomially in $1/\varepsilon$, and the term \emph{high accuracy} for results which scale polylogarithmically in $1/\varepsilon$; here, $\varepsilon$ is the desired target accuracy. These two regimes require different algorithms and analyses, as explained in the sequel.
% These two regimes comprise different classes of samplers which in turn generally require different sets of techniques for their analysis, as explained in the sequel.
}, since classical high-accuracy samplers such as the Metropolis-adjusted Langevin algorithm (MALA) and the Metropolized Hamiltonian Monte Carlo algorithm (MHMC) remain the de facto gold standard in practice. Yet the complexity for this high-accuracy setting has been particularly difficult to pin down, as we explain shortly.

The purpose of this paper is to develop faster high-accuracy samplers, and in doing so move towards a better understanding of the first-order complexity of sampling.
For simplicity of exposition, let us presently assume that $\pi$ is well-conditioned, since by the proximal reduction framework~\cite{leeshentian2021rgo,chenetal2022proximalsampler}, it is known that improvements to the complexity of well-conditioned sampling lead to improvements in more general settings such as when $\pi$ is (non-strongly) log-concave, or even non-log-concave but satisfies standard isoperimetric assumptions such as the log-Sobolev or Poincar\'e inequality. (Indeed, our results improve upon the state-of-the-art for all these settings.)

\paragraph*{The gap between low-accuracy and high-accuracy samplers.} A central motivation of this paper is the large gap between (our current understanding of) the complexity of low-accuracy samplers and high-accuracy samplers. To explain this gap, let us briefly provide relevant background on both classes of algorithms. 

Low-accuracy samplers arise as discretizations of stochastic processes with stationary distribution $\pi$, such as the Langevin diffusion~\parencite[the sampling analog of the gradient flow, see][]{jordan1998variational, wibisono2018sampling} or the underdamped Langevin diffusion~\parencite[the sampling analog of the accelerated gradient flow, see][]{ma2021there}.
Once discretized, however, the resulting discrete-time Markov chain is typically \emph{biased}, i.e., its stationary distribution is no longer equal to $\pi$. In order to control the size of the bias, the step size of the algorithm is chosen to scale polynomially with $\varepsilon$, and hence the overall running time scales polynomially with $1/\varepsilon$.
Despite this drawback, the discretization analysis is by now well-understood, with state-of-the-art results achieving a complexity of $\widetilde O(d^{1/3}/\varepsilon^{2/3})$~\parencite{shen2019randomized, foslyoobe21shiftedode, boumar22hamiltonian}; see~\cite{caoluwang2021uldlowerbd} for a discussion of tightness. 

\par High-accuracy samplers, in contrast, 
are typically designed in such a way that there is no bias. This is achieved by, e.g., appending a Metropolis{--}Hastings filter to each step (see Appendix~\ref{app:hiacc:mala} for background). Common examples of these algorithms include MALA and MHMC, which are routinely deployed in large-scale applications and are the default implementations of sampling routines in many modern software packages~\cite{gelleeguo15stan, abadi2016tensorflow}. However, the filter which debiases the algorithm also greatly complicates the analysis, and thus far the best complexity result for these algorithms\footnote{We discuss the result of~\cite{luwang2022zigzag} for the zigzag sampler further in \S\ref{ssec:intro:prior}.} is $\widetilde O(d\log^{O(1)}(1/\varepsilon))$~\cite{dwivedi2018log, chenetal2020hmc, leeshentian2020gradientconcentration}.
Note that the dimension dependence of this result is substantially worse than what is known in the low-accuracy regime and is at odds with the popularity of high-accuracy samplers in practice.

\paragraph*{The mystery of warm starts.} A promising first step towards resolving this gap was put forth in~\cite{chewi2021optimal} and later refined in~\cite{wuschche2022minimaxmala}: when initialized from a \emph{warm start} (i.e., a measure $\mu_0$ with $\chi^2(\mu_0 \mmid \pi) \le O(1)$), the complexity of MALA improves to $\widetilde O(d^{1/2} \log^2(1/\varepsilon))$ since it can safely take much larger step sizes (of size $d^{-1/2}$ rather than $d^{-1}$). This raises the natural question: is the warm start condition merely an artefact of the analyses? Rather surprisingly, it was shown in~\cite{lee2021lower} that there exist bad initializations for MALA for which the dimension dependence is at least $\widetilde \Omega(d)$.
Taken together, these results show that the complexity of MALA fundamentally hinges on the warmness of its initialization.
\par The key question is thus: can such a warm start be obtained algorithmically? Or more precisely:
\begin{center}
\emph{Is there an algorithm which makes} $\widetilde O(d^{1/2})$ \emph{queries to a first-order oracle for} $f$ \\ \emph{and outputs a measure} $\mu_0$ \emph{with} $\chi^2(\mu_0 \mmid \pi) \le O(1)$?
\end{center}
The requirement that the algorithm makes $\widetilde O(d^{1/2})$ queries is essential, else the cost of obtaining the warm start dominates the subsequent cost of running MALA. Yet this was the state of affairs---previously, the fastest algorithms took significantly longer to produce a warm start than to actually use it, defeating the purpose of the warm start. Resolving this discrepancy has been posed as an important question in many papers, e.g.,~\cite{chewi2021optimal, lee2021lower, chewi2021analysis, luwang2022zigzag, wuschche2022minimaxmala}.

\par The main challenge for answering this warm start question is that the chi-squared divergence
is quite a strong performance metric. (We emphasize that it is essential to obtain the warm start in the chi-squared divergence, or more generally in a R\'enyi divergence $\cR_q$ of order $q > 1$, rather than other common metrics such as total variation, Wasserstein, or KL divergence; see \S\ref{ssec:intro:tech} for an in-depth discussion.)
The aforementioned results in the low-accuracy regime fall short of achieving this goal, since they only hold in the Wasserstein metric (for which standard coupling arguments are readily available). Despite significant effort, the best known guarantee for producing a warm start---achieved by the Langevin Monte Carlo (LMC) algorithm~\cite{chewi2021analysis}---is far too costly as it requires $\widetilde O(d)$ queries, which defeats the purpose of the warm start.

\par Towards this hope of algorithmic warm starts, \cite{wuschche2022minimaxmala} made the promising empirical observation that MALA mixes much faster if it is initialized at the output of the \emph{underdamped} Langevin Monte Carlo (ULMC) algorithm. %
However, they left open the question of rigorously proving that this yields a warm start. While it is widely believed that ULMC is substantially faster than LMC\@, the previous best results for computing a warm start with ULMC had dimension dependence $\Otilde(d^{5/2})$ (implicit from~\cite{GaneshT20}) or very recently $\Otilde(d^{2})$ (implicit from~\cite{Matt23girsanov}), see the prior work section \S\ref{ssec:intro:prior} for details. 
We emphasize that this dimension dependence is not only a far cry from the elusive $\Otilde(d^{1/2})$ goal, but moreover
is even worse than known results for the simpler LMC algorithm. Unfortunately, any improvement to these ULMC warm start bounds appears to require overcoming fundamental difficulties with studying hypocoercive differential equations which remain unsolved today, despite being the focus of intensive research activity within the PDE community since the work of Kolmogorov~\cite{Kol34}.
For a further discussion of these technical obstacles, see \S\ref{ssec:intro:tech}.

\subsection{Contributions}\label{ssec:intro:cont}

In this paper, we develop techniques which
bypass longstanding challenges for analyzing hypocoercive dynamics, thereby establishing the first $\Otilde(d^{1/2})$ R\'enyi mixing results for ULMC\@.
This resolves the aforementioned warm start conjecture, which has been raised in a number of prior works, e.g.,~\cite{chewi2021optimal, lee2021lower, chewi2021analysis, luwang2022zigzag, wuschche2022minimaxmala}.
As discussed above, this enables us to design significantly faster high-accuracy samplers---both for the log-concave setting and far beyond. Finally, this also closes the long line of work devoted to understanding the complexity of MALA (see Table~\ref{tab:mala-progress}).
We present our results in more detail below, and then discuss our new techniques in \S\ref{ssec:intro:tech}.

\paragraph{Result 1: Algorithmic warm starts via ULMC\@.} Our first main result is an improvement of the state-of-the-art R\'enyi mixing bounds for ULMC from $\Otilde(d^2)$ to $\Otilde(d^{1/2})$. This resolves the warm start question in the affirmative.
We remark that although the warm start problem was stated above for $\chi^2$ convergence, our result actually holds more generally for R\'enyi divergences $\cR_q$ of any order $q\ge 1$, and thus we state it as such. (For the purpose of warm starts, it suffices to take $q=2$ since $\chi^2 = \exp(\cR_2) -1$ is of constant size when $\cR_2$ is.) Below, $\alpha$ and $\beta$ denote the strong log-concavity and log-smoothness bounds; their ratio $\kappa := \beta/\alpha$ is the condition number.

\begin{theorem}[R\'enyi guarantees for ULMC; informal version of Theorem~\ref{thm:loacc}]\label{thm:intro:ulmc}
    Consider the class of densities of the form $\pi \propto \exp(-f)$ on $\R^d$, where $\alpha I \preceq \nabla^2 f \preceq \beta I$ and $\kappa \deq \beta/\alpha < \infty$.
    The ULMC algorithm outputs a measure $\mu$ satisfying $\cR_q(\mu \mmid \pi) \le \varepsilon^2$ using $\widetilde O(\kappa^{3/2} d^{1/2} q^{1/2}/\varepsilon)$ first-order queries.
\end{theorem}

As we detail in \S\ref{ssec:intro:tech}, the main barrier to obtaining this result is that the underdamped Langevin dynamics falls within a class of PDEs known as hypocoercive equations, for which fundamental questions remain unresolved.\footnote{Implications of our techniques for the analysis of hypocoercive diffusions are explored in shortly forthcoming work.}

\begin{table}
    \centering
    \begin{tabular}{ccc}
    \textbf{Reference} & \textbf{Complexity} & \textbf{Algorithmically Achievable?} \\
         \cite{dwivedi2018log} & $\kappa d + \kappa^{3/2} d^{1/2}$ & No, requires a warm start \\
        \cite{chenetal2020hmc} & $\kappa d + \kappa^{3/2} d^{1/2}$ & Yes \\
        \cite{leeshentian2020gradientconcentration} & $\kappa d$ & Yes \\
        \cite{chewi2021optimal} & $\kappa^{3/2} d^{1/2}$ & No, requires a warm start \\
        \cite{wuschche2022minimaxmala} & $\kappa d^{1/2}$ & No, requires a warm start \\
        \textbf{\textcolor{magenta}{Theorem~\ref{thm:main-slc}}} & \textbf{\textcolor{magenta}{$\kappa d^{1/2}$}} & \textbf{\textcolor{magenta}{Yes, warm start provided in Theorem~\ref{thm:loacc}}}
    \end{tabular}
    \caption{\small 
    This table summarizes the community's progress towards non-asymptotic complexity bounds for MALA\@; 
    the asymptotic study of MALA is much more classical, and dates back to at least~\cite{robertsrosenthal1998optimalscaling}. The complexity bounds displayed are upper bounds; for brevity, we hide logarithmic factors as well as the dependence on $\eps$ since all results scale polylogarithmically in $1/\eps$.
    As discussed in the main text, Theorem~\ref{thm:main-slc} completes our understanding of MALA due to matching lower bounds in~\cite{chewi2021optimal, lee2021lower, wuschche2022minimaxmala}.
    }
    \label{tab:mala-progress}
\end{table}

\paragraph{Result 2: Faster high-accuracy log-concave sampling.} Theorem~\ref{thm:intro:ulmc} provides the first algorithm for computing warm starts that is not significantly slower than the use of the warm start. This enables us to exploit, for the first time, the recent breakthroughs on MALA~\cite{chewi2021optimal,wuschche2022minimaxmala} which improve the complexity of MALA from  $\Otilde(d)$ to $\Otilde(d^{1/2})$ from a warm start.\footnote{We remark that all of our results could replace MALA with the zigzag algorithm~\cite{luwang2022zigzag}. Indeed, the zigzag sampler has the same key issue as MALA: it requires a warm start in chi-squared divergence for the known $d^{1/2}$ mixing result to apply. However, we focus on MALA because MALA's robust empirical performance has made it a central focus of study in the MCMC literature for nearly three decades~\cite{besagetal1995bayesiancomp}.\label{fn:zigzzag}} By combining this with additional algorithmic tools for improving the dependence on the condition number, we obtain our second main result, which substantially advances the state-of-the-art for high-accuracy log-concave sampling.

\begin{theorem}[High-accuracy log-concave sampling; informal version of Theorem~\ref{thm:main-slc}]\label{thm:intro:main}
    Consider the class of densities of the form $\pi \propto \exp(-f)$ on $\R^d$, where $\alpha I \preceq \nabla^2 f \preceq \beta I$ and $\kappa \deq \beta/\alpha$.
    There is an algorithm which outputs a sample with law $\mu$ satisfying $\msf d(\mu,\pi) \le \varepsilon$, for any performance metric $\msf d \in \{\msf{TV}, \sqrt{\msf{KL}}, \sqrt{\chi^2}, \sqrt\alpha\,W_2\}$, after making $\kappa d^{1/2}\log^{O(1)}(\kappa d/\varepsilon)$ first-order queries.
    % in expectation.
\end{theorem}

The algorithmic warm start result of Theorem~\ref{thm:intro:ulmc} confirms the aforementioned empirical conjecture of~\cite{wuschche2022minimaxmala} and provides the final missing piece in our understanding of the complexity of MALA, closing the line of work developed in~\cite{robertsrosenthal1998optimalscaling,dwivedi2018log, chenetal2020hmc, leeshentian2020gradientconcentration, chewi2021optimal, lee2021lower, wuschche2022minimaxmala} (Table~\ref{tab:mala-progress}). Indeed, due to matching lower bounds in~\cite{chewi2021optimal, lee2021lower, wuschche2022minimaxmala}, the complexities $\Otilde(\kappa d)$ and $\Otilde(\kappa d^{1/2})$ with or without a warm start are known to be tight, and hence the key remaining question was whether a warm start is actually efficiently computable.

\par The complexity in Theorem~\ref{thm:intro:main} constitutes a natural barrier for high-accuracy sampling. 
Indeed, regarding the dimension dependence, any further progress beyond $\widetilde O(d^{1/2})$ would seem to require completely different algorithms---both for obtaining a warm start and also for exploiting a warm start. For example, the $\Otilde(d^{1/2})$ complexity of MALA is unimprovable even under arbitrarily warm starts~\cite{chewi2021optimal,wuschche2022minimaxmala}. And regarding the condition number dependence, any further progress beyond $\widetilde O(\kappa)$ in the high-dimensional regime\footnote{Analogous to classical optimization results, there are sampling algorithms which achieve logarithmic dependence on $\kappa$ at the expense of larger polynomial dependence on $d$. The open question mentioned here is really: can one improve the condition dependence beyond near-linear while also maintaining comparable dimension dependence?} would constitute a major breakthrough in the complexity of sampling since it is currently unknown whether an acceleration phenomenon holds in the sampling context.

\par More broadly, our result provides evidence of the potential for designing faster high-accuracy samplers by combining low-accuracy samplers for computing a warm start, together with improved high-accuracy mixing from the warm start. We believe that this research program may be crucial for future progress in high-accuracy sampling, since faster mixing from a warm start seems likely to hold for other Metropolized algorithms. See \S\ref{sec:discussion} for further discussion in this direction.

\paragraph{Result 3: Faster high-accuracy sampling beyond log-concavity.} 
High-accuracy log-concave sampling is the key to obtaining state-of-the-art complexity results for a wide class of distributions beyond log-concavity.
This is achieved
by using our faster log-concave sampler in Theorem~\ref{thm:intro:main} to improve the per-iteration complexity of the proximal sampler~\cite{leeshentian2021rgo, chenetal2022proximalsampler}. This approach is overviewed in the techniques section \S\ref{ssec:intro:tech}, and leads to the following result.

\begin{cor}[Sampling from other classes of distributions; informal version of results in \S\ref{ssec:hiacc-ext}]
    For each of the following classes of distributions, we obtain complexity bounds which improve by a factor of $d^{1/2}$ over the state-of-the-art results in~\cite{chenetal2022proximalsampler}:
    \begin{itemize}
        \item $\pi$ is log-smooth and weakly log-concave.
        \item $\pi$ is log-smooth and satisfies a log-Sobolev inequality.
        \item $\pi$ is log-smooth and satisfies a Poincar\'e inequality.
    \end{itemize}
\end{cor}

The latter two assumptions of log-Sobolev and Poincar\'e---called \emph{functional inequalities}---capture strictly richer classes of target distributions than strong-log-concavity.
There are two major motivations for studying the complexity of sampling in this setting. First, functional inequalities are quite flexible, as they are preserved under common operations such as bounded perturbations and Lipschitz mappings (see \S\ref{ssec:prelim-functional}).
Consequently, they often capture the breadth of settings encountered in practice, including non-log-concave settings.
Second, these functional inequalities classically imply convergence of diffusions in continuous time,
making them natural assumptions under which to study the corresponding discretizations.

Despite the appeal of this program, proving sampling guarantees under functional inequalities introduces a number of additional technical complications and was only accomplished recently, starting with~\cite{VempalaW19} and continued in the works~\cite{wibisono2019proximal, li2020riemannian, ma2021there, chewi2021analysis}.
Our result continues this line of work, and in particular highlights the use of high-accuracy samplers for well-conditioned distributions as a powerful algorithmic tool for the broader problem of sampling under isoperimetry.

\subsection{Challenges and techniques}\label{ssec:intro:tech}

\subsubsection{Challenges for warm starts: R\'enyi divergence and hypocoercivity} 

\paragraph*{Why R\'enyi?} 
To explain what properties are needed for a warm start requires first explaining why a warm start helps. Briefly, the complexity of MALA is governed by the largest possible step size for which the algorithm still accepts a reasonable fraction of the proposals (see Appendix~\ref{app:hiacc:mala} for background on MALA).
The basic reason why we might expect to improve the complexity of MALA from $\Otilde(d)$ to $\Otilde(d^{1/2})$
is that
\emph{at the stationary distribution $\pi$}, 
the step size can be increased significantly from $d^{-1}$ to $d^{-1/2}$ while keeping the acceptance probability high. 
More precisely, with step size $d^{-1/2}$,
the acceptance probability is large from a typical point from $\pi$; however, it
can be exponentially small in regions that are atypical (i.e., exponentially rare under $\pi$). The existence of such regions implies that there are ``bottlenecks'' in the state space which take exponentially long to traverse. The role of a warm start initialization is to avoid such bottlenecks.

In other words, a key property that a warm start $\mu_0$ must satisfy is that if $\pi$ assigns exponentially small probability to an event, then so must $\mu_0$. Crucially, this property does not hold if $\mu_0$ is only known to be close to $\pi$ in common probability metrics such as total variation, Wasserstein, or KL divergence---but this property \emph{does} hold if $\mu_0$ is close to $\pi$ in the chi-squared divergence, or more generally any R\'enyi divergence $\cR_q$ of order $q > 1$.\footnote{This is the same reason why differential privacy requires guarantees in R\'enyi divergences~\cite{mironov2017renyi}.}

\paragraph*{The key to warm starts: low-accuracy algorithms.}
In the preceding discussion, taking large step sizes from a non-warm initialization was problematic due to the rejections in the Metropolis--Hastings filter step. A natural idea, then, is to remove the filter for the initial stage of the algorithm and later reinstate it when the law of the iterate is closer to the target $\pi$.
Since the proposal of MALA is just one step of the LMC algorithm, this amounts to using LMC to procure the warm start.
More generally, we can consider using any low-accuracy sampler as our warm start algorithm, and indeed, as we discuss next, it will be crucial to consider ULMC instead of LMC in order to achieve the desired $\Otilde(d^{1/2})$ dimension dependence. 
% Recall that from the preceding paragraph that throughout this discussion, we must work in the chi-squared or R\'enyi divergence.

At a high level, if we discretize a diffusion with step size $h$ for continuous time $T$, then the total number of iterations is $N = T/h$. In order to understand the dimension dependence of the algorithm, one must therefore understand both $h$ and $T$. These two terms reflect two distinct aspects of mixing analysis: the discretization bias and the convergence time.

The first part---the discretization bias---is now relatively well-understood (see the prior work discussions in \S\ref{ssec:intro:prior}), even for the chi-squared divergence and more general R\'enyi divergences. In particular, it is known that the R\'enyi bias of LMC is controlled for step sizes $h \lesssim 1/(dT)$, and the R\'enyi bias of ULMC is controlled for step sizes $h\lesssim 1/\sqrt{dT}$. (In fact, we streamline arguments in the literature in order to provide a shorter and simpler proof of this in Appendix~\ref{app:loacc-girsanov}.)
Since the Langevin diffusion does not reach approximate stationarity until time $T\ge \Omega(\log d)$, it follows that LMC requires at least $N = T/h = \widetilde{\Omega}(d)$ iterations, which is too slow for our purposes.

ULMC is more promising, as the discretization bounds lead to iteration complexity bounds of $N = T/h = d^{1/2} T^{3/2}$. However, in order to reach our warm start goal of $N = \Otilde(d^{1/2})$, this means that the convergence time $T$ must be nearly dimension-free, i.e., of size $\Otilde(1)$.

\paragraph*{Why are nearly dimension-free convergence rates possible in continuous time?}
Since the R\'enyi divergence to $\pi$ initially scales as $\widetilde \Theta(d)$, in order to obtain nearly dimension-free bounds on $T$, we require the diffusion to converge to stationarity in R\'enyi divergence with an \emph{exponential} rate.
This is a strong property of the diffusion, which we call \emph{\hyper{}}.

Hyperequilibration was not even known for the simpler (standard, overdamped) Langevin diffusion (LD) until quite recently~\cite{caolulu2019renyi, VempalaW19}.
While a spectral gap for LD (or equivalently, a Poincar\'e inequality for $\pi$) classically implies exponential decay of the chi-squared divergence, this is far weaker than \hyper{}. Indeed, \hyper{} requires exponential decay of $\cR_2$, which amounts to \emph{doubly} exponential decay of the chi-squared divergence, 
since $\cR_2 = \log(\chi^2 + 1)$.
% since $\chi^2 = \exp(\cR_2) - 1$.
Under the stronger assumption of a log-Sobolev inequality for $\pi$, it is well-known that the KL divergence decays exponentially fast, but it was unclear that the same holds for the \emph{R\'enyi divergence} which, as discussed above, is crucial for warm starts.
It was only through the inspired semigroup calculations of~\cite{caolulu2019renyi, VempalaW19} that we now know this to be true, namely, a log-Sobolev inequality implies \hyper{} for LD\@.\footnote{This explains our choice of the terminology \emph{\hyper}: it is inspired by the analogy to the classical property of \emph{hypercontractivity}, which is equivalent to the logarithmic Sobolev inequality (LSI)~\cite{grosslsi}.}

Recall, though, that the LD incurs too much discretization bias.
To obtain sufficient control over both the discretization bias and the convergence time, we therefore need to establish \hyper{} for the \emph{underdamped} Langevin diffusion (ULD). However, this question brings us to longstanding challenges from the theory of hypocoercive PDEs.

\paragraph*{Hypocoercivity: a fundamental barrier for underdamped analysis.} 
To recap: for LD, we have exponential decay of the chi-squared divergence under a Poincar\'e inequality, exponential decay of the KL divergence under a log-Sobolev inequality, and finally \hyper{} under a log-Sobolev inequality. What, then, are the analogous results for ULD?
Since its introduction in the 1930s by Kolmogorov~\cite{Kol34}, the regularity and convergence of ULD have been the focus of intensive research.
It took nearly half a century to establish mixing~\cite{tropper1977ergodic}, and a further 30 years and Villani's ``slightly miraculous-looking computations''~\cite[pg.\ 42]{villani2009hypocoercivity} to prove exponential decay of the KL divergence under a log-Sobolev inequality.
Establishing \hyper{} for ULD remains out of reach for existing techniques.

The reason for this sudden jump in difficulty from the overdamped to the underdamped diffusions is
due to a fundamental issue: the \emph{degeneracy} of ULD\@.
In brief, whereas LD is driven by a full-dimensional Brownian motion, ULD is driven by a degenerate one which is only added to a subset of the coordinates.
For sampling purposes, this degeneracy is a desirable feature as it leads to smoother sample paths and smaller discretization error; however, this same degeneracy is also the source of deep questions in PDE theory which have motivated research in that field for nearly a century.
The key challenge here is that the standard tools of Markov semigroup theory---which provide the backbone of the analysis for LD---completely break down for ULD\@.
To address this difficulty, the theory of \emph{hypocoercivity}, inspired by H\"ormander's groundbreaking work on hypoellipticity~\cite{Hor67}, was laid down by Villani in the monograph~\cite{villani2009hypocoercivity} as a principled framework for the study of degenerate diffusions.
However, this is still a relatively nascent area of PDE and many important questions remain wide open; see the prior work in \S\ref{ssec:intro:prior} for further background.

\par In contrast, we note that it is well-known how to obtain fast rates of convergence in the Wasserstein metric via standard coupling arguments.
Consequently, the state-of-the-art $\Otilde(d^{1/2})$ guarantees for the ULMC algorithm hold in the Wasserstein metric or the KL divergence~\cite{cheng2018underdamped, shen2019randomized, dalalyanrioudurand2020underdamped, ma2021there,Matt23girsanov}, whereas for R\'enyi divergence bounds, it was previously unknown how to obtain rates which are better than even $\Otilde(d^2)$.

\subsubsection{Settling the warm start conjecture: regularization via privacy}

\paragraph*{Our approach to \hyper.} To settle the warm start conjecture, we adopt a fundamentally different perspective. Namely, instead of trying to directly establish \hyper{} via hypocoercivity techniques, we ask whether it can be deduced from simpler Wasserstein coupling arguments.
At the heart of this approach is the fact that diffusions often enjoy strong \emph{regularizing} properties, which allow for bounding stronger metrics (e.g., R\'enyi) in terms of weaker ones (e.g., Wasserstein).
Such regularization results are typically established for continuous-time diffusions via abstract calculus methods, such as the theory of Markov semigroups~\cite{bakry2014analysis}. However, as discussed above, these techniques do not extend to ULD due to the fundamental issue of degeneracy.

Our key insight is to prove a regularization result for the \emph{discrete-time} algorithm directly. This is enabled by the fact that although the noise added to each iteration of ULMC is \emph{nearly} degenerate---and indeed degenerates as the step size $h \searrow 0$, as it must because ULD is degenerate---this ULMC noise remains non-degenerate for any positive step size $h > 0$.
Hence, we can expect some mild amount of regularization for ULMC, a fact that we establish for the first time.
On a technical level, we accomplish this via a more sophisticated version of techniques from the differential privacy literature---namely, the shifted R\'enyi analysis---which we describe next.

\paragraph*{R\'enyi divergences with Orlicz--Wasserstein shifts.}
The regularization result we seek is of the following form: if we initialize two copies of our process of interest at the distributions $\mu_0$, $\nu_0$, and arrive at distributions $\mu_n$, $\nu_n$ respectively at iteration $n$, we wish to control $\cR_q(\mu_n \mmid \nu_n)$ in terms of an initial Wasserstein distance $W(\mu_0, \nu_0)$.
In our application, the process of interest---namely ULMC---is an instance of what is sometimes called a ``contractive noisy iteration'' (CNI): an algorithm that interleaves Lipschitz mappings with (Gaussian) noise convolution steps.
This notion of a Contractive Noisy Iteration is of broad interest as it captures algorithms in differential privacy (e.g., noisy optimization algorithms) and in sampling (e.g., discretizations of diffusions), and we therefore place our results in a framework which encompasses these various use cases.

A generalization of the regularization result we seek is to prove that for a CNI\@,
\begin{align}\label{eq:techniques-shift-red}\tag{$\star$}
    \cR_q(\mu_n \mmid \nu_n)
    &\lesssim \cR_q^{(w)}(\mu_0 \mmid \nu_0) + [\text{error term depending on}~w]\,,
\end{align}
where $\cR_q^{(w)}$ is the \emph{shifted R\'enyi divergence}, defined as
\begin{align*}
    \cR_q^{(w)}(\mu \mmid \nu)
	\deq \inf_{\mu'\;\mathrm{s.t.}\;W(\mu,\mu') \le w} \cR_q(\mu' \mmid \nu)\,,
\end{align*}
see \S\ref{sec:pabi} for details.
Indeed, if we take $w = W(\mu_0,\nu_0)$ in~\eqref{eq:techniques-shift-red}, then the term $\cR_q^{(w)}(\mu_0 \mmid \nu_0)$ vanishes, and we will have controlled $\cR_q(\mu_n \mmid \nu_n)$ in terms of $W(\mu_0,\nu_0)$ as desired.
However,~\eqref{eq:techniques-shift-red} is more general, as it allows for carefully tracking the shift parameter $w$ throughout.
This proof technique, called \emph{shifted divergence analysis},
was first introduced in the context of differential privacy by~\cite{pabi} for the purpose of establishing Privacy Amplification by Iteration, 
and was recently honed into a form amenable to sampling analyses in~\cite{AltTal22dp,AltTal22mix}.

\par A subtle yet essential technical issue that arises in establishing~\eqref{eq:techniques-shift-red} is: which Wasserstein metric $W$ do we use?
All previous versions of~\eqref{eq:techniques-shift-red} required the $W_{\infty}$ metric, which is problematic for our setting as the $W_\infty$ metric is infinite at initialization.
Here, our main insight is to use a non-standard Wasserstein metric, called the \emph{Orlicz--Wasserstein metric}, based on the sub-Gaussian Orlicz norm.
As we discuss in Remark~\ref{rmk:weaker-shifts}, this is exactly the right metric to use: in fact,~\eqref{eq:techniques-shift-red} cannot hold for any weaker metric (e.g., $W_p$ for any finite $p$), and the initialization bound cannot be finite for any stronger metric.
We then show that for Orlicz--Wasserstein shifts,~\eqref{eq:techniques-shift-red} indeed holds, with the caveat that the order of the shifted R\'enyi divergence on the right-hand side of~\eqref{eq:techniques-shift-red} is increased.
This increase in the order also means that additional care is required when applying~\eqref{eq:techniques-shift-red}, as the inequality cannot be iterated too many times, but we bypass this issue by showing that it suffices to only exploit the regularization from a \emph{single} step.

Finally, we note that our analysis answers the open question raised in~\cite{AltTal22mix} of how to use the shifted divergence technique in order to obtain sampling guarantees for discretized diffusions w.r.t.\ the true target distribution $\pi$, rather than w.r.t.\ the biased limit of the algorithm.

\subsubsection{From warm starts to faster high-accuracy samplers}

In light of the discussion thus far, combining our warm start result with the recent advances on MALA~\cite{chewi2021optimal, wuschche2022minimaxmala} immediately improves the dimension dependence of high-accuracy log-concave sampling to $\Otilde(d^{1/2})$. However, two further issues remain. First, thus far we have ignored the dependence on the condition number $\kappa$ for simplicity of exposition, but the combined approach of ULMC and MALA incurs suboptimal dependence on $\kappa$, namely $\kappa^{3/2}$ rather than $\kappa$.
% , which is undesirable as the condition number is also quite important for applications. %
Second, the result only holds for strongly log-concave targets. We address both of these issues simultaneously by adding a third algorithmic building block: the proximal sampler. Below, we briefly overview the proximal sampler and the final remaining technical challenges in its application.

\paragraph*{Algorithmic framework.}

The proximal sampler~\cite{titsias2018auxiliary, leeshentian2021rgo} is a Gibbs sampling method that can be viewed as the sampling analog of the proximal point method from optimization (see Appendix~\ref{app:hiacc:ps} for further background).
Each iteration requires sampling from a regularized distribution called the restricted Gaussian oracle (RGO), parametrized by $y\in\R^d$:
\begin{align*}
    \pi^{X\mid Y=y}(x)
    &\propto \exp\Bigl( - f(x) - \frac{1}{2h}\,\norm{y-x}^2\Bigr)\,.
\end{align*}
If $f$ is $\beta$-smooth, and the step size $h$ is chosen as $h \asymp \frac{1}{\beta}$, one can check that $\pi^{X\mid Y=y}$ is strongly log-concave and log-smooth with condition number $O(1)$.
Hence:
\begin{align*} 
\boxed{\begin{array}{c} \text{complexity of the} \\ \text{proximal sampler} \end{array}} = \boxed{\text{\# outer loops}} \times \boxed{\begin{array}{c} \text{complexity of sampling from} \\ O(1)\text{-conditioned distributions} \\ \text{to \textcolor{magenta}{\emph{high} accuracy}} \end{array}}
\end{align*}
The requirement of sampling from the RGO to \textcolor{magenta}{high accuracy} arises to avoid accumulation of the errors from inexact implementation of the RGO\@.

So far, we have not made use of any assumptions on $\pi$ beyond smoothness of $f$.
Additional assumptions on $\pi$, such as log-concavity, can then used to control the number of outer loops.
This program was carried out in~\cite{chenetal2022proximalsampler}, which carefully studied the outer loop complexity of the proximal sampler under a variety of assumptions on the target $\pi$ which, when combined with the implementation of the RGO via existing high-accuracy samplers, yielded state-of-the-art complexity bounds for sampling under those assumptions.
Our faster high-accuracy log-concave sampler provides a better implementation of the RGO\@, and hence we improve upon these prior results by a factor of roughly $d^{1/2}$ in each setting.
Moreover, in the strongly log-concave setting, the number of outer iterations of the proximal sampler is shown to be $\Otilde(\kappa)$~\cite{leeshentian2021rgo}, so using ULMC + MALA to implement the RGO boosts the condition number dependence of the overall sampler to near-linear.
This resolves the two issues described above, but in doing so we must also develop an inexact error analysis for the proximal sampler.

\paragraph{Inexact error analysis.}
In order to apply the proximal reduction framework, we must understand how the error from inexact implementation of the RGO propagates into the final sampling error.
This was carried out in~\cite{leeshentian2021rgo} for the TV distance via a simple coupling argument, which amounts to a union bound over failure events at each iteration.
Similarly, it is straightforward to carry out the inexact error analysis in the Wasserstein metric due to the availability of the triangle inequality. However, to establish our guarantees in \S\ref{sec:hiacc}, which hold also in the KL and $\chi^2$ metrics, we must perform an error analysis in $\chi^2$ (or equivalently, in R\'enyi).
This is also complicated by the fact while the outer loop of the proximal sampler converges exponentially fast in the strongly log-concave setting, which facilitates summing up the geometrically decaying errors from each iteration, the convergence in the weakly log-concave setting does not have an exponential rate and moreover uses a modified Lyapunov functional, changing the nature of the error analysis.
We remark that prior works such as~\cite{chenetal2022proximalsampler} did not encounter such issues, since their rejection sampling implementation of the RGO is \emph{exact}.
Therefore, we believe that our inexact error analysis will also be useful for any future applications of the proximal sampler.

We also remark that our application of the proximal sampler, and the ensuing need for careful inexact error analysis, resembles the use of the (accelerated) proximal point method in optimization, e.g.,~\cite{froetal15unreg, linmaihar15catalyst}.

\subsection{Related work}\label{ssec:intro:prior}

\paragraph*{Low-accuracy sampling and R\'enyi guarantees.}
R\'enyi guarantees for sampling are relatively recent. Indeed,~\cite{VempalaW19} proved fast R\'enyi mixing for LD and LMC to their respective stationary distributions, and this was translated into R\'enyi sampling guarantees for LMC in~\cite{GaneshT20, chewi2021optimal, chewi2021analysis, erdhoszha22chisq}, for the proximal sampler in~\cite{chenetal2022proximalsampler}, and for ULMC in~\cite{GaneshT20, Matt23girsanov}.
These lines of work have led to $\Otilde(d)$ dimension dependence for LMC and the proximal sampler, but for ULMC the rates are much worse, namely $\Otilde(d^{5/2})$ dependence~\cite{GaneshT20} and only very recently $\Otilde(d^2)$ dependence~\cite{Matt23girsanov}. In Theorem~\ref{thm:loacc}, we obtain the first $\Otilde(d^{1/2})$ rate in R\'enyi. %

In contrast, there are many more works which break the $\Otilde(d)$ barrier in the Wasserstein metric: the randomized midpoint discretization of Langevin~\cite{hebalasubramanianerdogdu2020randomizedmidpoint}, unadjusted Hamiltonian Monte Carlo (HMC)~\cite{chenvempala2019hmc}, ULMC~\cite{cheng2018underdamped, dalalyanrioudurand2020underdamped, monmarche2021high}, and more sophisticated discretizations of ULMC and HMC~\cite{shen2019randomized, foslyoobe21shiftedode, boumar22hamiltonian}.
Among these algorithms, at present we only understand how to perform R\'enyi discretization analysis for ULMC\@, but for ULMC it is the convergence of the corresponding \emph{continuous-time} diffusion which remains elusive, as we review next.

\paragraph*{Underdamped Langevin, hypoellipticity, and hypocoercivity.}
The underdamped Langevin diffusion has a rich history, dating back to Kolmogorov~\cite{Kol34}. The PDE governing the evolution of its marginal density is referred to as the kinetic Fokker--Planck equation.
Unlike the Langevin diffusion, which is driven by a full-dimensional Brownian motion and for which regularity and convergence fall within the purview of classical elliptic and parabolic PDE theory, the underdamped Langevin diffusion is the canonical example of a degenerate diffusion for which these and related questions remain active areas of research within PDE\@.
See \S\ref{ssec:ulmc-background} for background.
\par The question of regularity for these equations was largely
solved by H\"ormander~\cite{Hor67} in arguably one of the most influential breakthroughs in PDE theory of the last century through the introduction of the theory of \emph{hypoellipticity}.
In turn, it inspired Villani to coin the study of the convergence of such equations \emph{hypocoercivity} in his seminal monograph~\cite{villani2009hypocoercivity}.

While convergence of this diffusion has been studied for nearly a century, early convergence results were qualitative in nature. It took intensive developments in the PDE community to get to a point where quantitative rates could be extracted, beginning in the 1970s~\cite{tropper1977ergodic}.
We do not attempt to comprehensively survey the extensive literature here.
We refer to the monograph~\cite{villani2009hypocoercivity} for history; see also, e.g., the papers~\cite{dms, baudoin2017bakryemeryvillani, rousselstoltz2018langevin} for more modern references.
We also mention the recent space-time Poincar\'e approach of~\cite{albetal19kineticfp, caoluwang2020underdamped}, which is also directly inspired by H\"ormander's hypoelliptic theory.
As we discuss in \S\ref{ssec:intro:tech}, however, all of these approaches fall short of establishing the key property of \hyper{}.

\paragraph*{MALA.}
MALA has been intensely studied over the past three decades since its introduction in ~\cite{besagetal1995bayesiancomp}, in large part due to its strong practical performance---in fact, it and its variants comprise the default implementations of sampling routines in many modern software packages~\cite{gelleeguo15stan, abadi2016tensorflow}.
Many classical works studied the geometric ergodicity and asymptotic properties of MALA\@. With regards to the dimension dependence, particularly influential was the optimal scaling result of~\cite{robertsrosenthal1998optimalscaling}, which showed that taking step size $h\propto d^{-1/3}$ leads to a non-trivial diffusion limit for MALA as $d\to\infty$, at least for product measures $\pi$ satisfying strong regularity assumptions and when initialized at stationarity.
Modern analysis techniques have enabled an understanding of the \emph{non-asymptotic} complexity of MALA~\cite{dwivedi2018log, chenetal2020hmc, leeshentian2020gradientconcentration, chewi2021optimal, lee2021lower, wuschche2022minimaxmala}, see Table~\ref{tab:mala-progress} for a summary of the progress in this direction.
Our work closes this line of work by showing that the warm start rate of~\cite{wuschche2022minimaxmala}, which is tight due to their matching lower bound, is achievable. Moreover, our work provides theoretical justification for the improved empirical performance of MALA after using a low-accuracy algorithm for warm starts, as observed in~\cite{wuschche2022minimaxmala}.

\paragraph*{Proximal sampler.}
The proximal sampler is an algorithmic framework introduced in~\cite{titsias2018auxiliary, leeshentian2021rgo}.
% , see Appendix~\ref{app:hiacc:ps} for details. 
In~\cite{leeshentian2021rgo}, it was used as a mechanism for boosting the condition number dependence of any high-accuracy log-concave sampler to near-linear, which was then used to design samplers for composite and finite-sum potentials.
Then, in~\cite{chenetal2022proximalsampler}, it was shown that the proximal sampler reduces the problem of sampling from distributions satisfying weak log-concavity or functional inequalities to the problem of high-accuracy log-concave sampling. In this work, we exploit both these properties of the proximal sampler, and we contribute to its inexact error analysis (see \S\ref{ssec:intro:tech}).

We also mention that in recent work, the proximal sampler has been connected to stochastic localization, leading to recent progress on the KLS conjecture~\cite{klaput21spectral, cheneld22localization, klaleh22kls}, as well as to diffusion models~\cite{chenetal2023diffusionmodels}. There are also applications to sampling from semi-smooth or non-smooth potentials~\cite{liang2021proximal, liachen22proximal, liachen22proximalnoncvx}, and to differential privacy~\cite{gopleeliu2022expomech, gopietal23loglaplace, gopietal23privategeneralnorm}.

\paragraph*{Zigzag sampler.}
The zigzag sampler is an alternative high-accuracy sampler that was recently proposed in~\cite{biefearob19zigzag}. Instead of using a Metropolis--Hastings filter, the zigzag sampler is a piecewise deterministic Markov process which can be implemented without discretization bias. It was recently shown in~\cite{luwang2022zigzag} that similarly to MALA, the zigzag sampler has a dimension dependence of $\Otilde(d^{1/2})$ from a warm start.
Morever, in~\cite[Corollary 1.4]{luwang2022zigzag}, Lu and Wang show that by using LMC with a large step size to warm start the algorithm, one obtains a high-accuracy log-concave sampler with dimension dependence $\Otilde(d^{4/5})$.
Indeed, the same strategy can be used with the warm start results of~\cite{chewi2021optimal, wuschche2022minimaxmala} to obtain complexities strictly better than $\Otilde(d)$; however, it is clear that such an approach can never reach the desired complexity of $\Otilde(d^{1/2})$---and in fact there is a fundamental barrier even at $\Otilde(d^{3/4})$ because it is bottlenecked by the discretization bias of LMC.
The goal of this paper is achieving $\Otilde(d^{1/2})$ complexity as this is this a natural barrier for high-accuracy samplers given a warm start, and LMC cannot work for this goal.\footnote{With regards to dimension dependence, running LMC with step size $h$ for $1/h$ steps yields a distribution $\mu$ with $\log \chi^2(\mu \mmid \pi) \le \Otilde(dh)$~\cite{chewi2021analysis}. By optimizing the step size $h$ and combining this with the best known complexity $\Otilde(d^{1/2} \log^{3/2} \chi^2(\mu \mmid \pi))$ of the zigzag sampler, one obtains the final complexity $\Otilde(d^{4/5})$~\cite{luwang2022zigzag}. We point out that even if the complexity of the zigzag sampler were improvable to $\Otilde(d^{1/2} \log \chi^2(\mu \mmid \pi))$, the total complexity would still be at least  $1/h + d^{1/2}\,(dh) \geq  \widetilde\Omega(d^{3/4})$. Thus $d^{3/4}$ is a natural barrier for any warm starting approach using LMC.} In direct analogy to our use of MALA, our new complexity result for ULMC (Theorem~\ref{thm:loacc}) can also be used to warm start the zigzag sampler, leading to the same final complexity bound of $\Otilde(\kappa d^{1/2} \log^{O(1)}(1/\eps))$. This answers the open questions in~\cite{luwang2022zigzag} regarding warm starting the zigzag sampler.

\paragraph*{Differential privacy and sampling.} Sampling algorithms have been widely used in differential privacy ever since the invention of the exponential mechanism~\cite{McSherryT07}; for an exposition of the surrounding history and applications, see the textbook~\cite{dwork2014algorithmic}. Sampling-inspired analyses have also been recently used to prove privacy properties of optimization algorithms~\cite{ChourasiaYS21,RyffelBP22,YeS22}. Most related to this paper are connections in the other direction: the use of techniques from differential privacy in order to analyze sampling. 
There are two lines of work in this direction. One involves the technique of adaptive composition for R\'enyi divergences and its use for establishing R\'enyi bias bounds for LMC and ULMC~\cite{GaneshT20, erdhoszha22chisq, Matt23girsanov}.
The other involves the technique of Privacy Amplification by Iteration (PABI), which was originally used to bound the privacy loss of differentially private optimization algorithms~\cite{pabi, BalleBGG19, Asoodeh20, FeldmanKoTa20, Sordello21,AltTal22dp}, and its recent use for analyzing the mixing time of LMC to its biased stationary distribution~\cite{AltTal22mix}. In this paper, we build upon this technique in several key ways: we show how to improve mixing results for the biased distribution to mixing results for the target distribution, we show how to use these ideas for ULMC rather than LMC, and most importantly we overcome the key issue of unboundedness by replacing $W_{\infty}$ shifts by Orlicz--Wasserstein shifts, see \S\ref{ssec:intro:tech} for an overview.

\subsection{Simultaneous work}\label{ssec:intro:simultaneous}

While preparing a draft of this paper for submission, it came to our attention that another group was simultaneously working towards the same problem using completely different techniques and algorithms~\cite{Fan23dimension}. We are grateful to them for coordinating simultaneous arXiv submissions. We look forward to reading their paper after it is posted, and we will add a detailed comparison about the differences in a future revision.

\subsection{Organization}\label{ssec:intro:org}

We recall preliminaries in \S\ref{sec:prelim}, especially regarding R\'enyi divergences. We isolate in \S\ref{sec:pabi} our key new technique involving Orlicz--Wasserstein shifted R\'enyi divergences. We use this technique to obtain faster algorithmic warm starts in \S\ref{sec:ulmc}, and then use these warm starts to develop faster high-accuracy samplers in \S\ref{sec:hiacc}. We conclude in \S\ref{sec:discussion} by discussing several future research directions that are motivated by our results. For brevity, we defer various proofs and technical details to the appendices.

%% file: sections/prelim.tex
\section{Preliminaries}\label{sec:prelim}

\subsection{Notation}

Throughout, $\pi \propto \exp(-f)$ denotes the target density and $f : \R^d \to \R$ denotes the potential. We assume $f$ is twice continuously differentiable for simplicity. We reserve the symbol $N$ for the number of iterations that the Markov Chain Monte Carlo algorithm is run, $h$ for the discretization step size, and $T = Nh$ for the total elapsed time. We write $\cP(\R^d)$ to denote the space of probability distributions over $\R^d$, and we write $\cP_2(\R^d)$ to denote the subset of $\cP(\R^d)$ with finite second moment. All logarithms are natural.

\par For simplicity of exposition, we assume throughout that we have access to an algorithm for generating independent standard Gaussian random variables. We use the standard notation $\Otilde(g) = g \log^{O(1)}(g)$ to suppress low-order terms. Note that since our final results depend polynomially on the dimension $d$ and condition number $\kappa$, the $\Otilde$ hides polylogarithmic factors in these terms---on the other hand, since $\eps$ occurs only polylogarithmically in our high-accuracy results, we do not hide the polylogarithmic factors in $\eps$.

We say that $f$ is $\alpha$-strongly convex if $\nabla^2 f \succeq \alpha I_d$, and that $f$ is $\beta$-smooth if $\norm{\nabla^2 f}_{\rm op} \le \beta$. If $f$ is convex, then $f$ is $\beta$-smooth if and only if $\nabla^2 f \preceq \beta I_d$.
We always denote by $\kappa$ the condition number $\kappa \deq \beta/\alpha$.
If $\pi\propto\exp(-f)$ where $f$ is $\alpha$-strongly convex (resp.\ $\beta$-smooth), we say that $\pi$ is $\alpha$-strongly log-concave (resp.\ $\beta$-log-smooth). All other notation is introduced in the main text.

\subsection{R\'enyi divergences}\label{ssec:prelim-renyi}

\begin{defin}[R\'enyi divergence]
	The R\'enyi divergence of order $q \in (1, \infty)$ between probability measures $\mu$ and $\nu$ is defined as
	\[
		\cR_q(\mu \mmid \nu) = \frac{1}{q-1} \log \int \Bigl( \frac{\D\mu}{\D\nu} \Bigr)^q \,\D\nu
	\]
	if $\mu \ll \nu$, and otherwise is $\infty$. The R\'enyi divergences of order $q \in \{1,\infty\}$ are defined by continuity. 
\end{defin}

\begin{remark}[Special cases of R\'enyi divergence]\label{rem:renyi-cases}
	The R\'enyi divergence of order $q=1$ coincides with the KL divergence, i.e.,
	\[
		\cR_1(\mu \mmid \nu) = \KL(\mu \mmid \nu)\,.
	\]
	The R\'enyi divergence of order $q=2$ is related to the $\chi^2$ divergence via the formula
	\[
            \cR_2(\mu \mmid \nu)
		= \log\bigl(1+\chi^2(\mu \mmid \nu)\bigr)\,.
	\]
	The R\'enyi divergence of order $q =\infty$ is given by
	\[
		\cR_\infty(\mu \mmid \nu) = \log \esssup \frac{\D\mu}{\D\nu}\,.
	\] 
\end{remark}

Our analysis makes use of the following elementary properties of the R\'enyi divergence. Further details about these properties and their proofs can be found in, e.g., the R\'enyi divergence survey~\cite{van2014renyi} as Theorem 1, Theorem 3, Equation 10, and Remark 1, respectively.

\begin{lemma}[Post-processing inequality for R\'enyi divergences]\label{lem:renyi-process}
	For any R\'enyi order $q \geq 1$, any Markov transition kernel $P$, and any probability distributions $\mu,\nu$,
	\[
		\cR_q\left( \mu P \mmid \nu P \right) \leq \cR_q\left( \mu \mmid \nu \right)\,.
	\]
\end{lemma}

\begin{lemma}[Monotonicity of R\'enyi divergences]\label{lem:renyi-monotonicity}
	For any R\'enyi orders $q' \geq q \geq 1$, and any probability distributions $\mu,\nu$, 
	\[
	\cR_q(\mu \mmid \nu) \leq \cR_{q'}(\mu \mmid \nu)\,.
	\]
\end{lemma}

\begin{lemma}[R\'enyi divergence between isotropic Gaussians]\label{lem:renyi-gaussians}
    For any R\'enyi order $q \geq 1$, any variance $\sig^2 > 0$, and any means $x,y \in \R^d$, 
    \[
    	\cR_q\bigl(\cN(x,\sig^2 I_d) \bigm\Vert \cN(y,\sig^2 I_d)\bigr) = \frac{q\,\|x-y\|^2}{2\sig^2}\,.
    \]
\end{lemma}

\begin{lemma}[Relation to $f$-divergences]\label{lem:renyi-convex}
	For any R\'enyi order $q \in (1,\infty)$, the corresponding function $\exp((q-1) \; \cR_q(\cdot \mmid \cdot))$ is an $f$-divergence, and thus in particular is jointly convex in its arguments.
\end{lemma}

We end this section with one last property of R\'enyi divergences: the weak triangle inequality. The name of this property arises from the fact that although R\'enyi divergences do not satisfy the triangle inequality, they do satisfy a modified version of it in which the R\'enyi order is increased and the bound is weakened by a multiplicative factor. Since this property does not appear in the aforementioned survey~\cite{van2014renyi} on R\'enyi divergences, we provide a brief proof for completeness. It can also be found in, e.g.,~\cite[Proposition 11]{mironov2017renyi}.

\begin{lemma}[Weak triangle inequality for R\'enyi divergence]\label{lem:renyi-triangle}
	For any R\'enyi order $q > 1$, any  $\lambda \in (0, 1)$, and any probability distributions $\mu,\nu,\pi$,
	\begin{align*}
		\cR_q(\mu \mmid \pi)
		&\le \frac{q-\lambda}{q-1} \, \cR_{q/\lambda}(\mu \mmid \nu) + \cR_{(q-\lambda)/(1-\lambda)}(\nu \mmid \pi)\,.
	\end{align*}
\end{lemma}
\begin{proof}
	Expand $\cR_q(\mu \mmid \nu) = \frac{1}{q-1} \log \int fg$ where $f = \mu^q/\nu^{q-\lambda}$ and $g = \nu^{q - \lambda}/\pi^{q-1}$, and then apply H\"older's inequality $\int f g \leq (\int f^a)^{1/a}\, (\int g^b)^{1/b}$ using H\"older exponents $a = 1/\lambda$ and $b = 1/(1 - \lambda)$. 
\end{proof}

%% file: sections/pabi.tex
\section{Improved shifted divergence analysis}\label{sec:pabi}

In this section we isolate from our analysis a key new technique of independent interest. As overviewed in \S\ref{ssec:intro:tech}, this technique is a strengthening of the ``shifted divergence'' analysis, a.k.a., ``privacy amplification by iteration'' (PABI), in which we crucially 
improve the $\infty$-Wasserstein shift to an Orlicz--Wasserstein shift. This enables obtaining R\'enyi divergence bounds on the mixing of any Markov chain which interleaves Lipschitz mapping steps (e.g., gradient descent steps) and noise convolution steps (e.g., adding a Gaussian). This notion captures a variety of algorithms from the differential privacy and sampling communities, often called ``contractive noisy iterations''. 

The main result of this section is formally stated as follows. This result makes use of the Wasserstein metric $W_{\psi_2}$ that evaluates a coupling's quality via the sub-Gaussian Orlicz norm; see \S\ref{ssec:pabi-orlicz} for background on this notion. 

\begin{theorem}[Shifted R\'enyi divergence analysis]\label{thm:pabi-orlicz}
    Consider two Markov chains $\{\mu_n\}_{n \geq 0}$ and $\{\mu_n'\}_{n \geq 0}$ with possibly different initialization, but with the same update transitions 
    \begin{align*}
        \mu_{n+1}  &= (\mu_n P_n) \ast \cN(0,\sig^2 I_d) \\
        \mu_{n+1}' &= (\mu_n' P_n) \ast \cN(0,\sig^2 I_d)
    \end{align*}
    where $P_n$ is a Markov transition kernel that is $c$-Lipschitz in the $W_{\psi_2}$ metric.
    Then for any R\'enyi order $q \geq 1$,
	\begin{align}
		\cR_q(\mu_N \mmid \mu_N') \leq 
		c^{2N}\,
				 \frac{q\,W_{\psi_2}^2(\mu_0,\mu_0')}{2\sigma^2} \,,
		\label{eq:pabi-bound-new}
	\end{align}
	so long as $N \geq \log_{1/c} \Bigl( \frac{\sqrt{q\,(q-1)}\,W_{\psi_2}(\mu_0,\mu_0')}{\sigma\sqrt 2} \Bigr)$.
\end{theorem}

\par We remark that unlike previous versions of the shifted divergence technique, Theorem~\ref{thm:pabi-orlicz} requires a restriction on the number of iterations $N$.\footnote{This restriction comes from the fact that with this new Orlicz--Wasserstein shifted R\'enyi divergence, the new shift-reduction lemma (Lemma~\ref{lem:shift-reduction}) does not apply to arbitrarily large shifts.}
But this restriction is mild due to the logarithmic dependence. In fact, it is equivalent to requiring the upper bound in~\eqref{eq:pabi-bound-new} to be at most $1/(q-1)$.

\par The rest of this section is devoted to proving Theorem~\ref{thm:pabi-orlicz}. In \S\ref{ssec:pabi-orlicz} we define a new Lyapunov function, and in \S\ref{ssec:pabi-proof} we use it to prove Theorem~\ref{thm:pabi-orlicz}.

\subsection{Shifted R\'enyi divergence using Orlicz--Wasserstein shifts}\label{ssec:pabi-orlicz}

\par Key to our proof of Theorem~\ref{thm:pabi-orlicz} is a new Lyapunov function for tracking how indistinguishable the Markov chains become as they evolve. This new Lyapunov function is a shifted R\'enyi divergence, but unlike the standard shifted divergence technique, here we measure the shift using an ``Orlicz--Wasserstein metric'' rather than $W_{\infty}$.

\par We begin by recalling the definition of a sub-Gaussian Orlicz norm. For shorthand, we drop the adjective ``sub-Gaussian'' as this is the only Orlicz norm considered in this paper. For further background on Orlicz norms, we refer the reader to, e.g., the textbooks~\cite{Rao91book,vershynin2018high}, and we mention that the standard significance of this particular (sub-Gaussian) Orlicz norm is that a random variable is sub-Gaussian if and only if this norm is finite~\cite[Example 2.7.13]{vershynin2018high}.

\begin{defin}[Orlicz norm]\label{def:orlicz-norm}
	The Orlicz norm of a random variable $X$ is
	\begin{align*}
			\norm X_{\psi_2}
			&\deq \inf\biggl\{\lambda > 0 \; : \; \E \psi_2\Bigl( \frac{\norm X}{\lambda}\Bigr) \le 1\biggr\}\,,
		\end{align*}
	where the function $\psi_2 : \R \to \R$ is defined as $\psi_2(x) \defeq \exp(x^2) - 1$.
\end{defin}

Our proof of Theorem~\ref{thm:pabi-orlicz} uses the Orlicz norm for defining an optimal transport metric between probability distributions. In what follows, we write $\cC(\mu,\nu)$ to denote the set of couplings between $\mu$ and $\nu$; that is, the set of all jointly defined pairs $(X,Y)$ of random variables with first marginal $\law(X) = \mu$ and second marginal $\law(Y) = \nu$.

\begin{defin}[Orlicz--Wasserstein metric]\label{def:w-orlicz}
	The Orlicz--Wasserstein metric between distributions $\mu,\nu$ is 
	\[
	W_{\psi_2}(\mu,\nu) \deq \inf_{(X,Y) \in \cC(\mu,\nu)} \|X - Y\|_{\psi_2}\,.
	\]
\end{defin}

Note that since the Orlicz norm satisfies the triangle inequality, the standard gluing lemma from classical optimal transport theory shows that $W_{\psi_2}$ is indeed a metric; see, e.g.,~\cite[Chapter 6]{villani2009optimal}.
The Orlicz{--}Wasserstein metric has also been considered in prior works~\cite{sturm11orlicz, kell17interpcurv, guhahongu23gmm}, but to our knowledge this paper constitutes the first use of this metric for sampling analysis.

\begin{remark}[Comparison to $W_p$]\label{rem:wass-comparison}
	Let $W_p$ denote the standard $p$-Wasserstein metric on $\R^d$. Then 
	\[
	\frac{1}{\sqrt p}\, W_p \lesssim W_{\psi_2} \leq \frac{1}{\sqrt{\log2}}\, W_{\infty}\,.
	\]
 Intuitively, the first inequality is because a finite Orlicz norm implies sub-Gaussianity with a related bound on the concentration parameter, which implies related moment bounds; and the second inequality is because an $\esssup$ bound implies compact support, which implies sub-Gaussianity by Hoeffding's lemma. Proofs of these bounds are provided in Appendix~\ref{app:wass-comparison}. We also note that the reverse inequalities do not hold, even if one weakens them by an arbitrarily large amount. For example, for the first inequality,
    take $\mu = \delta_0$ and $\nu$ the Laplace distribution with density $\nu(x) = \frac{1}{2}\exp(-\abs x)$ on $\R$; then, $W_p(\mu,\nu)$ is finite for any $1 \le p < \infty$, but $W_{\psi_2}(\mu,\nu) = \infty$. And for the second inequality, take $\mu = \delta_0$ and $\nu = \mc N(0, 1)$; then  $W_{\psi_2}(\mu,\nu)$ is finite but $W_\infty(\mu,\nu) =\infty$.
\end{remark}

\begin{defin}[Orlicz--Wasserstein shifted R\'enyi divergence]
	For any R\'enyi order $q \geq 1$ and shift $w \geq 0$, the $W_{\psi_2}$-shifted R\'enyi divergence between probability distributions $\mu$ and $\nu$ is defined as
	\[
	\cR_q^{(w)}(\mu \mmid \nu)
	\deq \inf_{\mu'\;\mathrm{s.t.}\;W_{\psi_2}(\mu,\mu') \le w} \cR_q(\mu' \mmid \nu)\,.
	\]
\end{defin}

\subsection{Proof of Theorem~\ref{thm:pabi-orlicz}}\label{ssec:pabi-proof}

Here we describe how the standard shifted divergence analysis is modified when using shifts in $W_{\psi_2}$ rather than $W_{\infty}$, and how this modified argument leads to a proof of Theorem~\ref{thm:pabi-orlicz}.

At a high level, the shifted divergence technique---in both its original form and the new form here---is built upon two key lemmas. These two lemmas track how the shifted R\'enyi divergence evolves when both distributions are either (1) pushed forward through a Lipschitz map; or (2) convolved with Gaussian noise. These two lemmas are called the ``contraction-reduction lemma\footnote{Although we use this name to be consistent with the previous literature on the shifted divergence technique, we note that this map need not be a contraction, i.e., the Lipschitz constant can be greater than $1$.}'' and the ``shift-reduction lemma.''

\par The contraction-reduction lemma is the simpler of these two lemmas, and extends unchanged---in terms of both statement and proof---when the standard $W_{\infty}$ shift is replaced by our proposed $W_{\psi_2}$ shift. For completeness, we provide a brief proof.

\begin{lemma}[New contraction-reduction lemma, for Orlicz--Wasserstein shifted R\'enyi]\label{lem:contraction-reduction}
	For any R\'enyi order $q \geq 1$, any shift $w \geq 0$, any Markov transition kernel $P$ that is $W_{\psi_2}$-Lipschitz with parameter $c$, and any distributions $\mu,\nu \in \cP(\R^d)$, 
	\[
		\cR_q^{(w)}\left( \mu P \mmid \nu P \right)
		\leq
		\cR_q^{(w/c)} \left( \mu \mmid \nu \right)\,.
	\]
\end{lemma}
\begin{proof}
	Let $\mu'$ be the surrogate for $\mu$ in $\cR_q^{(w/c)} \left( \mu \mmid \nu \right)$. Then by definition, $W_{\psi_2}(\mu,\mu') \leq w/c$ and $\cR_q \left( \mu' \mmid \nu \right) = \cR_q^{(w/c)} \left( \mu \mmid \nu \right)$. Thus
	\[
		\cR_q^{(w)}( \mu P \mmid \nu P)
		\leq
		\cR_q( \mu' P \mmid \nu P)
		\leq 
		\cR_q ( \mu' \mmid \nu)
  =
		\cR_q^{(w/c)}( \mu \mmid \nu)\,,
	\]
	where the first step is because $W_{\psi_2}(\mu P,\mu'P) \leq c\,W_{\psi_2} (\mu, \mu') \leq c\,(w/c) = w$ by Lipschitzness of $P$; the second step is by the data-processing inequality for R\'enyi divergences (Lemma~\ref{lem:renyi-process}); and the third step is by construction of $\mu'$.
\end{proof}

The shift-reduction lemma, however, requires substantial modification.

\begin{lemma}[New shift-reduction lemma, for Orlicz--Wasserstein shifted R\'enyi]\label{lem:shift-reduction}
	For any R\'enyi order $q \geq 1$, any noise variance $\sig^2 > 0$, any initial shift $w \geq 0$, any shift increase 
			$\delta \leq \sig/\sqrt{(2q-1)\,(q-1)}$,
	and any distributions $\mu,\nu \in \cP(\R^d)$, 
		\begin{align*}
		\cR_q^{(w)}\bigl(\mu \ast \cN(0,\sig^2 I_d) \bigm\Vert \nu \ast \cN(0, \sig^2 I_d)\bigl)
		&\leq \cR_{2q-1}^{(w+\delta)}(\mu \mmid \nu) 
		+ 
		\frac{(2q-1)\,\delta^2}{2\sig^2}\log 2
		\,.
	\end{align*}
\end{lemma}

\begin{proof}
	\underline{Case 1: initial shift $w =0$.}
	For shorthand, let $\gamma$ denote $\cN(0,\sig^2 I_d)$.
	We bound the R\'enyi divergence between $\law(X+Z)$ and $\law(Y+Z)$, where $X\sim \mu$, $Y\sim \nu$, and $Z \sim \gamma$. Let $\mu'$ be the surrogate for $\cR_{2q-1}^{(\delta)}(\mu \mmid \nu)$, so that $\cR_{2q-1}^{(\delta)}(\mu \mmid \nu) = \cR_{2q-1}(\mu' \mmid \nu)$ and $W_{\psi_2}(\mu,\mu') \le \delta$. Let $X'\sim\mu'$ be optimally coupled with $X \sim \mu$ with respect to the Orlicz--Wasserstein metric $W_{\psi_2}(\mu,\mu')$ so that
	\begin{align}
		\iint \psi_2\Bigl( \frac{\|x-x'\|}{\delta} \Bigr) \, p_{X,X'}(\D x,\D x') \leq 1\,,
		\label{eq:shift-reduction-coupling}
	\end{align}
	where here and henceforth we write $p_{\eta}$ as shorthand for the law of a random variable $\eta$.
	\par Note that $X+Z$ and $Y+Z$ are the result of the same function applied to the tuples $(X', X-X'+Z)$ and $(Y, Z)$ respectively. Thus, by the data-processing inequality for R\'enyi divergences (Lemma~\ref{lem:renyi-process}),
	\begin{align*}
		\cR_q(\mu \ast \gamma \mmid \nu \ast \gamma)
		=
		\cR_q\bigl( \law(X+Z) \bigm\Vert \law(Y+Z) \bigr)
		&\le \cR_q\bigl(\law(X', X-X'+Z) \bigm\Vert \law(Y, Z)\bigr)\,.
	\end{align*}
		By expanding the definition of R\'enyi divergence and applying H\"older's inequality,
		we bound this by
	\begin{align*}
		\cdots 
		&= \frac{1}{q-1} \log \iint \Bigl( \frac{ p_{X',X-X'+Z}(x',z) }{ p_{Y,Z}(x',z) } \Bigr)^{q-1}\, p_{X',X-X'+Z}( \D x', \D z) \\
		&= \frac{1}{q-1} \log \iint \Bigl( \frac{p_{X'}(x')}{p_Y(x')}\, \frac{p_{X-X'+Z\mid X'=x'}(z)}{p_{Z\mid Y=x'}(z)}\Bigr)^{q-1} \, 
		p_{X-X'+Z \mid X'=x'}(\D z) 
		\, p_{X'}(\D x') \\
		& = \frac{1}{q-1} \log \iint \Bigl( \frac{\mu'(x')}{\nu(x')}\, \frac{p_{X-X'+Z\mid X'=x'}(z)}{\gamma(z)}\Bigr)^{q-1} \, \, 
		p_{X-X'+Z \mid X'=x'}(\D z) 
		\, \mu'(\D x')  \\
		&\le \underbrace{\frac{1}{2\,(q-1)} \log \int \Bigl( \frac{\mu'(x')}{\nu(x')} \Bigr)^{
				2\,(q-1)
		} \, \mu'(\D x')}_{\circled{1}} \\
		&\qquad{} + {\underbrace{\frac{1}{2\,(q-1)} \log \iint \Bigl( \frac{p_{X-X'+Z\mid X'=x'}(z)}{\gamma(z)}\Bigr)^{
				2\,(q-1)
			} \, 
			p_{X-X'+Z \mid X'=x'}(\D z) 
			\, \mu'(\D x')}_{\circled{2}}}\,.
	\end{align*}
	\par By definition of R\'enyi divergence and then the construction of $\mu'$, the first term $\circled{1}$ simplifies to
	\[
	\circled{1} = \cR_{2q-1}(\mu' \mmid \nu)  = \cR_{2q-1}^{(\delta)}(\mu \mmid \nu)\,.
	\]
		The second term $\circled{2}$ can be bounded as
	\begin{align*}
		\circled{2} &= \frac{1}{2\,(q-1)}\log \int \exp\bigl( 2\,(q-1) \, \cR_{2q-1}(p_{X-X'+Z\mid X'=x'} \mmid \gamma)\bigr) \, \mu'(\D x')
		\\ &\leq \frac{1}{2\,(q-1)}\log \iint \exp\bigl(2\,(q-1) \,\cR_{2q-1}(p_{x-x'+Z} \mmid \gamma)\bigr) \, p_{X,X'}(\D x,\D x') 
		\\ &= \frac{1}{2\,(q-1)} \log \iint \exp\bigl( 2\,(q-1) \,\cR_{2q-1}( \cN(x-x',\sig^2 I_d) \mmid \cN(0,\sig^2 I_d) )\bigr) \, p_{X,X'}(\D x,\D x') 
		\\ &= \frac{1}{2\,(q-1)}\log \iint \exp \Bigl( \frac{(q-1)\,(2q-1)\,\norm{x-x'}^2}{\sig^2} \Bigr) \, p_{X,X'}(\D x,\D x')
		\\ &\leq \frac{\delta^2\,(2q-1)}{2\sig^2}\log \iint \exp \Bigl( \,\frac{\norm{x-x'}^2}{\delta^2} \Bigr) \, p_{X,X'}(\D x,\D x')
		\\ &\leq \frac{\delta^2\,(2q-1)}{2\sig^2} \log 2\,.
	\end{align*}
	Above, the first step is by definition of the R\'enyi divergence; the second step is by noting that $p_{X-X'+Z\mid X'=x'}=\int p_{x-x'+Z} \, p_{X\mid X'}(\D x \mid x')$ and using  convexity of $f$-divergences (Lemma~\ref{lem:renyi-convex}); the fourth step is by the closed-form expression for the R\'enyi divergence between Gaussians (Lemma~\ref{lem:renyi-gaussians});
	the fifth step is by the assumption that $\rho \deq \delta^2\,(q-1)\,(2q-1)/\sig^2\le 1$, which enables us to use Jensen's inequality to bound $\E[R^{\rho}] \leq \E[R]^{\rho}$ where $R \defeq \exp(\|X-X'\|^2/\delta^2)$;
	and the final step is by the property~\eqref{eq:shift-reduction-coupling} of the coupling $(X,X')$. %
     Combining these bounds on $\circled{1}$ and $\circled{2}$ completes the proof for case $1$.
	\par \underline{Case 2: initial shift $w > 0$.} Let $\mu'$ denote the surrogate for $\cR_{2q-1}^{(w + \delta)}(\mu \mmid \nu)$, so that $\cR_{2q-1}^{(w+\delta)}(\mu \mmid \nu) = \cR_{2q-1}(\mu' \mmid \nu)$ and $W_{\psi_2}(\mu,\mu') \le w+\delta$. Let $X'\sim\mu'$ be optimally coupled with $X \sim \mu$ with respect to the Orlicz--Wasserstein metric $W_{\psi_2}(\mu,\mu')$ so that $\|X-X'\|_{\psi_2} = W_{\psi_2}(\mu,\mu') \leq w + \delta$. Decompose
	\[
	X' = \underbrace{\tau\, X + (1 - \tau)\, X' }_{X_1'} + \underbrace{\tau\, (X' - X)}_{X_2'}\,,
	\]
	where $\tau \deq \delta/(w + \delta)$. 
	Then
	\begin{align*}
		\cR_q^{(w)}(\mu \ast \gamma \mmid \nu \ast \gamma )
		&\leq \cR_{q}( p_{X_1'} \ast \gamma \mmid \nu \ast \gamma)
		\\ &\leq \cR_{q}^{(\delta)}( p_{X_1'} \mmid \nu ) + \frac{\delta^2\,(2q-1)}{2\sig^2} \log 2
		\\ &\leq \cR_{2q-1}(\mu' \mmid \nu) + \frac{\delta^2\,(2q-1)}{2\sig^2} \log 2
		\\ &= \cR_{2q-1}^{(w+\delta)}(\mu \mmid \nu) + \frac{\delta^2\,(2q-1)}{2\sig^2} \log 2\,.
	\end{align*}
	Above, the first step is by using $p_{X_1'}\ast \gamma$ as a surrogate for $\mu\ast \gamma$, which is allowed since
 \begin{align*}
     W_{\psi_2}(\mu \ast \gamma, \,p_{X_1'} \ast \gamma) \le W_{\psi_2}(\mu,p_{X_1'}) \leq \|X - X_1'\|_{\psi_2} = (1 - \tau)\, \|X - X'\|_{\psi_2} \leq (1 - \tau)\, (w + \delta) = w\,.
 \end{align*}
 The second step is by using the result from case 1; the third step is by using $\mu'$ as a surrogate for $p_{X_1'}$, which is allowed since $W_{\psi_2}(\mu',p_{X_1'}) \leq \|X' - X_1'\|_{\psi_2} = \|X_2'\|_{\psi_2} = \tau\, \|X - X'\|_{\psi_2} \leq \tau\, (w + \delta) = \delta$; and the final step is by construction of $\mu'$.
\end{proof}

\begin{remark}\label{rem:shift-reduction-holder}
	Lemma~\ref{lem:shift-reduction} can be generalized to 
	\begin{align}\label{eq:shift-reduction-holder}
			\cR_q^{(w)}\bigl(\mu \ast \cN(0,\sig^2 I_d) \bigm\Vert \nu \ast \cN(0, \sig^2 I_d)\bigr)
			&\leq \cR_{(q+\lambda-1)/\lambda}^{(w+\delta)}(\mu \mmid \nu) + 
			\frac{(q-\lambda)\,\delta^2}{2\,(1-\lambda) \,\sigma^2} \log2
		\end{align}
	for any $\lambda \in [0,1]$ and any $\delta \leq (1-\lambda)\,\sigma \sqrt{\tfrac{2}{(q-1)\,(q-\lambda)}}$.
	Choosing $\lambda = 1/2$ recovers Lemma~\ref{lem:shift-reduction}.
	Different choices of $\lambda$ enable trading off the increase in the R\'enyi order in the first term against the penalty in the second term. 
	The proof of this generalized bound is identical except for one change: replace the Cauchy--Schwarz inequality in the proof with H\"older's inequality $\int fg \leq (\int \abs f^{1/\lambda})^{\lambda} \,(\int \abs g^{1/(1-\lambda)})^{1-\lambda}$.
 
\par We conjecture that~\eqref{eq:shift-reduction-holder} is tight in that for any input parameters $q$, $\sig^2$, $w$, there exist distributions $\mu,\nu$ such that this bound holds with equality when it is optimized over the knobs $\delta, \lambda$. Details on this conjecture are provided in Appendix~\ref{app:shift-reduction-tight}.
\end{remark}

\begin{remark}[Failure of the shift-reduction lemma for weaker shifts]\label{rmk:weaker-shifts}
    The Orlicz--Wasserstein metric is the ``right'' metric to use for the shifted R\'enyi analysis in the sense that (1) Lemma~\ref{lem:shift-reduction} holds with this Orlicz--Wasserstein shifted R\'enyi divergence; and (2) the Orlicz--Wasserstein distance at initialization is bounded for sampling algorithms (shown in Lemma~\ref{lem:ulmc-init}). 
    \par In contrast, (2) fails for all previous versions of the shifted R\'enyi divergence analysis, since they use $W_{\infty}$ shifts. And (1) fails for other natural candidates of the Wasserstein metric for which the initialization distance is bounded. This includes the $W_p$ metric for any finite $p$, as well as the Orlicz--Wasserstein metric for any Orlicz norm that is weaker than sub-Gaussian. See Appendix~\ref{app:rmk-weaker-shifts} for details. Finally, we remark that this discussion is tailored to the fact that we are analyzing Markov chains with Gaussian noise; if for example, this were replaced by Laplacian noise, then the right notion of shift would be the Orlicz--Wasserstein metric with the sub-exponential Orlicz norm, and our techniques would extend in a straightforward way.
\end{remark}

We now use Lemmas~\ref{lem:contraction-reduction} and~\ref{lem:shift-reduction} to prove Theorem~\ref{thm:pabi-orlicz}.

\begin{proof}[Proof of Theorem~\ref{thm:pabi-orlicz}]
	Let $\delta = c^N\, W_{\psi_2}(\mu_0,\mu_0')$. By using, in order: the definition of the Markov chain update, Lemma~\ref{lem:shift-reduction} (in the form of Remark~\ref{rem:shift-reduction-holder} with $\lambda = 0$),  and then Lemma~\ref{lem:contraction-reduction}, we obtain
	\begin{align*}
		\cR_q( \mu_N \mmid \mu_N')
		&= 
		\cR_q\bigl( ( \mu_{N-1} P_{N-1} ) \ast \cN(0,\sig^2I_d)  \bigm\Vert ( \mu_{N-1}' P_{N-1} ) \ast \cN(0,\sig^2I_d) \bigr) 
		\\ &\leq
		\cR_\infty^{(\delta)}( \mu_{N-1} P_{N-1}  \mmid \mu_{N-1}' P_{N-1}) + \frac{q\delta^2}{2\sig^2} \log 2
		\\ &\leq
		\cR_\infty^{(\delta/c)}( \mu_{N-1} \mmid \mu_{N-1}' ) + \frac{q\delta^2}{2\sig^2} \log 2\,.
	\end{align*}
	Note that the use of Lemma~\ref{lem:shift-reduction} is valid since 
	$\delta \leq \sig \sqrt{\frac{2}{q\,(q-1)}}$
	by the assumption on $N$.
	\par It suffices to show that the R\'enyi term in the above display vanishes. To this end, let $Q_n$ denote the transition kernel for the $n$-th step of the Markov chain, i.e., $\rho Q_n = (\rho P_n) \ast \cN(0,\sig^2 I_d)$. Clearly $Q_n$ is $W_{\psi_2}$-Lipschitz with parameter $c$ since $W_{\psi_2}(\rho Q_n, \rho' Q_n) \leq W_{\psi_2}(\rho P_n, \rho' P_n) \leq c\, W_{\psi_2} (\rho,\rho')$ for any distributions $\rho,\rho'$. Thus we may apply Lemma~\ref{lem:shift-reduction} $N-1$ times to argue that
	\begin{align*}
		\cR_{\infty}^{(\delta/c)} (\mu_{N-1} \mmid \mu_{N-1}') 
		&=
		\cR_{\infty}^{(\delta/c)}(\mu_0 Q^{N-1} \mmid \mu_0' Q^{N-1})
		\\ & \leq 
		\cR_{\infty}^{(\delta/c^N)}( \mu_0 \mmid \mu_0' )\,
		\\ &=
		\cR_{\infty}^{(W_{\psi_2}(\mu_0,\mu_0'))}( \mu_0 \mmid \mu_0')
		\\ &\leq 
		\cR_{\infty}( \mu_0' \mmid \mu_0')
		\\ &= 0\,.
	\end{align*}
	Here, the third step is by the choice of $\delta$, and the fourth step is by definition of the shifted R\'enyi divergence. The proof is complete by combining the above displays.
\end{proof}

%% file: sections/ulmc.tex
\section{Low-accuracy sampling with \texorpdfstring{$O(\sqrt{d})$}{sqrt(d)} complexity}\label{sec:ulmc}

The main result of the section is the first R\'enyi convergence guarantee for log-concave sampling that requires a number of first-order queries that scales in the dimension $d$ only as $\sqrt{d}$. This improves over the state-of-the-art which has $d^2$ scaling. This result is formally stated as follows. 

\begin{theorem}[Low accuracy sampling with $O(\sqrt{d})$ complexity]\label{thm:loacc}
	Suppose that $\pi \propto \exp(-f)$ where $f$ is $\alpha$-strongly-convex and $\beta$-smooth, and let $0 < \varepsilon \lesssim \frac{1}{\sqrt q}$.
	There is a randomized algorithm that, given knowledge of the minimizer of $f$ and access to 
	\[
		N = \Otilde\Bigl( \frac{\kappa^{3/2} d^{1/2} q^{1/2}}{\eps} \Bigr)
	\]
	gradient queries for $f$, outputs a random point in $\R^d$ with $\mu$ satisfying
	\[
	\cR_{q} (\mu \mmid \pi ) \leq \eps^2\,.
	\]
\end{theorem}

\begin{remark}[Extension to arbitrary initialization]\label{rem:loacc-init}
	The algorithm in Theorem~\ref{thm:loacc} initializes at the Dirac distribution $ \delta_{x^*}$. This is reasonable because the cost of using gradient descent to compute $x^*$ approximately, using the same first-order oracle access, is dominated by the cost of subsequently running the sampling algorithm. If the algorithm is initialized at some other point $x$, the runtime only increases by a logarithmic factor of $\log W_{\psi_2}(\delta_{x},\pi)$, which is lower order unless $x$ is exponentially far from $x^*$, since $W_{\psi_2}(\delta_{x},\pi) \lesssim \sqrt{d/\alpha} + \|x - x^*\|$ by Lemma~\ref{lem:ulmc-init} and the triangle inequality.
\end{remark}

While the results of~\cite{cheng2018underdamped, dalalyanrioudurand2020underdamped, ma2021there, Matt23girsanov} have also shown iteration complexities that scale in the dimension $d$ as $\Otilde(d^{1/2})$, a key difference is that these results do not hold for R\'enyi divergences. In particular, past work has only proven weaker mixing results in the Wasserstein metric or the KL divergence. As discussed in the introduction, Wasserstein and KL guarantees are insufficient for the purpose of warm-starting high-accuracy sampling algorithms---for this, it is essential to have guarantees in the more stringent R\'enyi divergence. (Note that Wasserstein bounds are weaker than KL bounds by Talagrand's $T_2$ inequality, and moreover KL bounds are weaker than R\'enyi bounds by monotonicity of R\'enyi divergences.) See \S\ref{ssec:intro:tech} for a detailed discussion of this and of the many longstanding technical difficulties involved with establishing R\'enyi guarantees.

The algorithm we use is Underdamped Langevin Monte Carlo (ULMC) with certain parameters (stated explicitly in the proof in \S\ref{ssec:ulmc-proof}).
Background on this algorithm is recalled in \S\ref{ssec:ulmc-background}.
At a high level, the proof of Theorem~\ref{thm:loacc} uses the weak triangle inequality for R\'enyi divergences to decompose the sampling error of ULMC into the following two terms, both measured in R\'enyi divergence:
\begin{enumerate}
    \item The ``bias'' error between the stationary distribution of ULMC and the target distribution $\pi$.
    \item The ``discrete mixing'' error of ULMC to its biased stationary distribution.
\end{enumerate}
The bias error (1) is readily handled by recent results such as~\cite{GaneshT20,Matt23girsanov}. Bounding the discrete mixing error (2) is the key technical challenge; see \S\ref{ssec:intro:tech} for a detailed discussion of the technical obstacles related to this, and the connections to open problems about hypocoercivity in the PDE literature. The key contribution of this section is to bound this quantity in Theorem~\ref{thm:ulmc-discrete-mix} below. To do this, we use the new shifted divergence technique developed in \S\ref{sec:pabi}.
\par The section is organized as follows. In \S\ref{ssec:ulmc-background} we recall relevant background about ULMC, in \S\ref{ssec:ulmc-mix} we bound the discrete mixing error of ULMC, and in \S\ref{ssec:ulmc-proof} we use this to prove Theorem~\ref{thm:loacc}. Remark about notation in this section: the ULMC algorithm studied naturally operates on the augmented Hamiltonian state space $\R^{2d}$, so in this section we use boldface to denote probability distributions on $\R^{2d}$ from non-boldfaced distributions on $\R^d$. For example, we write $\bs{\pi}$ to denote the target distribution $\pi \otimes \cN(0,I_d)$ in this augmented space.

\subsection{Background on Underdamped Langevin Monte Carlo}\label{ssec:ulmc-background}

The exposition in this subsection is based on the corresponding section in~\cite{chewibook}; we refer the interested reader there for further details.

\paragraph*{Underdamped Langevin diffusion.} Studied since the work of Kolmogorov in the 1930s~\cite{Kol34}, the underdamped Langevin diffusion---sometimes also called the kinetic Langevin diffusion---is the solution to the stochastic differential equation
\begin{align}\label{eq:underdamped-diffusion}
\begin{aligned}
	\D X_t &= Y_t \, \D t\,, \\
	\D Y_t &= - \nabla f(X_t)\, \D t - \gamma Y_t \, \D t + \sqrt{2\gamma}\, \D B_t\,,
 \end{aligned}
\end{align}
where ${(B_t)}_{t\ge 0}$ is a standard $d$-dimensional Brownian motion. Analogous to the classical theory of convex optimization, here the auxiliary state variable $Y_t$ has the physical interpretation of momentum, and the linking parameter $\gamma$ has the physical interpretation of friction. A related interpretation of the underdamped Langevin diffusion is as a variant of the idealized Hamiltonian Monte Carlo algorithm, in which the momentum is refreshed continuously rather than periodically. The stationary distribution for this SDE is the joint distribution
\begin{align}
	{\bs \pi}(x,y) \propto \exp\Bigl( - f(x) - \frac{1}{2}\, \|y\|^2 \Bigr)\,.
	\label{eq:ulmc-pi}
\end{align}
In this section, we slightly abuse notation by using the boldface $\bs{\pi}$ to distinguish this joint distribution on $\R^{2d}$ from the target distribution $\pi \propto \exp(-f)$ on $\R^d$; of course the latter is the $x$-marginalization of the former. 

\par A major obstacle for analyzing the convergence of the underdamped Langevin diffusion is that this process exhibits hypocoercive dynamics, i.e., the standard Markov semigroup approach based on isoperimetric inequalities does not work. It is a longstanding question in PDE theory to develop general tools for establishing fast convergence of hypocoercive dynamics. See the discussion in \S\ref{ssec:intro:tech}. We bypass these issues by instead developing tools for analyzing a discrete-time version of this diffusion. (In forthcoming work, we detail the extent to which these discrete-time analyses enable analysis of continuous-time hypocoercive dynamics.) 

\paragraph*{Underdamped Langevin Monte Carlo.} There are several ways to discretize the underdamped Langevin diffusion. Perhaps the simplest way is the Euler--Maruyama discretization, as is standard for defining (unadjusted) Langevin Monte Carlo. However, for these underdamped Langevin dynamics, there is a better discretization which dates back at least to 1980~\cite{ermak1980numerical}, namely:
\begin{align*}
	\D X_t &= Y_t \, \D t\, \\
	\D Y_t &= -\nabla f(X_{nh}) \, \D t - \gamma Y_t \,\D t + \sqrt{2 \gamma} \,\D B_t\,,
\end{align*}
for $t \in [nh, (n+1)h]$. This process is called Underdamped Langevin Monte Carlo (ULMC) or Kinetic Langevin Monte Carlo; we use the former term in this paper. The point of this discretization is that since the gradient is refreshed periodically rather than continuously, the SDE is linear within these periods, and thus can be integrated exactly in closed form (see, e.g.,~\cite[Appendix A]{cheng2018underdamped}).
This is called an ``exponential integrator'' in the lingo of numerical analysis, and is formalized as follows.

\begin{lemma}[Explicit Gaussian law for ULMC iterates]\label{lem:ulmc-law}
	Conditioned on $(X_{nh}, Y_{nh})$, the law of $(X_{(n+1)h}, Y_{(n+1)h})$ is the Gaussian distribution $\cN(F(X_{nh}, Y_{nh}), \Sigma \otimes I_d)$ where 
	\begin{align*}
		F(x,y) \deq \left( 
		x + \gamma^{-1}\, (1 - a)\, y - \gamma^{-1}\, (h - \gamma^{-1}\,(1 - a)) \,\nabla f(x),\;
		a y - \gamma^{-1}\,(1 - a)\, \nabla f(x)\right)
\end{align*}
	and
	\begin{align*}
		\Sigma 
              \deq
		\begin{bmatrix}
			\frac{2}{\gamma}\, (h - \frac{2}{\gamma}\, (1 - a) + \frac{1}{2\gamma}\,(1 - a^2))
			&
			\frac{2}{\gamma}\, ( \frac{1}{2} - a + a^2) 
			\\ 
			\frac{2}{\gamma}\,( \frac{1}{2} - a + a^2)
			&
			1-a^2
		\end{bmatrix}\,.
	\end{align*}
	Above, we use the notational shorthand $a \deq \exp(-\gamma h)$.
\end{lemma}

In the rest of this section, we write $\bs P \deq \bs P_{h,\gamma}$ to denote the Markov transition kernel on $\R^{2d}$ that corresponds to an iteration of ULMC. We suppress the dependence of $\bs P$ on the parameters $h$ and $\gamma$ for simplicity of notation.

\subsection{Discrete mixing of Underdamped Langevin Monte Carlo}\label{ssec:ulmc-mix}

\begin{theorem}[Discrete mixing of ULMC]\label{thm:ulmc-discrete-mix}
    Suppose that $\pi \propto \exp(-f)$ where $f$ is $\alpha$-strongly-convex and $\beta$-smooth.
	Let $\bs P$ denote the Markov transition kernel for ULMC when run with friction parameter $\gamma = \sqrt{2\beta}$ and step size $h \lesssim 1/(\kappa \sqrt{\beta})$. Then, for any target accuracy $0 < \eps \le \sqrt{\frac{\log 2}{q-1}}$, any R\'enyi order $q \geq 1$, 
	and any two initial distributions 
 $\bs{\mu_0},\bs{\mu_0'} \in \cP(\R^{2d})$,
	\[
	\cR_{q}( \bs{\mu_0} \bs P^N \mmid \bs{\mu_0'} \bs P^N) \leq \eps^2 \,,
	\]
	if the number of ULMC iterations is
	\[
	N \gtrsim
	\frac{\sqrt{\beta}}{\alpha h} \log  \Bigl( \frac{q\, W_{\psi_2}^2( \cM_{\#} \bs{\mu_0}, \cM_{\#} \bs{\mu_0'})}{\beta^{1/2}\, \eps^2\, h^3} \Bigr)\,,
	\]
	where $\cM$ is the linear map defined in~\eqref{eq:ulmc-twisting}.
 \end{theorem}

 We remark that in a typical use case of this discrete mixing result, $\bs{\mu_0}$ is initialized at a product distribution of the form $\mu_0 \otimes \cN(0,I_d)$, and is compared to the target distribution $\bs{\mu_0'} = \bs{\pi} = \pi \otimes \cN(0,I_d)$. In this setting, the Orlicz--Wasserstein metric in the upper bound can be simplified to $W_{\psi_2}( \cM_{\#} \bs{\mu_0}, \cM_{\#} \bs{\mu_0'}) \leq 2 W_{\psi_2}(\mu_0, \pi)$.

To prove Theorem~\ref{thm:ulmc-discrete-mix}, we appeal to our new shifted divergence technique developed in \S\ref{sec:pabi}. 
This requires analyzing the ULMC iterates in a \emph{twisted norm}, since an iteration of the ULMC algorithm (or more precisely, the mean-shifting function $F$ defined in Lemma~\ref{lem:ulmc-law}) is not contractive with respect to the standard Euclidean norm. This twisted norm is the Euclidean norm after the change of coordinates
\begin{align}
	(u,v)
	&\deq \mc M(x, y)
	\deq \bigl(x, x + \frac{2}{\gamma}\,y\bigr)\,.
	\label{eq:ulmc-twisting}
\end{align}
In these new coordinates, the mean of the next iterate of ULMC started at $(u,v)$ is $\bar F(u,v)$, where $\bar F = \mc M \circ F\circ \mc M^{-1}$.
Since $\mc M^{-1}(u,v) = (u, \frac{\gamma}{2}\,(v - u))$, we can explicitly write
\begin{align}
	\bar F(u,v)
	&= \Bigl( u + \frac{1-a}{2}\,(v - u)- \frac{h - \gamma^{-1}\,(1-a)}{\gamma}\,\nabla f(u),\; u + \frac{1+a}{2}\,(v-u) - \frac{h+\gamma^{-1}\,(1-a)}{\gamma}\,\nabla f(u)\Bigr)\,.
	\label{eq:ulmc-barF}
\end{align}
By Lemma~\ref{lem:ulmc-law}, it follows that conditioned on $(U_{nh}, V_{nh})$, the law of $(U_{(n+1)h}, V_{(n+1)h})$ is the Gaussian distribution $\cN(\bar F (U_{nh}, V_{nh}), \bar \Sigma \otimes I_d)$ where
\begin{align}
	\bar \Sigma = \mc M \Sigma \mc M^\T\,.
		\label{eq:ulmc-barSigma}
\end{align}

We make use of the following two helper lemmas about the dynamics of ULMC in this twisted norm. The first helper lemma shows that the ULMC algorithm sends two iterates to Gaussians with means that are closer in the twisted norm than the original iterates. Since the two Gaussians have the same covariance $\bar \Sigma \otimes I_d$, this implies that the ULMC Markov transition kernel $P$ is contractive in $W_{\psi_2}$, which will allow us to use our new shifted divergence technique from \S\ref{sec:pabi}. This lemma first appeared in the recent paper~\cite{Matt23girsanov}; for completeness, a proof is provided in Appendix~\ref{app:pf-lem:ulmc-contract}.

\begin{lemma}[Contractivity of ULMC in the twisted  norm]\label{lem:ulmc-contract}
	Suppose that $\pi \propto \exp(-f)$ where $f$ is $\alpha$-strongly-convex and $\beta$-smooth.
	For step size $h \lesssim 1$ and friction parameter $\gamma = \sqrt{2\beta}$, the function $\bar{F} : \R^{2d} \to \R^{2d}$ defined in~\eqref{eq:ulmc-barF} is a contraction with
	\[
		\|\bar{F}\|_{\rm Lip} \leq 1 - \frac{\alpha}{\sqrt{2\beta}}\,h + O(\beta h^2)\,.
	\]
\end{lemma}

The second helper lemma estimates the noise of ULMC in this twisted norm. The proof is an explicit computation and is provided in Appendix~\ref{app:pf-lem:ulmc-noise}.

\begin{lemma}[Noise of ULMC in the twisted norm]\label{lem:ulmc-noise}
Suppose that $h \lesssim 1/\gamma$. Then the matrix $\bar{\Sigma}$ defined in~\eqref{eq:ulmc-barSigma} satisfies
	\[
		\lambda_{\min}( \bar\Sigma) =
		\frac{\gamma h^3}{6}\, \bigl( 1 - O(\gamma h)\bigr)\,.
	\]
\end{lemma}

Armed with Lemmas~\ref{lem:ulmc-contract} and~\ref{lem:ulmc-noise}, we are now ready to prove Theorem~\ref{thm:ulmc-discrete-mix}.

\begin{proof}[Proof of Theorem~\ref{thm:ulmc-discrete-mix}]
	Let $\{\bs{\mu_n}\}_{n \geq 0}$ and $\{\bs{\mu_n'}\}_{n \geq 0}$ denote the two processes $\bs{\mu_n} = \bs{\mu_0} \bs P^n$ and $\bs{\mu_n'} = \bs{\mu_0'} \bs P^n$ obtained by running ULMC from initialization distributions $\bs{\mu_0}$ and $\bs{\mu_0'}$, respectively. Define twisted processes $\{\bs{\nu_n}\}_{n \geq 0}$ and $\{\bs{\nu_n'}\}_{n \geq 0}$ by $\bs{\nu_n} = \cM_{\#} \bs{\mu_n}$ and $\bs{\nu_n'} = \cM_{\#} \bs{\mu_n'}$, where $\cM$ is the change-of-coordinates matrix defined in~\eqref{eq:ulmc-twisting}. Since $\cM$ is invertible, applying the data-processing inequality for R\'enyi divergences (Lemma~\ref{lem:renyi-process}) in both directions implies
	\[
	\cR_q(\bs{\mu_N} \mmid \bs{\mu_N'})
	=
	\cR_q( \bs{\nu_N} \mmid \bs{\nu_N'} )\,.
	\]
	\par We now show that the latter term is at most $\eps^2$. For shorthand, define $\lambda  \deq  \lambda_{\min}(\bar\Sigma)$ and define $\bs{Q}$ to be the Markov operator given by $\bs{\nu} \bs{Q} = \bar{F}_{\#} \bs{\nu} \ast \cN(0,\bar{\Sigma} \otimes I_d -\lambda I_{2d})$. Then by Lemma~\ref{lem:ulmc-law} and a change of measure, the law of $\bs{\nu_{n+1}}$ is 
	\[
	\bar{F}_{\#} \bs{\nu_n} \ast \cN(0, \bar\Sigma \otimes I_d)
	=
	\bs{\nu_n} \bs{Q} \ast \cN(0, \lambda I_{2d})\,.
	\] 
	Similarly, the law of $\bs{\nu_{n+1}'}$ is 
	\[
	\bar{F}_{\#} \bs{\nu_n'} \ast \cN(0, \bar\Sigma \otimes I_d)
	=
	\bs{\nu_n'} \bs{Q} \ast \cN(0, \lambda  I_{2d})\,.
	\]
	Thus, letting $c$ denote the Lipschitz constant of the Markov operator $\bs Q$ w.r.t.\ the $W_{\psi_2}$ metric, we may invoke\footnote{The application of this result requires the number of iterations $N$ to be large enough that the right hand side of~\eqref{eq:pf-ulmc-mixing:bound}, later set to $\eps^2$, is at most $(\log 2)/(q-1)$. But this holds by assumption.}
 the new shifted divergence result (Theorem~\ref{thm:pabi-orlicz}) to bound
	\begin{align}
		\cR_q( \bs{\nu_N} \mmid \bs{\nu_N'})
		\leq c^{2N}\, \frac{q}{2\lambda}\, W_{\psi_2}^2( \bs{\nu}_0,\bs{\nu}_0')\,.
		\label{eq:pf-ulmc-mixing:bound}
	\end{align}
	\par We now use the two helper lemmas to quantify the various terms in~\eqref{eq:pf-ulmc-mixing:bound}. First, by a simple coupling argument and then an application of Lemma~\ref{lem:ulmc-contract}, the Markov operator $\bs Q$ is $W_{\psi_2}$-contractive with parameter $c$, where
	\[
	c \leq \|\bar F\|_{\rm Lip} \leq \exp\Bigl( - 
	\Omega\Bigl( \frac{\alpha h}{\sqrt{\beta}} \Bigr) \Bigr)\,.
	\]
	Second, because ULMC is run with 
	step size $h \lesssim 1/\gamma$, Lemma~\ref{lem:ulmc-noise} implies
	\[
	\lambda = \Omega\bigl( \sqrt{\beta}\,h^3 \bigr)\,.
	\]
	\par Therefore, by combining the above displays, we conclude that
	\[
	\cR_q( \bs{\mu_N} \mmid \bs{\mu_N'} )
	=
	\cR_q( \bs{\nu_N} \mmid \bs{\nu_N'} )
	\lesssim
	\exp\biggl( - \Omega\Bigl( \frac{\alpha h}{\sqrt{\beta}}\, N \Bigr)\biggr) \,
	\frac{q}{\sqrt{\beta}\, h^3}\,  W_{\psi_2}^2( \cM_{\#} \bs{\mu_0},  \cM_{\#} \bs{\mu_0'})\,.
	\]
	Setting this bound to $\eps^2$ and solving for $N$ completes the proof.
\end{proof}

\subsection{Warm start with Underdamped Langevin Monte Carlo 
}\label{ssec:ulmc-proof}

We begin by bounding the distance to the target at initialization. 
We emphasize that this initialization is \emph{not} a warm start. 
Indeed, this initialization is even weaker than what is typically called a ``feasible start'' in the literature (namely $\cN(x^*,\beta^{-1} I_d)$), and moreover can be further relaxed to an arbitrary initialization $x_0$ so long as the distance between $x_0$ and the mode $x^*$ of the target distribution is sub-exponentially large (since our final bound depends only logarithmically on this distance).

\begin{lemma}[Orlicz--Wasserstein distance at initialization]\label{lem:ulmc-init}
	Suppose that $\pi \propto \exp(-f)$ where $f$ is $\alpha$-strongly convex.
	Let $x^*$ denote the minimizer of $f$. Then
	\[
	W_{\psi_2}(\delta_{\xopt},\pi) \leq 6 \sqrt{d/\alpha}\,.
	\]
\end{lemma}
\begin{proof}
	Let $X \sim \pi$ and define $Y  \deq  \|X - x^*\|$. By definition of the Orlicz--Wasserstein metric,
	\[
	W_{\psi_2}( \delta_{x^*}, \pi )
	=
	\|X - x^*\|_{\psi_2}
	=
	\|Y\|_{\psi_2}
	=
	\inf \biggl \{ \lambda > 0 \; : \; \E \exp\Bigl(  \frac{Y^2}{\lambda^2}\Bigr) \leq 2 \biggr \}\,.
	\]
	By the elementary inequality $(a+b)^2 \leq 2\,(a^2 + b^2)$, we can bound $Y^2 = (Y - \E Y + \E Y)^2  \leq 2\, \E[Y]^2 + 2\,(Y - \E Y)^2$. Thus
	\begin{align*}
		\E \exp \Bigl( \frac{Y^2}{\lambda^2} \Bigr)
		&\le  \exp\Bigl( \frac{2\,\E[Y]^2}{\lambda^2} \Bigr) \cdot \E \exp \Bigl( \frac{2\,(Y - \E Y)^2}{\lambda^2} \Bigr)\,.
	\end{align*}
	For the first term, use the basic inequality $\E[Y]^2 \leq \E[Y^2]$ and a standard second-moment-type bound for strongly log-concave distributions (Lemma~\ref{lem:lsi-second-moment}) to obtain
\[
\exp\Bigl( \frac{2\,\E[Y]^2}{\lambda^2} \Bigr)
\leq
\exp\Bigl( \frac{2d}{\alpha \lambda^2} \Bigr)\,.
\]
	For the second term, we use sub-Gaussian concentration. 
	Specifically, combining the log-Sobolev inequality on $\pi$ 
	with the fact that $Y = \|X - x^*\|$ is a $1$-Lipschitz function of $X \sim \pi$, it follows that $Y - \E Y$ is sub-Gaussian with variance proxy $1/\alpha$ (see, e.g.,~\cite[equation (5.4.1)]{bakry2014analysis} and~\cite[\S 2.3]{boucheronlugosimassart}). Therefore, by a standard moment generating function bound for the square of a sub-Gaussian random variable (see, e.g.,~\cite[\S 2.3]{boucheronlugosimassart}), 
	\[
		\E \exp \Bigl( \frac{2\,(Y - \E Y)^2}{\lambda^2} \Bigr)
		\leq
		2^{16/(\alpha \lambda^2)}\,,
	\]
	for any $\lambda^2 \geq 16/\alpha$. By combining the above displays and setting $\lambda = \sqrt{32d/\alpha}$, we conclude that
	\[
			\E \exp \Bigl( \frac{Y^2}{\lambda^2} \Bigr)
			\leq
			\exp\Bigl( \frac{2d}{\alpha \lambda^2} + \frac{16\log 2}{\alpha \lambda^2} \Bigr)
			\leq
			\exp(
				\log 2
			)
			=2\,.
	\]
	Therefore this choice of $\lambda$ is an upper bound on the Orlicz--Wasserstein norm $W_{\psi_2}(\delta_{x^*}, \pi)$. 
\end{proof}

\par The second lemma uses Girsanov's theorem to bound the bias of ULMC.
Here, we build upon recent advances in the literature on R\'enyi discretization of stochastic processes.
Beginning with the works~\cite{GaneshT20, chewi2021optimal, erdhoszha22chisq} and culminating in the paper~\cite{chewi2021analysis}, it is now understood that the Girsanov discretization technique leads to bias bounds for LMC in R\'enyi divergence matching prior results which only held for weaker divergences, and yet remains flexible enough to cover varying assumptions.
The recent paper~\cite{Matt23girsanov} extends this technique for ULMC\@.
Since the results of~\cite{Matt23girsanov} hold under more general assumptions at the expense of a more involved analysis, and in the interest of keeping our derivations more self-contained, in Appendix~\ref{app:loacc-girsanov} we simplify and streamline the Girsanov argument of~\cite{Matt23girsanov} for our setting of interest. In order to clarify where the $d^{1/2}$ comes from, we write the final bound in terms of the total elapsed continuous time $T = Nh$ rather than the number of iterations $N$.

\begin{lemma}[Bias of ULMC]\label{lem:ulmc-bias}
	Suppose that $f$ is $\alpha$-strongly-convex and $\beta$-smooth.
	Let $\bs{\pi}(x,y) \propto \exp(-f(x) - \half \|y\|^2)$, and let $\bs P$ denote the Markov transition kernel corresponding to an iteration of ULMC with friction parameter $\gamma \asymp \sqrt\beta$ and step size $h \lesssim \frac{1}{\beta^{3/4} d^{1/2} q\, {(T\log N)}^{1/2}}$, where $N$ is the total number of iterations and $T = Nh$ is the total elapsed time.
 Then,
	\begin{align*}
		\cR_q(\bs{\pi} \bs P^N \mmid \bs{\pi} )
		&\lesssim \beta^{3/2} d h^2 q T\,.
	\end{align*}
\end{lemma}

\begin{proof}[Proof of Theorem~\ref{thm:loacc}]
	Consider the following algorithm: run ULMC 
	for $N$ iterations 
	from initialization $\bs{\mu_0} = \delta_{x^*} \otimes \cN(0,I_d)$ to obtain an iterate $(X,Y) \in \R^{2d}$, and then output $X$. We show that for a certain setting of the ULMC parameters, this algorithm satisfies the guarantees of the theorem.
	\par To this end, given a joint distribution $\bs{\nu} \in \cP(\R^{2d}) \cong \cP(\R^{d} \times \R^d)$, let $(\Pi^x)_\#\bs{\nu} \in \cP(\R^d)$ denote the marginal on the first $d$ coordinates, where $\Pi^x : \R^{2d}\to\R^d$ maps $(x,y)\mapsto x$. Then in particular, the law of the algorithm's output is $\mu  \deq (\Pi^x)_\# (\bs{\mu_0} \bs P^N)$, and the target distribution is $\pi = (\Pi^x)_\# \bs{\pi}$. Thus by the data-processing inequality for R\'enyi divergences (Lemma~\ref{lem:renyi-process}), the sampling error of the algorithm is bounded by
	\[
		\cR_q( \mu \mmid \pi)
		\leq
		\cR_q( \bs{\mu_0} \bs P^N \mmid \bs{\pi} )\,.
	\]
	By the weak triangle inequality for R\'enyi divergence (Lemma~\ref{lem:renyi-triangle}),  we can further bound this by
	\begin{align}
		\cR_q ( \bs{\mu_0} \bs P^N \mmid \bs{\pi} )
		\leq
		\frac{q-1/2}{q-1} \, \cR_{2q}( \bs{\mu_0} \bs P^N \mmid  \bs{\pi} \bs P^n )
		+
		\cR_{2q-1}( \bs{\pi} \bs P^N \mmid \bs{\pi} )\,.
		\label{eq:pf-ulmc:triangle}
	\end{align}
	The coefficient $(q-1/2)/(q-1)$ can be crudely bounded by $2$, say, since it suffices to bound the R\'enyi divergence error for $q \geq 3/2$ (indeed, monotonicity of R\'enyi divergences in the order $q$ then implies the same bound for $q < 3/2$). 
	\par Now, by combining our discrete mixing result for ULMC (Theorem~\ref{thm:ulmc-discrete-mix}), our initialization bound (Lemma~\ref{lem:ulmc-init}), and the ULMC bias bound (Lemma~\ref{lem:ulmc-bias}), we conclude that
	\begin{align}
		\cR_q ( \bs{\mu_0} \bs P^N \mmid \bs{\pi} ) \leq \eps^2\,,
	\end{align}
	if ULMC is run with friction parameter $\gamma$, step size $h$, and iteration complexity $N$ that satisfy:
	\[
		\gamma = \sqrt{2\beta}
		\qquad \text{and} \qquad 
		h \lesssim \frac{\eps}{\beta^{3/4} d^{1/2} q^{1/2} T^{1/2}}
		\qquad \text{and} \qquad 
		N \gtrsim \frac{\sqrt{\beta}}{\alpha h} \log \Bigl( \frac{dq}{ \alpha \beta^{1/2} \varepsilon^2 h^3} \Bigr)\,.
	\] 
	By recalling that $T \defeq Nh$, solving for these choices of parameters, and omitting logarithmic factors, we conclude that it suffices to run ULMC with the following choices of parameters:
	\[
		\gamma = \sqrt{2\beta}
		\qquad \text{and} \qquad 
		h = 
		\widetilde{\Theta}\Bigl( \frac{\eps \alpha^{1/2}}{\beta d^{1/2} q^{1/2}} \Bigr) 
		\qquad \text{and} \qquad 
		N = \widetilde{\Theta}\Bigl( \frac{\kappa^{3/2} d^{1/2} q^{1/2} }{\eps} \Bigr)\,. \qedhere
	\] 
\end{proof}

%% file: sections/hiacc.tex
\section{High-accuracy sampling with \texorpdfstring{$O(\sqrt{d})$}{sqrt(d)} complexity}\label{sec:hiacc}

Establishing fast mixing results for MALA is a longstanding problem. As detailed in \S\ref{sec:intro}, recent breakthroughs have made it clear that the key barrier for fast mixing of MALA is the question of warm starts. In this section, we use the faster low-accuracy sampling result developed in \S\ref{sec:ulmc} to efficiently warm start MALA. This leads to the fastest known high-accuracy sampling algorithms not only in strongly log-concave settings (details in \S\ref{ssec:hiacc-slc}), but also in weakly-log-concave and isoperimetric, non-log-concave settings (details in \S\ref{ssec:hiacc-ext}), for which we improve over state-of-the-art query complexity results by a factor of $\sqrt{d}$.

\subsection{Strongly-log-concave setting}\label{ssec:hiacc-slc}

Here we improve the query complexity for high-accuracy sampling from a strongly log-concave target distribution to $\Otilde(\kappa d^{1/2} \, \polylog(1/\eps))$. This result is formally stated as follows.

\begin{theorem}[Faster high-accuracy sampler for well-conditioned targets]\label{thm:main-slc}
	Suppose that $\pi \propto \exp(-f)$ where $f$ is $\alpha$-strongly-convex and $\beta$-smooth. There is an algorithm with randomized runtime that, given knowledge of the minimizer of $f$ and access to $N$ first-order queries for $f$, outputs a random point in $\R^d$ with law $\mu$ satisfying
	$
		 \msf{d}(\mu \mmid \pi) \leq \eps
	$
	for any of the following metrics:
	\[
		\msf{d} \in \{ \TV, \sqrt{\KL}, \sqrt{\chi^2}, \sqrt{\alpha}\, W_2 \}\,.
	\]
    Moreover, for any $\delta \in (0,1)$ with probability at least $1-\delta$, the number of queries made satisfies
    \begin{align*}
        N
        &\le \widetilde O\Bigl(\kappa d^{1/2}\, \log^4 \max\Bigl\{\frac{1}{\varepsilon}, \, \log \frac{1}{\delta}\Bigr\} \Bigr)\,.
    \end{align*}
\end{theorem}

In analogy to the familiar concept from algorithm design for deterministic problems~\cite{CLRS}, the algorithm in Theorem~\ref{thm:main-slc} may be called a ``Las Vegas'' algorithm because it has a randomized runtime which is small with high probability. The fact that this runtime is randomized is not an issue in practice because the iteration complexity depends on a quantity that is efficiently estimable during the execution of the algorithm.

In the rest of this subsection, we overview the algorithm in  Theorem~\ref{thm:main-slc} and its analysis; see Appendix~\ref{app:hiacc} for full technical details. This algorithm combines three algorithms as building blocks: ULMC, MALA, and the proximal sampler algorithm. Let us explain this by building up to the full complexity in two steps---both because this will motivate why all three algorithmic components are needed, and also because this is how our analysis actually proceeds.

\paragraph*{Weak version of Theorem~\ref{thm:main-slc} (full details in Appendix~\ref{app:highacc:slc-weak}).} First, consider the following simplified version of the algorithm in Theorem~\ref{thm:main-slc} which is only comprised of two algorithmic components: run ULMC, and then use this as a warm start for MALA. In order to argue that this two-phase algorithm mixes rapidly, we crucially use our result from \S\ref{sec:ulmc} to guarantee that ULMC mixes to constant R\'enyi divergence error in a number of iterations that scales in the dimension $d$ as $\Otilde(d^{1/2})$ rather than $\Otilde(d)$. This allows us to provide an algorithm for warm-starting MALA which is not significantly slower than MALA is when initialized from a warm start. In other words, this lets us exploit, for the first time, the recent breakthroughs on MALA~\cite{chewi2021optimal,wuschche2022minimaxmala} which show that the mixing time of MALA scales in the dimension as $\Otilde(d^{1/2})$ rather than $\Otilde(d)$ when it is initialized at a warm start rather than a feasible start. We recall that such an improvement \emph{cannot} be obtained without the warm start, due to the lower bound of~\cite{lee2021lower}.
This leads to a final runtime of roughly $\Otilde(\kappa^{3/2} d^{1/2} + \kappa d^{1/2} \, \polylog(1/\eps))$; a formal statement is given in  Theorem~\ref{thm:main-slc-weak}.
\par However, while this simple combination of ULMC and MALA achieves the desired dependence on the dimension $d$, it leads to a suboptimal dependence on the condition number $\kappa$, namely $\Otilde(\kappa^{3/2})$ rather than $\Otilde(\kappa)$. This worsened dependence in $\kappa$ arises from the state-of-the-art bounds on the discretization of ULMC~\cite{Matt23girsanov}. For full details on this weak version of Theorem~\ref{thm:main-slc}, see Appendix~\ref{app:highacc:slc-weak}.

\paragraph*{The full version of Theorem~\ref{thm:main-slc} (full details in Appendix~\ref{app:highacc:slc-main}).} In order to improve the condition number dependence of the weak version of Theorem~\ref{thm:main-slc}, we require an extra algorithmic component: the recently proposed proximal sampler algorithm. See Appendix~\ref{app:hiacc:ps} for background on this proximal sampler algorithm. Briefly, this algorithm reduces the problem of sampling a strongly log-concave distribution with condition number $\kappa$, to the the problem of sampling $\widetilde O(\kappa)$ related strongly log-concave distributions each with constant condition number. The upshot is that the latter can be accomplished in the desired $\Otilde(\kappa d^{1/2} \, \polylog(1/\eps))$ runtime by using the weak version of the algorithm since each sampling subproblem is well-conditioned. 

\paragraph*{Mixing in R\'enyi divergence.} While the main conceptual innovation here is the high-level strategy of combining these three algorithmic building blocks, we remark that an additional technical obstacle for proving Theorem~\ref{thm:main-slc} is showing mixing in more stringent notions of distance than TV. See the discussion in \S\ref{ssec:intro:tech}. Indeed, while our new ULMC result proves fast mixing in R\'enyi divergence, existing results on MALA and its combination with the proximal sampler are limited to TV. 
We boost this mixing in TV to R\'enyi divergences (and thus all the other desired metrics by standard comparison inequalities) using two additional ideas.

\par The first improves mixing bounds for the proximal sampler from TV to R\'enyi divergence.
To do this, we control the propagation of error when each step of the proximal sampler algorithm is performed approximately in R\'enyi divergence. As we show, this is readily accomplished by appealing to the ``strong composition rule'' of R\'enyi divergences from the differential privacy literature. 

\par The second improves mixing bounds for MALA from TV to R\'enyi divergence. We accomplish this by further exploiting the fact that MALA is warm started in R\'enyi divergence. Note that this means we use the R\'enyi warm start in two ways: first to show that MALA mixes fast in TV, which is what we can conclude from the above argument and appealing to~\cite{chewi2021optimal,wuschche2022minimaxmala}; and second, to boost the TV bound at the final iterate to a more stringent bound. We isolate this TV-to-R\'enyi boosting technique in the following simple lemma as it may be of independent interest: indeed, since it uses the TV mixing bound in an entirely black-box way, this lemma may be useful for establishing R\'enyi mixing of other warm-started algorithms.
This lemma improves the previous result of~\cite[Lemma 28]{chewi2021optimal} because that result required a warm-start in $\cR_{\infty}$ which is currently unavailable algorithmically, whereas this lemma here only requires the weaker condition of a warm start in a R\'enyi divergence of finite order (stated here with $q=3$ for simplicity).

\begin{lemma}[Boosting TV to R\'enyi mixing given R\'enyi warm start]\label{lem:warm-boost-tv}
	Let $P$ be a Markov transition kerrnel which has stationary distribution $\pi$. Consider running $P$ from any initialization distribution for $N$ steps to obtain a distribution $\mu_N \defeq \mu_0 P^N$. Then
	\[
	\chi^2(\mu_N \mmid \pi) \leq \sqrt{ \TV (\mu_N,\pi) \cdot \bigl( \exp(2 \cR_3(\mu_0 \mmid \pi)) + 1 \bigr) } \,.
	\]
\end{lemma}

See Appendices~\ref{app:hiacc:mala} and~\ref{app:hiacc:ps} for background on MALA and the proximal sampler algorithm, respectively; and see Appendices~\ref{app:highacc:slc-weak} and~\ref{app:highacc:slc-main} for proofs of the weak version and full version of Theorem~\ref{thm:main-slc}, respectively.

\subsection{Extensions to weakly-log-concave and non-log-concave settings}\label{ssec:hiacc-ext}

Our faster algorithm for sampling from well-conditioned targets (Theorem~\ref{thm:main-slc}) yields faster samplers for a variety of other settings, due essentially to the reductions in~\cite{chenetal2022proximalsampler}. We present here several such extensions that
concern target distributions which satisfy isoperimetric inequalities, 
which is quite flexible in the sense that this allows for non-log-concavity and also is preserved under, e.g., bounded perturbations and Lipschitz mappings. 
See \S\ref{ssec:prelim-functional} for background on these isoperimetric inequalities. 
\par A comment on notation for these isoperimetric settings: we still use the condition number $\kappa$ to denote the ratio $\kappa = \beta/\alpha$, but now $\alpha$ denotes the (inverse) parameter of an isoperimetric bound, rather than the parameter for strong convexity. The motivation behind this notation is that $\alpha$-strong-convexity implies the log-Sobolev inequality with parameter $1/\alpha$, which in turn implies the Poincar\'e inequality with parameter $1/\alpha$ (see Lemma~\ref{lem:bakry-emery}).

In this section, for simplicity of exposition, in addition to the first-order oracle for $f$ we also assume access to a prox oracle which can compute the proximal operator for $f$ with step size $h = \frac{1}{2\beta}$, namely for any $y \in \R^d$ the oracle returns return $\prox_{hf}(y) \deq \argmin_{x\in\R^d}\{f(x) + \frac{1}{2h}\,\norm{y-x}^2\}$.
Note that for this choice of step size, computing the proximal operator is a strongly convex and smooth optimization problem with condition number $O(1)$, so this can be done with off-the-shelf optimization methods such as gradient descent.
We emphasize, however, that this assumption is only made to ease the presentation of the results, and we provide more detailed results without this assumption in Appendix~\ref{app:hiacc} that more closely mirror Theorem~\ref{thm:main-slc}. 

\begin{theorem}[Faster high-accuracy sampling from LSI targets]\label{thm:main-lsi}
	Suppose that $\pi \propto \exp(-f)$ satisfies $1/\alpha$-LSI and that $f$ is $\beta$-smooth. There is an algorithm that, given access to a first-order + prox oracle for $f$ and initialized at $\mu_0$, outputs a random point in $\R^d$ with law $\mu$ satisfying
	$
		 \msf{d}(\mu \mmid \pi) \leq \eps
	$
	for any of the following metrics:
	\[
		\msf{d} \in \{ \TV, \sqrt{\KL}, \sqrt{\chi^2}, \sqrt{\alpha}\, W_2 \}\,,
	\]
 using at most
 \begin{align*}
N = \Otilde\biggl(\kappa d^{1/2} \, \log\Bigl(\frac{\cR_2(\mu_0 \mmid \pi)}{\eps^2}\Bigr)\,  \log^3\Bigl(\frac{1}{\eps}\Bigr) \biggr)\qquad\text{queries}\,.
 \end{align*}
\end{theorem}

\begin{theorem}[Faster high-accuracy sampling from PI targets]\label{thm:main-poincare}
    Suppose that $\pi \propto \exp(-f)$ satisfies $1/\alpha$-PI and that $f$ is $\beta$-smooth. There is an algorithm that, given access to a first-order + prox oracle for $f$ and initialized at $\mu_0$, outputs a random point in $\R^d$ with law $\mu$ satisfying
	$
		 \msf{d}(\mu \mmid \pi) \leq \eps
	$
	for any of the following metrics:
	\[
		\msf{d} \in \{ \TV, \sqrt{\KL}, \sqrt{\chi^2}, \sqrt{\alpha}\, W_2 \}\,,
	\]
 using at most
    \begin{align*}
        N
	\le \Otilde\biggl( \kappa d^{1/2} \, \log\Bigl( \frac{\chi^2(\mu_0 \mmid \pi)}{\varepsilon^2}\Bigr)\biggr)\qquad\text{queries}\,.
    \end{align*}
\end{theorem}

\begin{remark}[Extensions to Lata\l{}a--Oleszkiewicz targets]
    Just as in~\cite{chenetal2022proximalsampler}, we could also obtain a result for distributions satisfying a Lata\l{}a--Oleszkiewicz inequality, which interpolates between PI and LSI\@. In this setting, we again improve over the previous state-of-the-art bounds by a factor of $d^{1/2}$. However, for the sake of brevity, we omit this extension as it is conceptually similar but requires more involved technical details.
\end{remark}

These results are the direct analogs of the state-of-the-art complexity results in~\cite[Corollary 7]{chenetal2022proximalsampler}, but here with a dimension dependence that is improved by a factor of $d^{1/2}$. 

\par We mention another consequence of our improved high-accuracy sampler for the strongly-log-concave setting. Namely, via the same proximal reduction framework, this gives the following alternative complexity bound for target distributions which are (non-strongly) log-concave, sometimes called weakly-log-concave. This bound is a direct analog of~\cite[Corollary 6]{chenetal2022proximalsampler}, but here with a dimension dependence that is also improved by a factor of $d^{1/2}$. Note that this theorem is a low-accuracy guarantee; one can also obtain high-accuracy samplers from our results in this log-concave setting by using the fact that log-concavity implies a Poincar\'e inequality, albeit with a function-dependent constant~\cite{kannan1995isoperimetric}, and then appealing to Theorem~\ref{thm:main-poincare}. The resulting low-accuracy and high-accuracy results are incomparable in the sense that each can dominate in different settings---but in any case, our theorems for both settings yield improvements by a factor of $d^{1/2}$.

\begin{theorem}[Faster low-accuracy sampling from log-concave targets]\label{thm:main-logconcave}
    Suppose that $\pi \propto \exp(-f)$, where $f$ is convex and $\beta$-smooth. There is an algorithm that, given access to a first-order + prox oracle for $f$ and initialized at $\mu_0$, outputs a random point in $\R^d$ with law $\mu$ satisfying $\KL(\mu \mmid \pi) \le \eps^2$, using at most
    \begin{align*}
        N
    &\le \Otilde\biggl( \frac{\beta d^{1/2}\, W_2^2(\mu_0, \pi)}{\eps^2} \biggr)\qquad\text{queries}\,.
    \end{align*}
\end{theorem}

Proofs for the results in this section are provided in Appendix~\ref{app:hiacc}. At a high level, the proof of all these results use the same reduction to the problem of sampling from well-conditioned distributions. This reduction is based on the proximal sampler (described in Appendix~\ref{app:hiacc:ps}) and lets us apply our improved sampler for the well-conditioned case (Theorem~\ref{thm:main-slc}).
In each case, however, we must track the propagation of error due to the inexact implementation of the backwards step of the proximal sampler, which was not previously done in any work except for in the TV distance.

In Appendix~\ref{app:hiacc:explicit}, we provide more explicit, albeit more complicated, statements of these results to address the following two points.
(1) The above results depend on the initialization (through $W_2(\mu_0,\pi)$, $\cR_2(\mu_0\mmid\pi)$, or $\chi^2(\mu_0\mmid \pi)$) and it may be unclear how large these quantities are in a given application.
(2) We assumed that the algorithm has access to a stronger oracle than just a first-order oracle for $f$, namely, we also assumed access to a prox oracle for $hf$ with $h=\frac{1}{2\beta}$.
We address (1) by explicitly bounding these initialization quantities in terms of other, more easily computable problem parameters, and we address (2) by removing the assumption of a prox oracle.

%% file: sections/discussion.tex
\section{Discussion}\label{sec:discussion}

Here we mention several interesting questions for future research that are inspired by our results.

\begin{itemize}
    \item \textbf{Are warm starts essential for future progress in high-accuracy sampling?} Our work is the first to show the achievability of the faster rates proven for high-accuracy samplers under a warm start assumption. We do this by exhibiting an efficient algorithm for producing the warm start.
    We believe that this general strategy may be important for future progress in high-accuracy sampling.
    Indeed, the natural next candidate for improving upon MALA is Metropolized Hamiltonian Monte Carlo~\cite{neal2011mcmc}, or related variants.
    For such Metropolized algorithms, we suspect that much of the intuition from \S\ref{ssec:intro:tech} remains true; namely, that the algorithm can benefit from a more aggressive step size near stationarity.
    Hence, to extract the full potential of these algorithms, it seems likely that we must again pursue the dual plan of improving the rates under a warm start, and efficiently computing that warm start.
    Insofar as warm starts continue to play an important role in sampling analysis, the R\'enyi analysis techniques that we developed in \S\ref{sec:pabi} and \S\ref{sec:ulmc} could prove useful for future progress in this direction.
    \item \textbf{Can we leverage shifted divergence techniques for further advances in differential privacy and beyond?} Core to our results is an improved version of the shifted R\'enyi divergence technique that uses Orlicz--Wasserstein shifts rather than $W_{\infty}$ shifts.
    Since their introduction, shifted divergences have been instrumental for advances in differentially private optimization (see the prior work discussion in \S\ref{ssec:intro:prior}), and also very recently in the context of sampling (\cite{AltTal22mix} and Theorem~\ref{thm:loacc}).
    We believe that we are only scratching the surface of potential applications, extensions, and refinements of this technique, and we are optimistic that a deeper understanding of our R\'enyi analysis toolbox will have implications far beyond.
    %Since the general technique of shifted R\'enyi divergences has been key to many recent developments in the field of differentially private optimization, it is natural to ask if this improved version of the technique can lead to further progress in differential privacy.     
\end{itemize}

%% file: sections/app_background.tex
\section{Background on functional inequalities}\label{ssec:prelim-functional}

In this section, we collect relevant background material for the convenience of the reader.

To go beyond the strongly log-concave case, we can instead assume that the target distribution $\pi$ satisfies a \emph{functional inequality}, which encodes geometric information (e.g., isoperimetric properties) about $\pi$.
We focus primarily on the two most well-studied functional inequalities in this context, namely the \emph{log-Sobolev inequality} and the \emph{Poincar\'e inequality}.
For each of these two functional inequalities, the class of distributions satisfying this assumption not only includes all strongly log-concave distributions
 (see Lemma~\ref{lem:bakry-emery} below), but also includes many more examples because of closedness properties of these functional inequalities under operations such as bounded perturbations of the potential, pushforwards via Lipschitz mappings, or taking suitable mixtures.
 We refer to~\cite{bakry2014analysis} for many of these properties, and to~\cite{chenchewinilesweed2021dimfreelsi} for the closure under taking mixtures.

 \begin{defin}
    A probability distribution $\pi$ on $\R^d$ satisfies a \emph{log-Sobolev inequality} (LSI) with constant $C_{\msf{LSI}}$ if for all smooth and compactly supported functions $\phi : \R^d\to\R$,
    \begin{align*}
        \on{ent}_\pi(\phi^2)
        &\deq \E_\pi\Bigl[ \phi^2 \log \frac{\phi^2}{\E_\pi[\phi^2]}\Bigr]
        \le 2C_{\msf{LSI}} \,\E_\pi[\norm{\nabla \phi}^2]\,.
    \end{align*}
 \end{defin}

 \begin{defin}
     A probability distribution $\pi$ on $\R^d$ satisfies a \emph{Poincar\'e inequality} (PI) with constant $C_{\msf{PI}}$ if for all smooth and compactly supported functions $\phi : \R^d\to\R$,
     \begin{align*}
         \var_\pi(\phi)
         &\le C_{\msf{PI}} \,\E_\pi[\norm{\nabla \phi}^2]\,.
     \end{align*}
 \end{defin}

We next collect together key facts about these functional inequalities.
The following results show that the class of distributions satisfying these inequalities is larger than the class of strongly log-concave distributions; see~\cite[Proposition 5.1.3 and Corollary 5.7.2]{bakry2014analysis}.

\begin{lemma}[Strong log-concavity implies LSI implies PI]\label{lem:bakry-emery}
    Let $\pi$ be a distribution on $\R^d$.
    \begin{enumerate}
        \item (\emph{Bakry--\'Emery theorem}) If $\pi$ is $\alpha$-strongly log-concave, then it satisfies an LSI with constant at most $1/\alpha$.
        \item If $\pi$ satisfies an LSI with constant $C_{\msf{LSI}}$, then it also satisfies a PI with constant at most $C_{\msf{LSI}}$.
    \end{enumerate}
\end{lemma}

A useful consequence of the LSI is the following sub-Gaussian concentration inequality, typically established via the Herbst argument; see~\cite[Proposition 5.4.1]{bakry2014analysis}.

\begin{lemma}[LSI implies sub-Gaussian concentration]
    Suppose that $\pi$ is a distribution on $\R^d$ satisfying an LSI with constant $C_{\msf{LSI}}$.
    Then, for any $L$-Lipschitz function $\phi : \R^d\to\R$ and any $\lambda \in \R$,
    \begin{align*}
        \E_\pi \exp\bigl(\lambda\,(\phi - \E_\pi \phi)\bigr)
        &\le \exp\Bigl( \frac{\lambda^2 C_{\msf{LSI}} L^2}{2}\Bigr)\,.
    \end{align*}
    Consequently, for all $\eta \ge 0$,
    \begin{align*}
        \pi\{\phi - \E_\pi \phi \ge \eta\}
        &\le \exp\Bigl( -\frac{\eta^2}{2C_{\msf{LSI}} L^2}\Bigr)\,.
    \end{align*}
\end{lemma}

Similarly, the PI implies subexponential concentration, see~\cite[\S 4.4.3]{bakry2014analysis}.

\begin{lemma}[PI implies subexponential concentration]\label{lem:pi-subexp}
    Suppose that $\pi$ is a distribution on $\R^d$ satisfying a PI with constant $C_{\msf{PI}}$.
    Then, for any $L$-Lipschitz function $\phi : \R^d\to\R$ and any $\eta \ge 0$,
    \begin{align*}
        \pi\{\phi - \E_\pi \phi \ge \eta\}
        &\le 3\exp\Bigl( -\frac{\eta}{\sqrt{C_{\msf{PI}}}\, L}\Bigr)\,.
    \end{align*}
\end{lemma}

Next we recall two comparison inequalities which enable proving sampling guarantees in Wasserstein distance as an immediate corollary of proving sampling guarantees in other metrics---namely KL divergence in the LSI setting, and chi-squared divergence in the PI setting. Such comparison inequalities are often called transport inequalities. Specifically, the first result, attributed to Otto and Villani~\cite{ottovillani2000lsi}, shows that under an LSI, a transportation inequality between Wasserstein and KL divergence holds (this inequality is often referred to as Talagrand's $T_2$ inequality).

\begin{lemma}[Otto--Villani theorem]\label{lem:otto-villani}
    Suppose that $\pi$ is a distribution on $\R^d$ satisfying an LSI with constant $C_{\msf{LSI}}$. Then, for all distributions $\mu \in \cP_2(\R^d)$,
    \begin{align*}
        W_2^2(\mu, \pi)
        &\le 2C_{\msf{LSI}}\, \KL(\mu \mmid \pi)\,.
    \end{align*}
\end{lemma}

The second result shows a similar transport inequality in the PI setting~\cite{liu2020poincare}. Under a PI, Talagrand's $T_2$ inequality does not necessarily hold anymore. Nevertheless, a useful transport inequality still holds if one replaces the KL divergence by the chi-squared divergence.

\begin{lemma}[Quadratic transport-variance inequality]\label{lem:quadratic-transport-var}
    Suppose that $\pi$ is a distribution on $\R^d$ satisfying a PI with constant $C_{\msf{PI}}$.
    Then, for all distributions $\mu \in \cP_2(\R^d)$,
    \begin{align*}
        W_2^2(\mu,\pi)
        &\le 2C_{\msf{PI}}\,\chi^2(\mu \mmid \pi)\,.
    \end{align*}
\end{lemma}

Finally, we record the following standard second-moment-type bound for strongly log-concave measures; see, e.g.,~\cite[Proposition 2]{dalalyan2022bounding}. We give a short proof sketch for the convenience of the reader.

\begin{lemma}[Second-moment-type bound under strong log-concavity]\label{lem:lsi-second-moment}
	Suppose that the distribution $\pi \propto \exp(-f)$ is $\alpha$-strongly log-concave, and that $x^*$ is the minimizer of $f$.
 Then,
	\begin{align*}
		\E_{X \sim \pi}[ \|X - x^*\|^2] \leq d/\alpha\,.
	\end{align*}
\end{lemma}
\begin{proof}
    Integration by parts shows that for any smooth function $\phi : \R^d\to\R$ of controlled growth, it holds that $\E_\pi[\Delta \phi - \langle \nabla f,\nabla \phi\rangle] = 0$.
    We apply this to $\phi(x) \deq\half\, \norm{x-x^*}^2$, for which $\nabla \phi(x) = x-x^*$ and $\Delta \phi = d$.
    By strong convexity of $f$, $\langle \nabla f(x), x-x^*\rangle \ge \alpha\,\norm{x-x^*}^2$, and the result follows.
\end{proof}

%% file: sections/app_pabi.tex
\section{Deferred details for \S\ref{sec:pabi}}\label{app:pabi}

\subsection{Proof for Remark~\ref{rem:wass-comparison}}\label{app:wass-comparison}

Here, we prove the inequality in Remark~\ref{rem:wass-comparison}, repeated here for convenience:
\begin{align}
	\frac{1}{\sqrt p}\, W_p \lesssim W_{\psi_2} \leq \frac{1}{\sqrt{\log 2}}\, W_{\infty}\,.
	\label{eq:wass-comparison}
\end{align}

\paragraph*{Bounding $W_{p}$ by $W_{\psi_2}$.}
By~\cite[Proposition 2.5.2]{vershynin2018high}, if $(X, Y)$ is an optimal coupling of $\mu$ and $\nu$ for the $W_{\psi_2}$ distance, then
\[
W_p(\mu,\nu)
\le {\E[\norm{X-Y}^p]}^{1/p}
\lesssim \sqrt p\, \norm{X-Y}_{\psi_2}
= W_{\psi_2}(\mu,\nu)\,.
\]

\paragraph*{Bounding $W_{\psi_2}$ by $W_{\infty}$.} Observe that for any random variable $Z$, if we denote $\|Z\|_{\infty}  \deq  \esssup \|Z\|$, then we can bound
\[
\E\biggl[ \psi_2\Bigl( \frac{Z}{\|Z\|_{\infty} / \sqrt{\log 2}}\Bigr) \biggr]
=
\E\biggl[ \exp\Bigl( \frac{\|Z\|^2}{\|Z\|^2_{\infty} } \cdot \log 2\Bigr) - 1\biggr]
\leq
\exp(\log 2) - 1
= 1\,,
\]
and therefore $\|Z\|_{\psi_2} \leq \|Z\|_{\infty} / \sqrt{\log 2}$ by the definition of the Orlicz norm. Now, applying this bound to the random variable $Z = X-Y$, we conclude the desired inequality
\[
W_{\psi_2}(\mu,\nu)
= \inf_{(X,Y) \in \cC(\mu,\nu)} \|X-Y\|_{\psi_2}
\leq \frac{1}{\sqrt{\log 2}}\, \inf_{(X,Y) \in \cC(\mu,\nu)} \|X-Y\|_{\infty}
= \frac{1}{\sqrt{\log 2}}\, W_{\infty}(\mu,\nu)\,.
\]

\subsection{Remarks on tightness of Lemma~\ref{lem:shift-reduction}}\label{app:shift-reduction-tight}

Here we remark that, conditional on the following plausible conjecture, the generalized version of Lemma~\ref{lem:shift-reduction} (as stated in Remark~\ref{rem:shift-reduction-holder}) is tight. This conjecture states that the shifted R\'enyi divergence between two isotropic Gaussians with same covariance is achieved by a deterministic shift. Understanding this simple case could be more broadly helpful for understanding tightness of other inequalities and analyses using the shifted R\'enyi divergence.

\begin{conj}\label{conj:renyi-shift-gaussians}
	For any R\'enyi order $q \geq 1$, noise variance $\sig^2 > 0$, shift $w \geq 0$, and mean $x \in \R^d$,
	\[
	\cR_q^{(w)} \bigl( \cN(x, \sig^2 I_d) \bigm\Vert \cN(0, \sig^2 I_d)\bigr)
	=
	\cR_q \bigl( \cN(cx, \sig^2 I_d) \bigm\Vert \cN(0, \sig^2 I_d)\bigr)
	=
	\frac{c^2 q\,\|x\|^2}{2\sig^2}\,,
	\]
	where $c  \deq  \max(0, 1 - w\sqrt{\log 2}/\|x\|)$.
\end{conj}

Of course, the conjecture here is the first equality (the second equality is just the closed-form expression for the R\'enyi divergence between Gaussians in Lemma~\ref{lem:renyi-gaussians}). The direction ``$\leq$'' is clear because $\cN(cx, \sig^2 I_d)$ satisfies
	\[
	W_{\psi_2} \bigl( \cN(x, \sig^2 I_d), \, \cN(cx, \sig^2 I_d) \bigr) \leq \|x - cx\|_{\psi_2} = \frac{1 - c}{\sqrt{\log 2}} \,\|x\| \leq w\,,
	\]
and therefore is feasible for the optimization problem defining the shifted R\'enyi divergence. The direction ``$\geq$'' is the one requiring justification.

\par In the rest of this subsection, we show the claimed tightness assuming Conjecture~\ref{conj:renyi-shift-gaussians}. Fix any R\'enyi order $q \geq 1$, noise variance $\sig^2 > 0$, and initial shift $w \geq 0$. Consider distributions $\mu = \delta_{a e_1}$ and $\nu = \delta_0$, where $a > w \sqrt{\log 2}$. We claim that the bound in Remark~\ref{rem:shift-reduction-holder} holds with equality when its parameters $\delta, \lambda$ are optimized; that is,
\begin{align}
	\cR_q^{(w)}\bigl(\mu \ast \cN(0,\sig^2 I_d) \bigm\vert \nu \ast \cN(0, \sig^2 I_d)\bigr)
	&= \inf_{\delta > 0,\, \lambda \in [0,1]}{\Bigl[ \cR_{(q+\lambda-1)/\lambda}^{(w+\delta)}(\mu \mmid \nu) + 
	\frac{(q-\lambda)\,\delta^2}{2\,(1-\lambda)\, \sigma^2} \log 2 \Bigr]}
	\,.
	\label{eq-app:shift-reduction-tight}
\end{align}
\par To this end, supposing Conjecture~\ref{conj:renyi-shift-gaussians} holds, the left hand side of~\eqref{eq-app:shift-reduction-tight} is equal to
\begin{align*}
	\cR_q^{(w)}\bigl(\cN(ae_1,\sig^2 I_d) \bigm\Vert \cN(0, \sig^2 I_d)\bigr)
	&=  	\cR_q\bigl(\cN((a-w\sqrt{\log 2})\, e_1,\sig^2 I_d) \bigm\Vert \cN(0, \sig^2 I_d)\bigr) \\
	&= \frac{q\,(a - w\sqrt{\log 2})^2}{2\sig^2}\,.
\end{align*}
\par On the other hand, note that $\cR_{(q+\lambda-1)/\lambda}^{(w+\delta)}(\mu \mmid \nu)$ is equal to $0$ if $a \leq (w + \delta)\sqrt{\log 2}$, and otherwise is equal to $\infty$. This means that the optimal value of $\delta$ is $a/\sqrt{\log 2} - w$. Thus the right hand side of~\eqref{eq-app:shift-reduction-tight} simplifies to
\begin{align*}
	\inf_{\lambda \in [0,1]} 
	\frac{(q-\lambda)\,(a - w\sqrt{\log 2})^2}{2\,(1-\lambda)\, \sigma^2}
	=
	\frac{q\,(a - w\sqrt{\log 2})^2}{2\sig^2}\,,
\end{align*}
where the final step is because the optimal value of $\lambda$ is at $\lambda = 0$. We conclude that the left- and right-hand sides of~\eqref{eq-app:shift-reduction-tight} indeed match, as desired.

\subsection{Proof for Remark~\ref{rmk:weaker-shifts}}\label{app:rmk-weaker-shifts}

Here we provide details for why Lemma~\ref{lem:shift-reduction} fails if the shifted R\'enyi divergence is defined using the $W_p$ metric for any $p$ finite, or alternatively with the Orlicz--Wasserstein metric with any Orlicz norm that is weaker than sub-Gaussian (i.e., with Orlicz function $\psi_b(x) \deq \exp(x^b) - 1$ for $b < 2$). Specifically, for any of these Wasserstein metrics $W$, we claim that if $\mu$ is the distribution on $\R$ with density proportional to $\exp(-|\cdot|^a)$ for an appropriate constant $a < 2$, and $\nu = \delta_0$, then:
\begin{itemize}
    \item [(i)] $W(\mu,\nu) < \infty$.
    \item [(ii)] $\cR_q(\mu \ast \cN(0,\sig^2) \mmid \nu \ast \cN(0,\sig^2)) = \infty$.
\end{itemize}
This comprises a counterexample to Lemma~\ref{lem:shift-reduction} by taking $w = 0$ and $\delta = W(\mu,\nu)$ and noting that if $\delta > \sig / \sqrt{(2q-1)(q-1)}$, then we can simply dilate the space (i.e., replace $\mu(x)$ by $\mu_b(x) \propto \mu(bx)$ for a sufficiently large $b > 0$) and repeat the same argument. 

\paragraph*{Proof of (i).} In the case that $W = W_p$, then $W_p(\mu,\nu)$ is equal to the $p$-th norm of $\mu$, which is finite for any $p < \infty$. In the case that $W$ is an Orlicz--Wasserstein norm with Orlicz norm weaker than sub-Gaussian, then $W(\mu,\nu)$ is equal to the Orlicz norm of $\nu$, which is finite if we choose $a=b$.

\paragraph*{Proof of (ii).} Note that $\mu \ast \cN(0,\sig^2)$ is not sub-Gaussian, yet $\nu \ast \cN(0,\sig^2) = \cN(0,\sig^2)$ is sub-Gaussian. We may therefore appeal to the fact that the R\'enyi divergence is infinite whenever the first argument is not sub-Gaussian, but the second argument is. For a proof of this fact in the case that $q=2$, see~\cite[Lemma 21]{chewi2021analysis}; this proof readily extends to any finite $q \in (1,\infty)$ by replacing the Cauchy--Schwarz inequality by H\"older's inequality. 

%% file: sections/app_loacc.tex
\section{Deferred details for \S\ref{sec:ulmc}}\label{app:loacc}

\subsection{Proof of Lemma~\ref{lem:ulmc-contract}}\label{app:pf-lem:ulmc-contract}

	Compute the partial derivatives
	\begin{align*}
		\partial_u {\bar F(u,v)}_u
		&= \frac{1+a}{2}\, I_d - \frac{h-\gamma^{-1}\,(1-a)}{\gamma}\,\nabla^2 f(u)\,, \\
		\partial_u {\bar F(u,v)}_v
		&= \frac{1-a}{2}\, I_d - \frac{h+\gamma^{-1}\,(1-a)}{\gamma}\,\nabla^2 f(u)\,, \\
		\partial_v {\bar F(u,v)}_u
		&= \frac{1-a}{2}\, I_d\,, \\
		\partial_v {\bar F(u,v)}_v
		&= \frac{1+a}{2} \, I_d\,.
	\end{align*}
	Since $\tfrac{1}{\gamma} (h-\gamma^{-1}\,(1-a)) = O(h^2)$, we have
	\begin{align*}
		\norm{\nabla \bar F(u,v)}_{\rm op}
		&\le \frac{1}{2}\, \Bigl\lVert \underbrace{\begin{bmatrix} (1+a)\,I_d & (1-a)\,I_d - b\,\nabla^2 f(u) \\ (1-a)\, I_d & (1+a)\, I_d \end{bmatrix}}_{\eqqcolon A}\Bigr\rVert_{\rm op} + O(\beta h^2) \,,
	\end{align*}
	where we use the notational shorthand $b \deq \frac{2}{\gamma}\,(h+\gamma^{-1}\,(1-a))$.

	\par To bound this operator norm, we compute
	\begin{align*}
		AA^\T
		&= \begin{bmatrix} {(1+a)}^2\,I_d + {((1-a)\,I_d-b\,\nabla^2 f(u))}^2 & 2\,(1-a^2)\,I_d - (1+a)\,b\,\nabla^2 f(u) \\ 2\,(1-a^2)\,I_d - (1+a)\,b\,\nabla^2 f(u) & ((1-a)^2 + (1+a)^2)\,I_d\end{bmatrix}\,.
	\end{align*}
	Since $1-a = O(\gamma h)$ and $b = O(h/\gamma)$, we can approximate $AA^\T$ by the following matrix $B$ with error
	\begin{align*}
		&\Bigl\lVert AA^\T - 2\,\underbrace{\begin{bmatrix} (1+a^2)\,I_d & (1-a^2)\,I_d - b\,\nabla^2 f(u) \\ (1-a^2)\, I_d - b\,\nabla^2 f(u) & (1+a^2)\, I_d \end{bmatrix}}_{\eqqcolon B}\Bigr\rVert_{\rm op}
		= O\Bigl( \frac{\beta^2 h^2}{\gamma^2} + \beta h^2 \Bigr)\,.
	\end{align*}
	By a direct computation, the eigenvalues of $B$ are $1+a^2 \pm (1-a^2 - b\lambda)$, where $\lambda$ ranges over the eigenvalues of $\nabla^2 f(u)$. The strong-convexity and smoothness of $f$ implies that $\lambda \in [\alpha,\beta]$. Thus
	\begin{align*}
		\norm B_{\rm op}
		&\le \max\{2a^2 + \beta b, 2-\alpha b\}\,.
	\end{align*}
	We note that
	\begin{align*}
		2a^2 + \beta b
		= 2\exp(-2\gamma h) + \frac{2 \beta \,(h+\gamma^{-1}\,(1-\exp(-\gamma h)))}{\gamma}
		= 2\,\Bigl(1 - 2\gamma h + \frac{2\beta h}{\gamma} + O(\gamma^2 h^2 + \beta h^2)\Big)\,.
	\end{align*}
	In order for this to be strictly smaller than $2$, we must take $\gamma > \sqrt \beta$.
	We choose $\gamma = \sqrt{2\beta}$, whereby
	\begin{align*}
		\norm B_{\rm op}
		&\le 2 \max\Bigl\{1- h\sqrt{2 \beta},\; 1- \alpha h \sqrt{2/\beta}\Bigr\} + O(\beta h^2) 
		= 2\,\Bigl(1-\alpha h \sqrt{2/\beta} \Bigr) + O(\beta h^2)\,.
	\end{align*}
	We deduce that
	\begin{align*}
		\norm{AA^\T}_{\rm op}
		&\le 4\,\Bigl(1-\alpha h \sqrt{2/\beta} \Bigr) + O(\beta h^2)
	\end{align*}
	and therefore
	\begin{align*}
		\norm{\nabla \bar F(u,v)}_{\rm op}
		&\le \sqrt{1-\alpha h \sqrt{\frac{2}{\beta}}} + O(\beta h^2)
		\le 1-\frac{\alpha }{\sqrt{2\beta}}\,h + O(\beta h^2)\,. \qedhere
	\end{align*}

\subsection{Proof of Lemma~\ref{lem:ulmc-noise}}\label{app:pf-lem:ulmc-noise}

By definition of $\cM$ and $\Sigma$, 
\[
\bar{\Sigma}
= \cM \Sigma \cM^T
=
\frac{1}{\gamma^2} \,\begin{bmatrix}
	2\gamma h + 4 a - a^2 - 3 & 2\gamma h + a^2 - 1 \\
	2\gamma h + a^2 - 1 & 2\gamma h + 5 - a^2 - 4a 
\end{bmatrix}\,.
\]
The smallest eigenvalue of this matrix is
\begin{align*}
	\lambda_{\min}(\bar\Sigma) 
	&=
	\half \,\Bigl( \Tr(\bar \Sigma) - \sqrt{\Tr(\bar{\Sigma})^2 - 4 \det (\bar \Sigma)} \Bigr)
	\\ &=
	\frac{1}{\gamma^2}\, \Bigl( 2 \gamma h + 1 - a^2 - \sqrt{17 - 32a + 14a^2 + a^4 - 4\gamma h + 4 a^2 \gamma h+ 4\gamma^2 h^2}
	\Bigr)
	\\ &= \frac{\gamma h^3}{6}\, \Bigl(1 - \frac{\gamma h}{2} + \frac{\gamma^2 h^2}{240} - \cdots \Bigr)\,.
\end{align*}
Above, the first step is by the explicit formula for the eigenvalues of a $2 \times 2$ matrix; the second step is by plugging in the entries of $\bar \Sigma$ and simplifying; and the third step is by performing a Taylor expansion in the variable $\gamma h$.

\subsection{Proof of Lemma~\ref{lem:ulmc-bias}}\label{app:loacc-girsanov}

We invoke the following result, which appears as Lemma 26 in~\cite{Matt23girsanov}.

\begin{prop}[Movement bound for underdamped Langevin]\label{prop:ulmc-movement}
    Let ${(X_t, Y_t)}_{t\ge 0}$ denote the continuous-time underdamped Langevin diffusion~\eqref{eq:underdamped-diffusion} with potential $f$ that is $\beta$-smooth and minimized at $x^*$.
    Assume that $0 < h \lesssim \frac{1}{\sqrt\beta \vee \gamma}$ and $0 \le \lambda \lesssim \frac{1}{\gamma dh^3}$. Then, conditioned on $(X_0, Y_0)$,
    \begin{align*}
        \log \E \exp\bigl(\lambda \sup_{t \in [0,h]}{\norm{X_t - X_0}^2}\bigr)
        &\lesssim \bigl(\beta^2 h^4 \,\norm{X_0 - x^\star}^2 + h^2 \, \norm{Y_0}^2 + \gamma dh^3\bigr)\,\lambda\,.
    \end{align*}
\end{prop}
\begin{proof}
    This result can be easily adapted from the proof of~\cite[Lemma 26]{Matt23girsanov}, noting that in our situation the bound simplifies as we are assuming $\nabla f$ is $\beta$-Lipschitz rather than merely H\"older continuous.
\end{proof}

Let $\patha_T$, $\pathb_T$ denote the path measures (i.e., probability measures over $\mc C([0,T]; \R^d)$) for the discretized and continuous underdamped Langevin processes respectively, both started at the stationary measure $\bs \pi$.
Girsanov's theorem~\cite[Theorem 5.22]{legall2016stochasticcalc} yields
\begin{align*}
    \frac{\D\patha_T}{\D \pathb_T}
    &= \exp \sum_{k=0}^{N-1}\Bigl( \frac{1}{\sqrt{2\gamma}}\int_{kh}^{(k+1)h} \langle \nabla f(X_t) - \nabla f(X_{kh}), \D B_t \rangle - \frac{1}{4\gamma} \int_{kh}^{(k+1)h} \norm{\nabla f(X_t) - \nabla f(X_{kh})}^2 \, \D t\Bigr) \\
    &\eqqcolon \exp M_T \,,
\end{align*}
provided that Novikov's condition (see~\cite[Theorem 5.23]{legall2016stochasticcalc}) holds:
\begin{align}\label{eq:novikov}
    \E_{\pathb_T} \exp\sum_{k=0}^{N-1} \frac{1}{4\gamma} \int_{kh}^{(k+1)h} \norm{\nabla f(X_t) - \nabla f(X_{kh})}^2 \, \D t < \infty\,.
\end{align}
Assuming for the moment that~\eqref{eq:novikov} is indeed verified, It\^o's formula yields
\begin{align*}
    \E_{\pathb_T}\bigl[ \bigl( \frac{\D \patha_T}{\D \pathb_T} \bigr)^q\bigr] - 1
    &= \frac{q\,(q-1)}{4\gamma}\, \E_{\pathb_T}\sum_{k=0}^{N-1} \int_{kh}^{(k+1)h}\exp(qM_t)\, \norm{\nabla f(X_t) -\nabla f(X_{kh})}^2 \, \D t \\
    &\le \frac{q^2}{4\gamma}\sum_{k=0}^{N-1}\int_{kh}^{(k+1)h} \sqrt{\E_{\pathb_T}\exp(2qM_t)\, \E_{\pathb_T}[\norm{\nabla f(X_t) - \nabla f(X_{kh})}^4]} \, \D t\,.
\end{align*}

\paragraph*{Bounding the first term.}
For $t \in [kh,(k+1)h]$, let $\Delta_t \deq \nabla f(X_t) - \nabla f(X_{kh})$.
By the Cauchy{--}Schwarz inequality,
\begin{align*}
    \E_{\pathb_T}\exp(2qM_t)
    &\le \underbrace{\sqrt{\E_{\pathb_T} \exp \sum_{k=0}^{N-1} \int_{kh}^{(k+1)h \wedge t} \Bigl( \frac{2\sqrt 2\,q}{\sqrt\gamma}\, \langle \Delta_s, \D B_s\rangle - \frac{4q^2}{\gamma} \,\norm{\Delta_s}^2 \, \D s \Bigr)}}_{(\dagger)} \\
    &\qquad{}\times \sqrt{\E_{\pathb_T} \exp \sum_{k=0}^{N-1} \int_{kh}^{(k+1)h\wedge t} \bigl( \frac{4q^2}{\gamma} - \frac{q}{\gamma} \bigr)\,\norm{\Delta_s}^2\,\D s}\,.
\end{align*}
We claim that the term marked $(\dagger)$ equals $1$; this would follow if the quantity inside the expectation is a martingale.
In general, it is only a local martingale, but it is a bona fide martingale provided that Novikov's condition holds: it suffices to have
\begin{align}\label{eq:novikov-2}
    \E_{\pathb_T} \exp \sum_{k=0}^{N-1} \frac{4q^2}{\gamma} \int_{kh}^{(k+1)h} \norm{\nabla f(X_t) - \nabla f(X_{kh})}^2 \, \D t < \infty\,.
\end{align}
Note that this condition is stronger than~\eqref{eq:novikov}.

Towards this end, we bound, for a parameter $p\ge 1$ to be chosen later,\footnote{This argument avoids the use of the ``conditioning lemma'' from~\cite{GaneshT20} (Lemma 23 in~\cite{Matt23girsanov}).}
\begin{align*}
    &\E_{\pathb_T} \exp \sum_{k=0}^{N-1} \frac{4q^2}{\gamma} \int_{kh}^{(k+1)h} \norm{\nabla f(X_t) - \nabla f(X_{kh})}^2 \, \D t
    \le \E_{\pathb_T} \exp \sum_{k=0}^{N-1} \frac{4\beta^2 q^2}{\gamma} \int_{kh}^{(k+1)h} \norm{X_t - X_{kh}}^2 \, \D t \\
    &\qquad \le \E_{\pathb_T} \exp \max_{k=0,1,\dotsc,N-1} \sup_{t\in [kh,(k+1)h]} \frac{4\beta^2 q^2 T}{\gamma} \, \norm{X_t - X_{kh}}^2 \\
    &\qquad \le \biggl\{\E_{\pathb_T} \exp \max_{k=0,1,\dotsc,N-1} \sup_{t\in [kh,(k+1)h]} \frac{4\beta^2 pq^2 T}{\gamma} \, \norm{X_t - X_{kh}}^2\biggr\}^{1/p} \\
    &\qquad \le N^{1/p}\, \biggl\{\max_{k=0,1,\dotsc,N-1}\E_{\pathb_T} \exp \sup_{t\in [kh,(k+1)h]} \frac{4\beta^2 pq^2 T}{\gamma} \, \norm{X_t - X_{kh}}^2\biggr\}^{1/p}\,.
\end{align*}
By conditioning on $(X_{kh}, Y_{kh})$ and applying Proposition~\ref{prop:ulmc-movement}, this is bounded by
\begin{align*}
    \cdots
    &\le N^{1/p}\,\biggl\{\max_{k=0,1,\dotsc,N-1} \E_{\pathb_T}\exp\biggl( O\Bigl( \frac{\beta^2 pq^2 T}{\gamma}\,(\beta^2 h^4 \,\norm{X_{kh} - x^*}^2 + h^2\,\norm{Y_{kh}}^2 + \gamma dh^3) \Bigr)\biggr) \biggr\}^{1/p}
\end{align*}
provided that $h \lesssim \frac{1}{\beta^{2/3}d^{1/3} p^{1/3} q^{2/3} T^{1/3}}$.
We choose $p\asymp \log N$ so that $N^{1/p} \asymp 1$. We now need tail bounds for $\norm{X_{kh} - x^*}$ and $\norm{Y_{kh}}$.

Using the argument in the proof of Lemma~\ref{lem:ulmc-init}, for $c > 0$,
\begin{align*}
    \E_{\pathb_T} \exp(c\,\norm{X_{kh}-x^*}^2)
    &\le \exp(2c\,\E_{\pathb_T}[\norm{X_{kh} - x^*}^2])\,\E_{\pathb_T} \exp\bigl(2c\, {(\norm{X_{kh}-x^*} - \E_{\pathb_T}\norm{X_{kh}-x^*})}^2\bigr) \\
    &\le \exp\bigl( \frac{2cd}{\alpha}\bigr) \, \Bigl( \E_{\pathb_T} \exp\Bigl( \frac{2\, {(\norm{X_{kh}-x^*} - \E_{\pathb_T}\norm{X_{kh}-x^*})}^2}{36/\alpha}\Bigr) \Bigr)^{36c/\alpha} \\
    &\le \exp\Bigl( O\bigl( \frac{cd}{\alpha} \bigr)\Bigr)
\end{align*}
provided that $c\le \alpha/36$.
Therefore,
\begin{align*}
    \E_{\pathb_T} \exp\Bigl( O\Bigl( \frac{\beta^4 h^4 q^2 T\log N}{\gamma} \, \norm{X_{kh} - x^*}^2 \Bigr)\Bigr)
    &\le \exp\Bigl( O\Bigl( \frac{\beta^4 dh^4 q^2 T\log N}{\alpha\gamma}\Bigr)\Bigr)
\end{align*}
provided that $h \lesssim \frac{\alpha^{1/4} \gamma^{1/4}}{\beta q^{1/2} \, {(T \log N)}^{1/4}}$.

The same argument applied to $\norm{Y_{kh}}$ yields
\begin{align*}
    \E_{\pathb_T} \exp\Bigl( O\Bigl( \frac{\beta^2 h^2 q^2 T\log N}{\gamma} \,\norm{Y_{kh}}^2\Bigr) \Bigr)
    \le \exp\Bigl( O\Bigl( \frac{\beta^2 dh^2 q^2 T \log N}{\gamma} \Bigr)\Bigr)
\end{align*}
provided that $h\lesssim \frac{\gamma^{1/2}}{\beta q \, {(T\log N)}^{1/2}}$.

If we put these bounds together and take $\gamma \asymp \sqrt\beta$, we deduce that if $h\lesssim \frac{\gamma^{1/2}}{\beta d^{1/3} q \, {(T\log N)}^{1/2}}$ and $T \gtrsim \frac{\sqrt \beta}{\alpha}$, then it holds that
\begin{align*}
    \E_{\pathb_T} \exp \sum_{k=0}^{N-1} \frac{4q^2}{\gamma} \int_{kh}^{(k+1)h} \norm{\nabla f(X_t) - \nabla f(X_{kh})}^2 \, \D t
    &\le \exp\bigl( O(\beta^{3/2} dh^2 q^2 T\log N)\bigr)\,.
\end{align*}
This verifies~\eqref{eq:novikov} and~\eqref{eq:novikov-2}, and moreover shows that for $h\lesssim \frac{1}{\beta^{3/4} d^{1/2} q\,{(T\log N)}^{1/2}}$,
\begin{align*}
    \sup_{t\in [0, T]} \E_{\pathb_T}\exp(2qM_t) 
    \lesssim 1\,.
\end{align*}

\paragraph*{Bounding the second term.}
Next,
\begin{align*}
    \sqrt{\E_{\pathb_T}[\norm{\nabla f(X_t) - \nabla f(X_{kh})}^4]}
    &\le \beta^2 \sqrt{\E_{\pathb_T}[\norm{X_t - X_{kh}}^4]}\,.
\end{align*}
Applying Proposition~\ref{prop:ulmc-movement} and the above tail estimates,
\begin{align*}
    \cdots
    &\lesssim \beta^2 \,\E_{\pathb_T}[\beta^2 h^4 \, \norm{X_{kh} - x^*}^2 + h^2 \, \norm{Y_{kh}}^2 + \gamma dh^3]
    \lesssim \frac{\beta^4 dh^4}{\alpha} + \beta^2 dh^2 + \beta^2 \gamma dh^3
    \lesssim \beta^2 dh^2\,.
\end{align*}

\paragraph*{Concluding the proof of Lemma~\ref{lem:ulmc-bias}.}
In summary, we have shown
\begin{align*}
    \E_{\pathb_T}\bigl[ \bigl( \frac{\D \patha_T}{\D \pathb_T} \bigr)^q\bigr] - 1
    &\lesssim \frac{\beta^2 dh^2 q^2 T}{\gamma}\,.
\end{align*}
Hence, by definition of the R\'enyi divergence, we conclude that
\begin{align*}
    \cR_q(\patha_T \mmid \pathb_T)
    &\lesssim \frac{\beta^2 dh^2 q T}{\gamma}\,.
\end{align*}

%% file: sections/app_hiacc.tex
\section{Deferred details for \S\ref{sec:hiacc}}\label{app:hiacc}

\subsection{Background on the Metropolis-Adjusted Langevin Algorithm}\label{app:hiacc:mala}

Let $\pi\propto\exp(-f)$ be a target probability distribution on $\R^d$, and let $Q$ be a \emph{proposal kernel} (that is, a Markov transition kernel). Then, the Metropolis-adjusted algorithm with proposal $Q$ is the following Markov chain ${(X_n)}_{n\in\N}$: for $n=0,1,2,\dotsc$,
\begin{enumerate}
    \item \textbf{Propose} a new point $Y_n \sim Q(X_n, \cdot)$;
    \item \textbf{Accept} the new point with probability $1 \wedge \frac{\pi(Y_n)\,Q(Y_n, X_n)}{\pi(X_n)\, Q(X_n, Y_n)}$.
    That is, if the proposal is accepted, we set $X_{n+1} \deq Y_n$, otherwise we set $X_{n+1} \deq X_n$.
\end{enumerate}
It can be shown that under mild conditions on the proposal $Q$, the resulting Metropolis-adjusted algorithm is a reversible Markov chain with stationary distribution $\pi$.
Since Markov chains often exhibit geometrically ergodicity, the family of Metropolis-adjusted chains is often used to design high-accuracy samplers.

In this work, we are primarily considered with the Metropolis-adjusted Langevin algorithm, in which the proposal kernel is taken to be one step of the discretized Langevin algorithm, i.e.,
\begin{align*}
    Q(x,\cdot)
    &\deq \mc N\bigl(x-h\,\nabla f(x), \, 2hI_d\bigr)\,.
\end{align*}
For this choice of proposal kernel, generating a proposal and computing the acceptance ratio can be carried out with a constant number of evaluations of $f$ (zeroth-order queries) and $\nabla f$ (first-order queries).
Actually, in order for the mixing time results we invoke to be valid, we must instead consider the $\frac{1}{2}$-lazy version of the chain, in which each proposal is discarded with probability $\frac{1}{2}$. Since this only affects the mixing time bounds by a factor of $2$, we henceforth ignore this distinction.

The Metropolis-adjusted Langevin algorithm (MALA) has been studied for nearly three decades since~\cite{besagetal1995bayesiancomp}, especially within an asymptotic framework; see, e.g., the influential work of~\cite{robertsrosenthal1998optimalscaling}.
By way of contrast, the non-asymptotic study of the complexity of MALA over the class of strongly log-concave and log-smooth targets is relatively recent.
Through the work of~\cite{dwivedi2018log, chenetal2020hmc, leeshentian2020gradientconcentration, chewi2021optimal, lee2021lower, wuschche2022minimaxmala}, we now know that when initialized at a suitable (and computable) Gaussian initialization, MALA outputs a sample whose law is $\varepsilon$-close to $\pi$ in total variation distance using $\Otilde(\kappa d\,\polylog(1/\eps))$ queries, and when initialized at a warm start the complexity improves to $\Otilde(\kappa d^{1/2}\,\polylog(1/\eps))$; moreover, these rates are tight due to matching lower bounds.
See Table~\ref{tab:mala-progress} and the prior works section \S\ref{ssec:intro:prior} for more discussion.

\subsection{Background on the proximal sampler}\label{app:hiacc:ps}

The next component of our combined high-accuracy sampler (Theorem~\ref{thm:main-slc}) is the \emph{proximal sampler}, first proposed in~\cite{titsias2018auxiliary, leeshentian2021rgo}.
In this algorithm, we first augment the target distribution $\pi$ to a joint distribution $\bs\pi$ on $\R^d\times\R^d$ via
\[
\bs{\pi}(x,y) \propto \exp\Bigl( -f(x) - \frac{1}{2h}\,\|x-y\|^2 \Bigr)\,.
\]
Next, we perform Gibbs sampling on the augmented target; namely, we alternate between the following steps.
\begin{itemize}
	\item \emph{Forwards step (convolution):} Sample $Y_k \mid X_k \sim \pi^{Y\mid X}(\cdot \mid X_k) = \cN(X_k, h I_d)$.
	\item \emph{Backwards step (RGO):} Sample $X_{k+1} \mid Y_k \sim \pi^{X\mid Y}(\cdot \mid Y_k)$. 
\end{itemize}
The backwards step is known as the \emph{restricted Gaussian oracle} (RGO), and any successful use of the proximal sampler hinges upon an efficient implementation of the RGO\@.

Just as the Langevin diffusion is recognized as the Wasserstein gradient flow of the KL divergence w.r.t.\ $\pi$, thanks to the seminal result of~\cite{jordan1998variational}, the proximal sampler admits an appealing interpretation as a proximal discretization of that same gradient flow.
Indeed, the RGO can be interpreted as the proximal operator for the KL divergence over the Wasserstein space, evaluated at a Dirac measure; moreover, the convergence rates established for the proximal sampler exactly match the classical rates for the proximal point method (PPM) in Euclidean space, a fact which led to a new sharp analysis of the PPM under a Polyak--\L{}ojasiewicz inequality. For this and further discussion, we refer to~\cite{chenetal2022proximalsampler}.

From an algorithmic standpoint, we are interested in two key features of the proximal sampler: (1) the proximal reduction framework which allows us to boost the condition number dependence of any high-accuracy sampler to near-linear~\cite{leeshentian2021rgo}, and (2) its uses for sampling beyond strong log-concavity~\cite{chenetal2022proximalsampler}.

\paragraph*{Boosting the condition number dependence of any high-accuracy sampler to near-linear.} When $f$ is $\beta$-smooth, then one can check that the RGO is strongly log-concave and log-smooth, with condition number $O(1)$, as soon as the step size $h$ is taken to be $h\lesssim 1/\beta$. This observation was the starting point of~\cite{leeshentian2021rgo}, who used the proximal sampler to great effect as a \emph{proximal reduction framework}.
Namely, they showed that when $\pi$ is $\alpha$-strongly log-concave, the proximal sampler (with perfect implementation of the RGO), outputs a sample whose law is $\varepsilon$-close in total variation to $\pi$ in $\Otilde(\frac{1}{\alpha h} \log(1/\varepsilon))$ iterations.
This is at most $\Otilde(\kappa \log(1/\varepsilon))$ when we take $h\asymp 1/\beta$.
Moreover, they showed that when the RGO is only implemented \emph{approximately}, 
then the additional TV error incurred can be tracked and controlled, provided that the RGO is implemented via a \emph{high-accuracy} sampler, such as MALA\@.
Hence, the original problem of sampling from the distribution $\pi$ with condition number $\kappa$ is reduced to solving $\widetilde O(\kappa \,\polylog(1/\eps))$ subproblems, each with condition number $O(1)$, to high accuracy.
In particular, when passed through this reduction framework, the condition number dependence of any high-accuracy sampler is automatically boosted to near-linear in $\kappa$.
We take advantage of this framework in passing from the weak version of Theorem~\ref{thm:main-slc} to the full version. (Although one important difference is that we improve this reduction framework to control the error propagation from approximate RGO implementation in each iteration, in order to obtain final guarantees in metrics stronger than TV, namely R\'enyi divergence.)

\paragraph*{Reducing sampling in general settings to the strongly-log-concave setting.} The other virtue of the proximal sampler exploited in this work is its application to sampling beyond the strongly log-concave setting, as was put forth in~\cite{chenetal2022proximalsampler}.
In this paper, the authors analyzed the \emph{outer loop} complexity of the proximal sampler under weak log-concavity and under a variety of standard isoperimetric assumptions (e.g., the log-Sobolev or Poincar\'e inequality).
In all of these cases, however, as soon as $f$ is $\beta$-smooth, the approximate implementation of the RGO can be handled, as before, by a high-accuracy sampler for strongly log-concave and log-smooth distributions.
Hence, the results of~\cite{chenetal2022proximalsampler} effectively reduce the problem of sampling from each of the aforementioned classes of distributions to the problem of high-accuracy strongly log-concave and log-smooth sampling.
They used this reduction to provide new state-of-the-art guarantees for sampling from these classes, and in \S\ref{ssec:hiacc-ext} we leverage our faster implementation of the RGO to improve each of their complexity bounds by a factor of $\sqrt d$.

\subsection{Weak version of Theorem~\ref{thm:main-slc}}\label{app:highacc:slc-weak}

Here, we show the following weaker version of Theorem~\ref{thm:main-slc} as it provides the key building block to prove it. Like Theorem~\ref{thm:main-slc}, this result here shows that from a feasible start, the query complexity of high-accuracy sampling from a strongly-log-concave distribution scales in the dimension $d$ as $\Otilde(d^{1/2})$ rather than $\Otilde(d)$. The difference from Theorem~\ref{thm:main-slc} is that the weaker result here has a suboptimal dependence on the condition number $\kappa$, namely $\Otilde(\kappa^{3/2})$ rather than $\Otilde(\kappa)$. This dependence will later be boosted using the proximal sampler in \S\ref{app:highacc:slc-main}, allowing us to prove the full Theorem~\ref{thm:main-slc}.

\begin{theorem}[Weak version of Theorem~\ref{thm:main-slc}]\label{thm:main-slc-weak}
	Suppose that $\pi \propto \exp(-f)$ where $f$ is $\alpha$-strongly-convex and $\beta$-smooth. Let $x^*$ denote the minimizer of $f$. For any sampling error $\eps > 0$ and any initialization distribution $\delta_{x_0} \in \cP(\R^d)$, there is an algorithm that uses
	\[
	N = \Otilde\bigl( \kappa^{3/2}
	d^{1/2} \, \log( \|x_0 - x^*\| )
	\; + \;
	\kappa d^{1/2} \, \log^3(1/\eps)
	\bigr)
	\]
	first-order queries for $f$ to output a random point in $\R^d$ with law $\mu$ satisfying
	$
	\on{d}(\mu \mmid \pi) \leq \eps
	$
	for any of the following metrics:
	\[
	\on{d} \in \{ \TV, \sqrt{\KL}, \sqrt{\chi^2}, \sqrt{\alpha}\, W_2 \}\,.
	\]
\end{theorem}

\begin{remark}[Handling other metrics]\label{rem:other-metrics}
To prove convergence in the various metrics, due to standard comparison inequalities it usually suffices to prove a convergence result in the strongest metric, namely, the chi-squared divergence.
Indeed, convergence in the KL divergence follows from the monotonicity of R\'enyi divergences (Lemma~\ref{lem:renyi-monotonicity}) and convergence in the TV distance follows from Pinsker's inequality.
If $\pi$ satisfies an LSI, then convergence in $W_2$ follows from Talagrand's $T_2$ inequality (Lemma~\ref{lem:otto-villani}); otherwise, if $\pi$ only satisfies a PI, then convergence in $W_2$ follows from the quadratic transport-variance inequality (Lemma~\ref{lem:quadratic-transport-var}).
\end{remark}

The proof has three steps:
\begin{enumerate}
	\item Run ULMC from this arbitrary initialization $\delta_{x_0}$ to obtain a $\cR_3$ warm start (which also implies a $\chi^2$ warm start). 
	\item Use the $\chi^2$ warm start to argue that MALA mixes rapidly in TV. 
	\item Use the $\cR_3$ warm start to argue that the TV mixing guarantee implies mixing guarantees in $\chi^2$ (and therefore also the other desired metrics by Remark~\ref{rem:other-metrics}).
\end{enumerate}
Our main technical contribution here is the ability to implement step 1 in $\Otilde(d^{1/2})$ queries---this is an immediate application of our result on ULMC (Theorem~\ref{thm:loacc}). Step 2 follows from known results about MALA~\cite{chenetal2022proximalsampler,wuschche2022minimaxmala}, with only minor modification as described below, and step 3 follows from the helper Lemma~\ref{lem:warm-boost-tv}. 
\par The rest of this Appendix section is organized as follows. Step 2 is described in \S\ref{app:highacc:slc-weak-mala}, step 3 is proved in \S\ref{app:highacc:slc-weak-boost}, and then we combine these to prove Theorem~\ref{thm:main-slc-weak} in \S\ref{app:highacc:slc-weak-proof}.

\subsubsection{Rapid mixing of MALA from a R\'enyi warm start}\label{app:highacc:slc-weak-mala}

\par We first formally state the result in step 2, namely, that MALA mixes rapidly in TV from a $\chi^2$ warm start. This is~\cite[Theorem 1]{wuschche2022minimaxmala}, except with a less stringent assumption on the initialization $\mu_0$. Specifically,~\cite[Theorem 1]{wuschche2022minimaxmala} assumes that $\mu_0$ is an $M$-warm start with respect to the target $\pi$ (or equivalently, $\cR_\infty(\mu_0 \mmid \pi) \le \log M$), whereas the following lemma only assumes that $\mu_0$ has bounded $\chi^2$ (or equivalently, bounded $\cR_2$) distance to $\pi$.
This requires only very minor modification to their analysis, but is essential for our purposes, since our result in \S\ref{sec:ulmc} can only produce warm starts in R\'enyi divergences of finite order.

\begin{theorem}[Runtime of MALA from warm start; implicit from~\cite{wuschche2022minimaxmala}]\label{thm:mala-from-warm}
	Let $\pi \propto \exp(-f)$ where $f$ is $\alpha$-strongly-convex and $\beta$-smooth. 
	For any error $\eps \in (0,1)$, the $1/2$-lazy MALA algorithm with appropriate step size requires
	\[
	N = \Otilde\Bigl(\kappa d^{1/2}\, \log^3\Bigl( \frac{\chi^2(\mu_0 \mmid \pi)}{\eps^2} \Bigr) \Bigr)
	\]
	first-order queries to $f$ to output a random point in $\R^d$ whose law $\mu_N$ satisfies $\TV(\mu_N,\pi) \leq \eps$. 
\end{theorem}
\begin{proof}
	We first recall the proof of~\cite[Theorem 1]{wuschche2022minimaxmala}. The majority of their analysis is to bound the $s$-conductance $\Phi_s$. They then appeal to the classical result~\cite[Corollary 1.6]{lovasz1993random}, which states
	\begin{align}
		\|\mu_N - \pi\|_{\TV} \leq H_s + \frac{H_s}{s}\, e^{-N\Phi_s^2/2} \,,
		\label{eq:mala-weak:ls-orig}
	\end{align}
	where $H_s  \deq  \sup \{  \abs{\mu_0(A) - \pi (A)} : \pi(A) \leq s \}$. A consequence of $\mu_0$ being an $M$-warm start with respect to $\pi$ is that $H_s \leq Ms$. By plugging this into the right hand side of~\eqref{eq:mala-weak:ls-orig}, they conclude that the TV error is at most $\eps$ if they set $s = \tfrac{\eps}{2M}$ and $N \geq \frac{2}{\Phi_s^2} \log \frac{2M}{\eps}$. Plugging in their bound~\cite[equation (39)]{wuschche2022minimaxmala} on $\Phi_s$ finishes their proof.
	\par The proof here differs only in that we bound the right hand side of~\eqref{eq:mala-weak:ls-orig} by using the $\chi^2$ warm start. Specifically, observe that for any set $A$, the Cauchy--Schwarz inequality implies
	\[
	\abs{\mu_0(A) - \pi(A)}
	=
	\Bigl\lvert \int \mathds{1}_{A}\, \Bigl( \frac{\D \mu_0}{\D \pi} -1 \Bigr)\, \D \pi\Bigr\rvert
	\leq
	\sqrt{ \int \mathds{1}_A \,\D \pi \cdot \int  \Bigl( \frac{\D \mu_0}{\D \pi} -1 \Bigr)^2\,  \D\pi}
	=
	\sqrt{ \pi(A) \, \chi^2(\mu_0 \mmid \pi)}\,.
	\]
	It follows that $H_s \leq \sqrt{s \, \chi^2(\mu_0 \mmid \pi)}$, and thus we obtain the following modified version of~\eqref{eq:mala-weak:ls-orig}: 
	\begin{align}
		\|\mu_N - \pi\|_{\TV} \leq \sqrt{s\, \chi^2(\mu_0 \mmid \pi)} +
		\sqrt{ \frac{\chi^2(\mu_0 \mmid \pi)}{s}}\, e^{-N\Phi_s^2/2}\,.
		\label{eq:mala-weak:ls-new}
	\end{align}
	Thus the TV error is at most $\eps$ if we set $s = \frac{\eps^2}{4\chi^2(\mu_0 \mmid \pi)}$ and $N = \frac{2}{\Phi_s^2} \log ( \frac{8\chi^2(\mu_0\mmid \pi)}{\eps^2} )$. Plugging in their bound~\cite[equation (39)]{wuschche2022minimaxmala} on $\Phi_s$ completes the proof.\footnote{In fact, one can simply set $M = 2 \chi^2(\mu_0 \mmid \pi)/\eps$ in their final bounds to obtain Theorem~\ref{thm:mala-from-warm}.}
\end{proof}

\subsubsection{Proof of Lemma~\ref{lem:warm-boost-tv}}\label{app:highacc:slc-weak-boost}

By the stationarity property of $P$ and the data-processing inequality for R\'enyi divergences (Lemma~\ref{lem:renyi-process}), we have
$
\cR_3( \mu_n \mmid \pi )
=
\cR_3( \mu_0 P^n \mmid \pi P^n )
\leq
\cR_3( \mu_0 \mmid \pi )
$.
It now suffices to argue that the following inequality holds for any distributions $\mu,\pi$:
\begin{align}
	\chi^2(\mu \mmid \pi) \leq \sqrt{ \TV (\mu,\pi) \cdot \bigl( \exp(2 \cR_3(\mu \mmid \pi)) + 1 \bigr) } \,.
	\label{eq:warm-boost-tv:ineq}
\end{align}
To prove~\eqref{eq:warm-boost-tv:ineq}, we use the Cauchy--Schwarz inequality to bound
\[
\chi^2(\mu \mmid \pi)
=
\int \Bigl\lvert\frac{\D \mu}{\D \pi} - 1\Bigr\rvert^2 \,\D \pi
=
\int \Bigl\lvert \frac{\D\mu}{\D \pi} - 1\Bigr\rvert^{1/2}\, \Bigl\lvert\frac{\D\mu}{\D\pi} - 1\Bigr\rvert^{3/2}\, \D \pi
\leq
\sqrt{ \int \Bigl\lvert\frac{\D\mu}{\D \pi} - 1\Bigr\rvert \,\D \pi \; \cdot \int \Bigl\lvert\frac{\D\mu}{\D\pi} - 1\Bigr\rvert^{3} \,\D \pi }\,.
\]
The first integral is precisely $\TV(\mu, \pi)$. The second integral can be bounded by
\[
\int \Bigl\lvert \frac{\D\mu}{\D\pi} - 1\Bigr\rvert^{3} \,\D \pi
\leq
\int \Bigl\lvert\frac{\D\mu}{\D\pi}\Bigr\rvert^{3} \,\D \pi + 1
=
\exp( 2 \cR_3(\mu \mmid \pi) ) + 1\,,
\]
where above, the first step is by the elementary inequality $\abs{a-1}^3 \leq a^3 + 1$, which holds for all $a \geq 0$; and the second inequality is by the definition of R\'enyi divergence. This completes the proof of~\eqref{eq:warm-boost-tv:ineq} and thus also the proof of the lemma.

\subsubsection{Proof of Theorem~\ref{thm:main-slc-weak}}\label{app:highacc:slc-weak-proof}

By Theorem~\ref{thm:loacc}---or rather the extension in Remark~\ref{rem:loacc-init} to arbitrary initial distributions---ULMC outputs a distribution $\nu$ satisfying $\cR_3(\nu \mmid \pi) \leq \log 2$, say, using 
\[
\Otilde\bigl(\kappa^{3/2} d^{1/2} \, \log \|x_0 - x^*\| \bigr)
\]
gradient queries, where $\delta_{x_0}$ is its initial distribution. By monotonicity of R\'enyi divergences (Lemma~\ref{lem:renyi-monotonicity}) and the identity between $\chi^2$ and $\cR_2$ (Remark~\ref{rem:renyi-cases}), this ULMC guarantee implies $\chi^2(\nu \mmid \pi) = \exp(\cR_2(\nu \mmid \pi)) - 1 \leq \exp(\cR_3(\nu \mmid \pi)) - 1 \leq 1$, so $\nu$ is a warm start in $\chi^2$ divergence. Thus we may invoke Theorem~\ref{thm:mala-from-warm} to run MALA from initialization $\nu$ in order to produce a distribution $\mu$ satisfying $\TV(\mu, \pi) \leq \eps^4/5$, say, using
\[
\Otilde\bigl( \kappa d^{1/2} \, \log^3(1/\eps)\bigr)
\]
first-order queries. Now by Lemma~\ref{lem:warm-boost-tv}, we can use the warm start property of $\nu$ to boost the $\TV$ guarantee on MALA's output $\mu$ to the following $\chi^2$ guarantee:
\[
\chi^2(\mu \mmid \pi) 
\leq 
\sqrt{ \TV(\mu,\pi) \cdot \bigl(\exp(2\cR_3 (\nu \mmid \pi)) + 1\bigr)}
\leq 
\eps^2\,.
\]
This implies the desired $\chi^2$ mixing bound.
Mixing in the other metrics then follows from Remark~\ref{rem:other-metrics}.

\subsection{Proof of Theorems~\ref{thm:main-slc} and~\ref{thm:main-lsi}}\label{app:highacc:slc-main}

Here we prove our main results about faster high-accuracy sampling algorithms in the setting that the target distribution $\pi$ is strongly-log-concave (Theorem~\ref{thm:main-slc}) or satisfies a log-Sobolev inequality (Theorem~\ref{thm:main-lsi}). Since our analysis only relies upon the LSI property, we are able to prove both theorems simultaneously. (Indeed, recall that strong-log-concavity implies a log-Sobolev inequality by the Bakry--\'Emery theorem, see the first part of Lemma~\ref{lem:bakry-emery}). See \S\ref{ssec:hiacc-ext} for a high-level overview of the algorithm and analysis.

\par We begin with a helper lemma, which is similar to the Orlicz--Wasserstein initialization bound for $\pi$ in Lemma~\ref{lem:ulmc-init}, but now generalized to the RGO $\pi^{X|Y=y} \propto \exp(-f - \frac{1}{2h}\,\|\cdot-y\|^2)$ that is used in the backwards step of the proximal sampler. 

\begin{lemma}[Orlicz--Wasserstein distance at initialization of RGO step]\label{lem:rgo-init}
	Suppose that $\pi \propto \exp(-f)$ where $f$ is $\beta$-smooth. Let $x^*$ denote the mode of $\pi$. Then for any $y \in \R^d$ and any proximal step size $h \leq 1/(2\beta)$, 
	\[
	W_{\psi_2}(\delta_{y},\pi^{X|Y=y}) 
     \le 9\sqrt{dh} + 3\beta h\,\norm{y-x^*}\,.
	\]
\end{lemma}
\begin{proof}
	Let $x_y$ denote the mode of $\pi^{X|Y=y}$. By the triangle inequality,
	\[
	W_{\psi_2}( \pi^{X|Y=y}, \delta_{y})
	\leq
	W_{\psi_2}( \pi^{X|Y=y} , \delta_{x_y} ) +
	W_{\psi_2}( \delta_{x_y}, \delta_{y} )\,.
	\]
	The former term is bounded above by $9\sqrt{dh}$ by an application of Lemma~\ref{lem:ulmc-init} and the observation that $\pi^{X|Y=y} \propto \exp(-f - \frac{1}{2h}\,\|\cdot-y\|^2)$ is strongly-log-concave with parameter $-\beta + \frac{1}{h} \geq \frac{1}{2h}$. Next, we bound the latter term
	\[
	W_{\psi_2}(\delta_{x_y}, \delta_{y})
	=
	\|x_y - y\|_{\psi_2}
	=
	\frac{\|x_y - y\|}{\sqrt{\log 2}}\,.
	\]
	Since $x_y$ is the mode of $\pi^{X|Y=y}$, it is the minimizer of the convex log-density, thus by first-order optimality conditions we have $0 = \nabla f(x_y) + \frac{1}{h}\,(x_y - y)$. By rearranging this identity, using the smoothness of $f$, and then using the triangle inequality,
	\begin{align*}
		\|x_y - y\|
		= h\, \|\nabla f(x_y) - \nabla f(x^*) \|
		\leq \beta h\, \|x_y - x^*\|
		\leq \beta h\, \bigl( \|x_y - y\| + \|y - x^*\| \bigr)\,.
	\end{align*}
	Now by the assumption on the step size, $\beta h \leq 1/2$. Plugging this in and re-arranging yields
	\[
	\| x_y - y \| \leq 2\beta h\, \|y - x^*\|\,.
	\]
	Combining the above displays completes the proof.
\end{proof}

\par Armed with this initialization lemma, we are now ready to prove the main results of this section.

\begin{proof}[Proof of Theorems~\ref{thm:main-slc} and~\ref{thm:main-lsi}]
    Recall from the discussion at the beginning of this subsection that it suffices to prove Theorem~\ref{thm:main-lsi}. Hence, in this proof we assume that $\pi$ is $1/\alpha$-LSI but do not necessarily assume that it is $\alpha$-strongly-log-concave. We prove the mixing time for the $\chi^2$ divergence, which suffices by Remark~\ref{rem:other-metrics}.
	\par First, suppose that the RGO in the proximal sampler algorithm is implemented exactly. Let $X_n$ and $Y_n$ denote the proximal sampler iterates at iteration $n$; and let $\mu_n^X$ and $\mu_n^Y$ denote their respective laws. Then after initializing at $\mu_0^X \defeq \cN(x^*, (2\beta)^{-1} I_d)$, the laws of the iterates are given by $\mu_n^{Y} = \mu_n^X \ast \cN(0, hI_d)$, and $\mu_{n+1}^X = \int \pi^{X\mid Y}(\cdot \mid y) \,\mu_n^{Y}(\D y)$. By analyzing the simultaneous heat flow, it was shown in~\cite[Appendix A.4]{chenetal2022proximalsampler} that the forwards step of the proximal algorithm is a contraction in R\'enyi divergence, in the sense that
	\begin{align}
		\cR_{q}( \mu_n^Y \mmid \pi^Y) \leq \frac{1}{{(1 + \alpha h)}^{1/q}} \; \cR_q( \mu_n^X \mmid \pi^X)\,.
		\label{eq:lsi-proof:forwards}
	\end{align}
	\par Now suppose that we have oracle access to an approximate RGO in the sense that given any point $y \in \R^d$, we can sample from a distribution $\tilde{\pi}^{X|Y=y}$ that satisfies
	\begin{align}
		\cR_q( \tilde{\pi}^{X|Y=y}  \mmid \pi^{X|Y=y}) \leq \epsrgo^2\,,
		\label{eq:lsi-proof:approx-rgo}
	\end{align}
	using $\Nrgo(y)$ first-order queries to $f$. 
 Let $\tilde X_n$, $\tilde Y_n$ denote the iterates with inexact implementation of the RGO\@, and let $\tilde \mu_n^X$, $\tilde \mu_n^Y$ denote their laws respectively.
	\par We can bound the error of a backwards step using this approximate RGO as follows. Let $\tilde Y_n \sim \tilde \mu_n^Y$, so that $\tilde{X}_{n+1}$ is a sample from the approximate RGO $\tilde{\pi}^{X|Y=\tilde Y_n}$. Then
	\begin{align}
		\cR_q ( \tilde{\mu}_{n+1}^X \mmid \pi^X)
		&\leq
		\cR_q \bigl( \law(\tilde{X}_{n+1}, \tilde Y_n) \bigm\Vert \bs{\pi} \bigr) \nonumber
		\\ &\leq
		\cR_q(\tilde \mu_n^Y \mmid \pi^Y) + \sup_{y_n\in\R^d} \cR_q( \tilde{\pi}^{X|Y=y_n} \mmid \pi^{X|Y=y_n}) \nonumber
		\\ &\leq \cR_q(\tilde \mu_n^Y \mmid \pi^Y) + \epsrgo^2\,.
		\label{eq:lsi-proof:approx-backwards}
	\end{align}
	Above, the first step is by the data-processing inequality for R\'enyi divergences (Lemma~\ref{lem:renyi-process}); the second step is by the ``strong composition rule'' for R\'enyi differential privacy (this lemma has appeared in many equivalent forms, see, e.g.,~\cite{abadi2016deep, dwork2016concentrated, mironov2017renyi}; here we apply the version from~\cite[Lemma 2.9]{AltTal22dp}); and the final step is by the guarantee~\eqref{eq:lsi-proof:approx-rgo} of the approximate RGO.
	\par By combining the error bounds~\eqref{eq:lsi-proof:forwards} and~\eqref{eq:lsi-proof:approx-backwards} for the forward step and approximate backwards step of the proximal sampler, we conclude the one-iteration bound
	\begin{align}
		\cR_q( \tilde{\mu}_{n+1}^X \mmid \pi^X)
		\leq
		\frac{1}{{(1 + \alpha h)}^{1/q}}\, \cR_q(\tilde \mu_n^X \mmid \pi^X) + \epsrgo^2\,.
		\label{eq:lsi-proof:approx-1step}
	\end{align} 
	Iterating this bound $\Nprox$ times gives the following R\'enyi divergence bound on the mixing error of the proximal sampler when using this approximate RGO:
	\begin{align}
		\cR_q( \tilde{\mu}_{\Nprox}^X \mmid \pi^X \big)
		\leq
		\frac{1}{{(1 + \alpha h)}^{\Nprox/q}}\, \cR_q(\mu_0^X \mmid \pi^X) + \epsrgo^2 \sum_{n=0}^{\Nprox-1} \frac{1}{{(1 + \alpha h)}^{n/q}} \,.
		\label{eq:inexact_renyi_bd}
	\end{align}
	This error is at most $\eps^2$ if we run the proximal sampler with step size $h \asymp 1/\beta$ for
	\[
	\Nprox \asymp \kappa q \log \frac{\cR_q ( \mu_0^X \mmid \pi^X )}{\eps^2}
	\]
	iterations and perform each approximate RGO to accuracy 
	\begin{align}\label{eq:inexact_rgo_acc}
		\epsrgo \asymp \frac{\eps}{\sqrt{\kappa q}}\,.
	\end{align}
	\par Henceforth, consider $q = 2$, so that $\cR_2 \leq \chi^2$ (see Remark~\ref{rem:renyi-cases}). Observe that if the step size $h < 1/(2\beta)$, say, then the RGO is strongly-log-concave and has condition number of size at most
	\[
	\frac{\beta + 1/h}{-\beta + 1/h} = \frac{1 + \beta h}{1 - \beta h} = \Theta(1)\,.
	\]

    We next consider the complexity of implementing the RGO\@.
    If, as we have assumed in \S\ref{sec:hiacc}, we can compute the proximal operator for $hf$ exactly, then we can compute the mode $x({\tilde Y_n})$ of $\pi^{X\mid Y=\tilde Y_n}$ and initialize at $\nu_n = \delta_{x(\tilde Y_n)}$.
    Otherwise, we initialize at $\nu_n = \delta_{\tilde Y_n}$.
    In either case, by Theorem~\ref{thm:main-slc-weak}, we can implement the approximate RGO $\tilde{\pi}^{X|Y = \tilde Y_n}$ in the $n$-th iteration by using $\Nrgo(\tilde Y_n)$ first-order queries, where 
	\begin{align}
		\Nrgo(\tilde Y_n)
		=
		\Otilde\Bigl( d^{1/2} \, \log^3\Bigl( \frac{ W_{\psi_2} (\nu_n, \; \pi^{X|Y=\tilde Y_n})}{\epsrgo} \Bigr) \Bigr)\,.
		\label{eq:lsi-proof:N-RGO}
	\end{align}
	By Lemma~\ref{lem:rgo-init}, $W_{\psi_2}(\delta_{x(\tilde Y_n)},\;\pi^{X\mid Y=\tilde Y_n}) \lesssim \sqrt{d/\beta}$ and $W_{\psi_2} (\delta_{\tilde Y_n}, \; \pi^{X|Y=\tilde Y_n}) \lesssim \sqrt{d/\beta} + \|\tilde Y_n-x^*\|$.
     In the former case, we conclude that the total number of gradient queries required by this inexact proximal sampler algorithm is
	\begin{align}
		N 
		=
		\sum_{n=0}^{\Nprox-1} \Nrgo(\tilde Y_n)
		=
		\Otilde\biggl(\kappa d^{1/2} \, \log\Bigl(\frac{\cR_2(\mu_0^X \mmid \pi^X)}{\eps^2}\Bigr)\,  \log^3\Bigl(\frac{1}{\eps}\Bigr) \biggr)\,.
		\label{eq:N_rgo_with_prox}
	\end{align}
    Otherwise, the number of gradient queries is
	\begin{align}
		N 
		=
		\sum_{n=0}^{\Nprox-1} \Nrgo(\tilde Y_n)
		=
		\Otilde\biggl(\kappa d^{1/2} \, \log\Bigl(\frac{\cR_2(\mu_0^X \mmid \pi^X)}{\eps^2}\Bigr)\,  \log^3\Bigl(\frac{\max_{n=0,1,\dotsc,\Nprox-1}{\norm{\tilde Y_n - x^*}}}{\eps}\Bigr) \biggr)\,.
		\label{eq:N_rgo}
	\end{align}
    The latter expression will be made more explicit in Appendix~\ref{app:hiacc:explicit}, and upon doing so it leads to the final statement of Theorem~\ref{thm:main-slc}.
\end{proof}

\subsection{Proof of Theorem~\ref{thm:main-poincare}}\label{app:high-acc:poincare}

We prove the $\chi^2$ mixing bound; the other desired mixing bounds then follow immediately due to standard comparison inequalities (see Remark~\ref{rem:other-metrics}). We consider the same inexact RGO algorithm as in the LSI setting (see Appendix~\ref{app:highacc:slc-main}). Under the present Poincar\'e assumption, the forwards step of the proximal algorithm is known to be a contraction in $\chi^2$---in direct analog to~\eqref{eq:lsi-proof:forwards}. Specifically, by analyzing the simultaneous heat flow, it was shown in~\cite[Appendix A.4]{chenetal2022proximalsampler} that
\begin{align}
	\chi^2( \mu_n^Y \mmid \pi^Y) \leq \frac{1}{1 + \alpha h} \; \chi^2 ( \mu_n^X \mmid \pi^X)\,.
	\label{eq:poincare-proof:forwards}
\end{align}
The bound~\eqref{eq:lsi-proof:approx-backwards} on the error of a backwards step of the proximal sampler using an approximate RGO~\eqref{eq:lsi-proof:approx-rgo} remains unchanged (as it never uses the LSI assumption). This R\'enyi bound is equivalent to the $\chi^2$ bound
\begin{align}
	\chi^2( \tilde{\mu}_{n+1}^X \mmid \pi^X)
	\leq 
	e^{\epsrgo^2}\,
	\chi^2(\tilde \mu_n^Y \mmid \pi^Y)\,,
	\label{eq:poincare-proof:approx-backwards}
\end{align}
by using the relationship $\cR_2 = \log(1+\chi^2)$ between the chi-squared and R\'enyi divergences (see Remark~\ref{rem:renyi-cases}). By combining the above two displays, we obtain the following convergence bound for one full step of the proximal sampler: 
\begin{align}
	\chi^2( \tilde{\mu}_{n+1}^X \mmid \pi^X)
	\leq
	e^{- \Theta(\alpha h) } \,
	\chi^2(\tilde \mu_n^X \mmid \pi^X) \,,
	\label{eq:poincare-proof:approx-1step}
\end{align} 
if we solve each approximate RGO to accuracy
\begin{align}
	\epsrgo \lesssim \sqrt{\alpha h} \asymp \frac{1}{\sqrt\kappa} \,.
\end{align}
By iterating this one-step bound~\eqref{eq:poincare-proof:approx-1step}, we conclude that the final mixing error $\chi^2(\tilde{\mu}_{\Nprox}^X \mmid \pi^X)$ of this inexact proximal sampler is at most $\eps^2$ if it is run for $\Nprox$ iterations, where 
\[
\Nprox \asymp \kappa \log  \Bigl( \frac{\chi^2 ( \mu_0^X \mmid \pi^X )}{\eps^2} \Bigr)\,.
\]
\par Now by the same argument as the LSI setting (see Appendix~\ref{app:highacc:slc-main}), by the choice of the proximal step size $h$, the RGO $\tilde{\pi}^{X\mid Y=\tilde Y_n}$ has condition number $\Theta(1)$, and thus can be implemented using $\Nrgo(\tilde Y_n)$ gradient queries by Theorem~\ref{thm:main-slc-weak}, where $\Nrgo(\tilde Y_n)$ is the quantity in~\eqref{eq:lsi-proof:N-RGO}. Therefore the total number of gradient queries required by this inexact proximal sampler algorithm is
\begin{align}
	N 
	=
	\sum_{n=0}^{\Nprox-1} \Nrgo(\tilde Y_n)
	=
	\Otilde\biggl( \kappa d^{1/2} \, \log\Bigl( \frac{\chi^2(\mu_0^X \mmid \pi^X)}{\varepsilon^2}\Bigr)\, \log^3\max_{n=0,1,\dotsc,\Nprox-1} {\|\tilde Y_n - x^*\|}  \biggr)\,.
	\label{eq:N_rgo-poincare}
\end{align}
The term $\log^3\max_{n=0,1,\dotsc,\Nprox-1} {\|\tilde Y_n - x^*\|}$ vanishes if we have access to the proximal operator for $hf$; otherwise, it is made more explicit in Appendix~\ref{app:hiacc:explicit}.

\subsection{Proof of Theorem~\ref{thm:main-logconcave}}\label{app:high-acc:logconcave}

The proof for the weakly convex case is similar to the proofs of Theorems~\ref{thm:main-slc},~\ref{thm:main-lsi}, and~\ref{thm:main-poincare}, in that we carefully keep track of the error from inexact implementation of the RGO\@, but the proof requires key modifications.
It was shown in~\cite[Appendix A.3]{chenetal2022proximalsampler} that along the simultaneous heat flow,
\begin{align}\label{eq:prox-sampler-lc-step}
    \frac{1}{\KL(\mu_n^Y \mmid \pi^Y)} \ge \frac{1}{\KL(\mu_n^X \mmid \pi^X)} + \frac{h}{2W_2^2(\mu_n^X, \pi^X)}\,.
\end{align}
Let us assume that the RGO is implemented inexactly, so that for each $y\in\R^d$ we sample from $\tilde \pi^{X\mid Y=y}$ satisfying
\begin{align}\label{eq:inexact-rgo-lc}
    \KL(\tilde\pi^{X\mid Y=y} \mmid \pi^{X\mid Y=y}) \le \epsrgo^2\,, \qquad W_2(\tilde \pi^{X\mid Y=y}, \, \pi^{X\mid Y=y}) \le \sqrt{2\beta}\,\epsrgo\,.
\end{align}
In fact, the second guarantee follows from the first together with Talagrand's $T_2$ inequality (see Lemma~\ref{lem:otto-villani}) if we choose step size $h = \frac{1}{2\beta}$, because the RGO is then $\beta$-strongly log-concave.

By applying~\eqref{eq:lsi-proof:approx-backwards} for the proximal sampler with inexact RGO implementation, convexity of the map $x\mapsto 1/x$, and~\eqref{eq:prox-sampler-lc-step}, we deduce that
\begin{align}
    \frac{1}{\KL(\tilde\mu_{n+1}^X \mmid \pi^X)}
    \ge \frac{1}{\KL(\tilde \mu_n^Y \mmid \pi^Y) + \epsrgo^2}
    &\ge \frac{1}{\KL(\tilde\mu_n^Y \mmid \pi^Y)} - \frac{\epsrgo^2}{{\KL(\tilde\mu_n^Y \mmid \pi^Y)}^2} \nonumber \\
    &\ge \frac{1}{\KL(\tilde\mu_n^X \mmid \pi^X)} + \frac{h}{2W_2^2(\tilde\mu_n^X,\pi^X)} - \frac{\epsrgo^2}{{\KL(\tilde\mu_n^Y \mmid \pi^Y)}^2}\,. \label{eq:prox-sampler-lc-inexact-step}
\end{align}

We now split into two cases.
In the first case, suppose that $\KL(\tilde\mu_n^Y \mmid \pi^Y) \le \sqrt{\epsrgo}$ for some $n=0,1,\dotsc,\Nprox-1$.
By repeatedly applying~\eqref{eq:lsi-proof:approx-backwards} and the data-processing inequality (Lemma~\ref{lem:renyi-process}), we obtain
\begin{align*}
    \KL(\tilde\mu_{\Nprox}^X \mmid \pi^X)
    &\le \KL(\tilde \mu_{\Nprox-1}^Y \mmid \pi^Y) + \epsrgo^2 \\
    &\le \KL(\tilde\mu_{\Nprox-1}^X \mmid \pi^X) + \epsrgo^2
    \le \cdots \\
    &\le \KL(\tilde\mu_n^Y \mmid \pi^Y) + \Nprox\,\epsrgo^2
    \le \sqrt{\epsrgo} + \Nprox\,\epsrgo^2\,.
\end{align*}

For the other case, suppose that $\KL(\tilde\mu_n^Y \mmid \pi^Y) \ge \sqrt{\epsrgo}$ for all $n=0,1,\dotsc,\Nprox-1$.
Then, from~\eqref{eq:prox-sampler-lc-inexact-step}, we have
\begin{align*}
    \frac{1}{\KL(\tilde\mu_{n+1}^X \mmid \pi^X)}
    &\ge \frac{1}{\KL(\tilde\mu_n^X \mmid \pi^X)} + \frac{h}{2W_2^2(\tilde\mu_n^X,\pi^X)} - \epsrgo\,.
\end{align*}
Iterating this,
\begin{align*}
    \frac{1}{\KL(\tilde\mu_{\Nprox}^X \mmid \pi^X)}
    \ge \frac{1}{\KL(\mu_0^X \mmid \pi^X)} + \frac{h}{2} \sum_{n=0}^{\Nprox-1} \frac{1}{W_2^2(\tilde\mu_n^X, \pi^X)} - \Nprox\, \epsrgo\,.
\end{align*}
Moreover, from the second condition in~\eqref{eq:inexact-rgo-lc}, a standard coupling argument (see, e.g.,~\cite[Appendix A.2]{chenetal2022proximalsampler}), and Wasserstein contractivity of the exact proximal sampler steps under log-concavity~\cite[Theorem 1]{chenetal2022proximalsampler}, we obtain
\begin{align*}
    W_2(\tilde\mu_{n+1}^X, \pi^X)
    &\le W_2\Bigl(\tilde\mu_{n+1}^X, \,\int \pi^{X\mid Y=y}\, \tilde\mu_n^Y(\D y)\Bigr) + W_2\Bigl(\int \pi^{X\mid Y=y}\, \tilde\mu_n^Y(\D y), \, \pi^X\Bigr) \\
    &\le \sqrt{2\beta}\, \epsrgo + W_2(\tilde\mu_n^Y, \pi^Y) \\
    &\le \sqrt{2\beta}\, \epsrgo + W_2(\tilde\mu_n^X, \pi^X)
    \le \cdots \\
    &\le \sqrt{2\beta}\,\Nprox\,\epsrgo + W_2(\mu_0^X, \pi^X)\,.
\end{align*}
Therefore, we obtain
\begin{align*}
    \frac{1}{\KL(\tilde\mu_{\Nprox}^X \mmid \pi^X)}
    \ge \frac{1}{\KL(\mu_0^X \mmid \pi^X)} + \frac{\Nprox\,h}{2\,{(W_2(\mu_0^X, \pi^X) + \sqrt{2\beta}\,\Nprox\,\epsrgo)}^2} - \Nprox\, \epsrgo\,.
\end{align*}
Let us assume that $\epsrgo \lesssim W_2(\mu_0^X,\pi^X)/(\sqrt\beta\,\Nprox)$ and $\epsrgo \lesssim h/W_2(\mu_0^X,\pi^X)$.
This reads
\begin{align*}
    \frac{1}{\KL(\tilde\mu_{\Nprox}^X \mmid \pi^X)}
    \ge \frac{1}{\KL(\mu_0^X \mmid \pi^X)} + \Omega\Bigl( \frac{\Nprox\,h}{2W_2^2(\mu_0^X, \pi^X)}\Bigr)\,.
\end{align*}
Upon rearranging this and taking $h=\frac{1}{2\beta}$, it implies
\begin{align*}
    \KL(\tilde\mu_{\Nprox}^X \mmid \pi^X)
    &\lesssim \frac{\beta W_2^2(\mu_0^X, \pi^X)}{\Nprox}\,.
\end{align*}
Therefore, we obtain $\KL(\tilde\mu_{\Nprox}^X \mmid \pi^X) \le \eps^2$ if we take
\begin{align*}
    \Nprox
    &\asymp \frac{\beta W_2^2(\mu_0^X, \pi^X)}{\eps^2}\,.
\end{align*}

To summarize, after considering both cases, we obtain $\KL(\tilde\mu_{\Nprox}^X \mmid \pi^X) \le \eps^2$ if we take
\begin{align*}
    \Nprox
    \asymp \frac{d+\beta \mf m^2}{\eps^2} \qquad\text{and}\qquad \epsrgo \lesssim \min\Bigl\{ \eps^4, \, \frac{\eps^4}{\beta W_2^2(\mu_0^X, \pi^X)}, \, \frac{\eps^2}{\beta^{3/2} W_2(\mu_0^X,\pi^X)}, \, \frac{1}{\beta W_2(\mu_0^X,\pi^X)}\Bigr\}\,.
\end{align*}

By invoking Theorem~\ref{thm:main-slc-weak}, the total number of first-order queries is
\begin{align*}
    N
    &= \Otilde\biggl( \frac{\beta d^{1/2}\, W_2^2(\mu_0^X, \pi^X)}{\eps^2} \log^3 \max_{n=0,1,\dotsc,\Nprox-1}{\norm{\tilde Y_n - x^*}} \biggr)\,.
\end{align*}
The term $\log^3\max_{n=0,1,\dotsc,\Nprox-1} {\|\tilde Y_n - x^*\|}$ vanishes if we have access to the proximal operator for $hf$; otherwise, it is made more explicit in Appendix~\ref{app:hiacc:explicit}.

\subsection{Explicit bounds}\label{app:hiacc:explicit}

Here, we make the statements of the results in \S\ref{ssec:hiacc-ext} more explicit by bounding the initialization quantities in terms of other, more easily computable problem parameters. Moreover, we carry through the analysis without assuming access to a prox oracle for $hf$.

To do so, we instead assume that the algorithm has access to a stationary point $x^*$ of $f$, which is realistic even for the non-convex potentials considered in \S\ref{ssec:hiacc-ext}, and we let $\mf m \deq \E_\pi\norm{\cdot - x^*}$.
Our bounds will depend on $\mf m$ as well as on the objective gap $\Delta \deq f(x^*) - \min f$, and we will provide comments on these assumptions after stating the results.

We remark that these explicit bounds also involve randomized runtimes. The reason for this is that in each iteration of the proximal sampler, the implementation of the RGO takes a number of queries which depends on the size of the proximal sampler iterate; in turn, this is a random quantity.

\subsubsection{Initialization}

We make use of the following R\'enyi divergence bound at initialization when using a Gaussian ``feasible start''.  This bound is imported from~\cite[Lemma 30]{chewi2021analysis}\footnote{Their lemma has an extra $\beta$ because it is written for the general setting of any H\"older smoothness exponent $s$. Lemma~\ref{lem:renyi-gaussian-init} is obtained by setting $s=1$, in which case the extra $\beta$ trivially drops in the second line of their proof.} and can be thought of as a tighter, more explicit version of~\cite[Lemma 4]{VempalaW19}. Observe that this lemma does not require convexity of the potential $f$, which makes it applicable to the LSI and PI settings in \S\ref{sec:hiacc}. If $f$ is assumed convex, then every stationary point $x_0$ is a minimizer of $f$, hence the upper bound improves since $f(x_0) - \min f$ vanishes.

\begin{lemma}[R\'enyi divergence bound from Gaussian initialization]\label{lem:renyi-gaussian-init}
	Suppose $\pi \propto \exp(-f)$ where $f$ is $\beta$-smooth. Denote $\mf m = \int \|\cdot\| \,\D\pi$. For any stationary point $x_0$ of $f$,
	\[
	\cR_{\infty}\bigl(\mc N(x_0, (2\beta)^{-1} I_d) \bigm\Vert \pi \bigr) 
	\leq
	2 + f(x_0) - \min f + \frac{d}{2} \log (2 \beta \mf m^2)\,.
	\]
 In particular, if $f$ is also convex, then
 \[
	\cR_{\infty}\bigl(\mc N(x_0, (2\beta)^{-1} I_d) \bigm\Vert \pi \bigr) 
	\leq
	2 + \log (2 \beta \mf m^2)\,.
	\]
\end{lemma}

\subsubsection{Bounding the size of the proximal sampler iterates}

In our high-accuracy sampling results in \S\ref{app:highacc:explicit:results}, the algorithm runtimes will be random because they depend on the largest size $\max_{n=0,1,\dotsc,\Nprox-1}{\norm{\tilde Y_n - x^*}}$ of the iterates of the proximal sampler.
The following lemma provides a general bound on this quantity to ensure that the runtime of the overall randomized algorithm is typically small.
Since the final dependence on this quantity is only polylogarithmic, we will only focus on obtaining crude polynomial bounds.

\begin{lemma}[Size of the proximal sampler iterates]\label{lem:size-prox-iterates}
    Suppose that we run the proximal sampler algorithm for $\Nprox$ iterations with inexact implementation of the RGO\@.
    Namely, for $n=0,1,\dotsc,\Nprox-1$, we set $\tilde Y_n \sim \cN(\tilde X_n, hI_d)$ and $\tilde X_{n+1} \sim \tilde \pi^{X\mid Y=\tilde Y_n}$, where for each $y\in\R^d$, $\tilde\pi^{X\mid Y=y}$ satisfies
    \begin{align*}
        \cR_2(\tilde\pi^{X\mid Y=y} \mmid \pi^{X\mid Y=y}) \le \epsrgo^2\,.
    \end{align*}

    Assume that the target $\pi \propto\exp(-f)$ satisfies $1/\alpha$-PI\@, that $x^*$ is a stationary point of $f$, and that $\mf m \deq \E_\pi\norm{\cdot - x^*}$.
    Then, for $\delta \in (0,1)$, with probability at least $1-\delta$, it holds that
    \begin{align*}
        \max_{n=0,1,\dotsc,\Nprox-1}{\norm{\tilde Y_n - x^*}}
        &\le \on{poly}\Bigl(\epsrgo, \, C_{\msf{PI}}, \, \Nprox, \,d,\, h,\, \mf m,\,\cR_3(\mu_0^X \mmid \pi^X), \,\log \frac{1}{\delta}\Bigr)\,.
    \end{align*}
\end{lemma}
\begin{proof}
    Let $\tilde\mu_n^X \deq \law(X_n)$ and $\tilde\mu_n^Y \deq \law(Y_n)$. We first give a crude bound on $\cR_{3/2}(\tilde\mu_n^Y \mmid \pi^Y)$.
    Let $X_n$, $Y_n$ denote the iterates of the proximal sampler with exact RGO implementation at iteration $n$, with respective laws $\mu_n^X$, $\mu_n^Y$.
    By the weak triangle inequality (Lemma~\ref{lem:renyi-triangle}) with $q=3/2$ and $\lambda = 3/4$, together with the data-processing inequality (Lemma~\ref{lem:renyi-process}),
    \begin{align*}
        \cR_{3/2}(\tilde\mu_n^Y \mmid \pi^Y)
        &\le \frac{3}{2}\, \cR_2(\tilde\mu_n^Y \mmid \mu_n^Y) + \cR_3(\mu_n^Y \mmid \pi^Y) \\
        &\le \frac{3}{2}\, \cR_2(\tilde\mu_n^Y \mmid \mu_n^Y) + \cR_3(\mu_0^X \mmid \pi^X)\,.
    \end{align*}
    Next, by the strong composition rule in R\'enyi differential privacy,
    \begin{align*}
        \cR_2(\tilde\mu_{n+1}^Y \mmid \mu_{n+1}^Y)
        &\le \cR_2(\tilde\mu_{n+1}^X \mmid \mu_{n+1}^X)
        \le \cR_2(\tilde\mu_n^Y \mmid \mu_n^Y) + \epsrgo^2\,.
    \end{align*}
    By iterating this inequality, we deduce that for all $n\le \Nprox-1$,
    \begin{align}\label{eq:crude-renyi}
        \cR_{3/2}(\tilde\mu_n^Y \mmid \pi^Y)
        &\le \frac{3\Nprox\, \epsrgo^2}{2} + \cR_3(\mu_0^X \mmid \pi^X)\,.
    \end{align}
    This inequality will be used later for a change of measure argument.

    However, we must first investigate the concentration of $\norm{\cdot-x^*}$ under $\pi^Y$.
    Since $\norm{\cdot - x^*}$ is a $1$-Lipschitz function, then concentration under a Poincar\'e inequality (Lemma~\ref{lem:pi-subexp}) implies the following tail bound: for all $\eta \ge 0$,
    \begin{align*}
        \pi\bigl\{\norm{\cdot - x^*} \ge \mf m +\sqrt{C_{\msf{PI}}}\, \eta\bigr\}
        &\le 3\exp(-\eta)\,.
    \end{align*}
    Also, let $\rho_h = \mc N(0,hI_d)$. Standard concentration estimates for $\rho_h$ yield
	\begin{align*}
		\rho_h\{\norm \cdot \ge \sqrt{dh} + \eta\} \le \exp\bigl( -\frac{\eta^2}{2h}\bigr)\,.
	\end{align*}
	Since $\pi^Y = \pi^X * \rho_h$, a union bound yields
	\begin{align*}
		\pi^Y\bigl\{\norm{\cdot - x^*} \ge \mf m +\sqrt{dh} + \sqrt{C_{\msf{PI}}}\,\eta + \sqrt{2h\eta}\bigr\}
		\le 4\exp(-\eta)
	\end{align*}
     which implies the simpler bound
     \begin{align*}
         \pi^Y\bigl\{\norm{\cdot-x^*}\ge \mf m + 3\sqrt{dh} + \sqrt{2\,(C_{\msf{PI}} + h)}\,\eta \bigr\} \le 4\exp(-\eta)\,.
     \end{align*}
    
	\par We now adapt the change-of-measure argument from~\cite[Lemma 21]{chewi2021analysis}.
        Namely, let $E_\eta$ denote the event above.
        By H\"older's inequality,
        \begin{align*}
            \tilde\mu_n^Y(E_\eta)
            &= \int \one_{E_\eta}\,\frac{\D\tilde\mu_n^Y}{\D\pi^Y} \, \D \pi^Y
            \le \pi^Y(E_\eta)^{1/3} \,\Bigl( \int \bigl( \frac{\D\tilde\mu_n^Y}{\D\pi^Y}\bigr)^{3/2}\,\D\pi^Y\Bigr)^{2/3}
            = \pi^Y(E_\eta)^{1/3} \exp\Bigl( \frac{\cR_{3/2}(\tilde\mu_n^Y,\pi^Y)}{2}\Bigr)\,.
        \end{align*}
        This implies that
        \begin{align*}
            \tilde\mu_n^Y\Bigl\{\norm{\cdot-x^*} \ge \mf m + 3\sqrt{dh} + 3\sqrt{C_{\msf{PI}} + h}\,\cR_{3/2}(\tilde\mu_n^Y \mmid \pi^Y) + \sqrt{2\,(C_{\msf{PI}} + h)}\,\eta\Bigr\} \le 4\exp\Bigl(-\frac{\eta}{3}\Bigr)\,.
        \end{align*}
	After taking a union bound over $n=0,1,\dotsc,\Nprox -1$ and using~\eqref{eq:crude-renyi}, it shows that for $\delta \in (0,1)$, with probability at least $1-\delta$, it holds:
	\begin{align*}
		\max_{k=0,1,\dotsc,\Nprox-1}{\norm{\tilde Y_n - x^*}}
		&\lesssim \mf m + \sqrt{dh} + \sqrt{C_{\msf{PI}} + h} \, \Bigl( \Nprox \,\epsrgo^2 + \cR_3(\mu_0^X \mmid \pi^X) + \log \frac{\Nprox}{\delta}\Bigr)\,.
	\end{align*}
    This completes the proof.
\end{proof}

\subsubsection{Explicit versions of the results beyond strong log-concavity}\label{app:highacc:explicit:results}

Finally, we state and prove the more explicit versions of the results in \S\ref{ssec:hiacc-ext}.

\begin{theorem}[Faster low-accuracy sampling from log-concave targets, explicit]\label{thm:main-logconcave-explicit}
    Suppose that $\pi \propto \exp(-f)$, where $f$ is convex and $\beta$-smooth. There is an algorithm with randomized runtime that, given access to a minimizer $x^*$ of $f$ and to $N$ first-order queries for $f$, outputs a random point in $\R^d$ with law $\mu$ satisfying $\KL(\mu \mmid \pi) \le \eps^2$.
    Moreover, for any $\delta \in (0,1)$ with probability at least $1-\delta$, the number of queries made satisfies
    \begin{align*}
        N
        &\le \Otilde\biggl(\frac{d^{1/2}\,(d+\beta\mf m^2)}{\eps^2}\, \log^3 \log\frac{1}{\delta} \biggr)\,.
    \end{align*}
\end{theorem}

\begin{theorem}[Faster high-accuracy sampling from LSI targets, explicit]\label{thm:main-lsi-explicit}
	Suppose that $\pi \propto \exp(-f)$ satisfies $1/\alpha$-LSI and that $f$ is $\beta$-smooth. There is an algorithm with randomized runtime that, given access to a stationary point $x^*$ of $f$ and to $N$ first-order queries for $f$, outputs a random point in $\R^d$ with law $\mu$ satisfying
	$
		 \msf{d}(\mu \mmid \pi) \leq \eps
	$
	for any of the following metrics:
	\[
		\msf{d} \in \{ \TV, \sqrt{\KL}, \sqrt{\chi^2}, \sqrt{\alpha}\, W_2 \}\,.
	\]
    Moreover, for any $\delta \in (0,1)$ with probability at least $1-\delta$, the number of queries made satisfies
    \begin{align*}
        N
        &\le \Otilde\biggl(\kappa d^{1/2}\, \log^4 \max\Bigl\{\frac{1}{\varepsilon},\,\Delta,\, \mf m, \, \log \frac{1}{\delta}\Bigr\} \biggr)\,.
    \end{align*}
\end{theorem}

\begin{theorem}[Faster high-accuracy sampling from PI targets, explicit]\label{thm:main-poincare-explicit}
    Suppose that $\pi \propto \exp(-f)$ satisfies $1/\alpha$-PI and that $f$ is $\beta$-smooth. There is an algorithm with randomized runtime that, given access to a stationary point $x^*$ of $f$ and to $N$ first-order queries for $f$, outputs a random point in $\R^d$ with law $\mu$ satisfying
	$
		 \msf{d}(\mu \mmid \pi) \leq \eps
	$
	for any of the following metrics:
	\[
		\msf{d} \in \{ \TV, \sqrt{\KL}, \sqrt{\chi^2}, \sqrt{\alpha}\, W_2 \}\,.
	\]
    Moreover, for any $\delta \in (0,1)$ with probability at least $1-\delta$, the number of queries made satisfies
    \begin{align*}
        N
        &\le \widetilde O\biggl(\kappa d^{1/2}\,\max\Bigl\{\Delta, \,d, \,\log\frac{1}{\varepsilon}\Bigr\}\, \log^3 \max\Bigl\{\mf m,\,\log \frac{1}{\delta}\Bigr\} \biggr)\,.
    \end{align*}
\end{theorem}

\begin{proof}[Proofs of Theorems~\ref{thm:main-logconcave-explicit},~\ref{thm:main-lsi-explicit}, and~\ref{thm:main-poincare-explicit}]
    For the proofs of Theorems~\ref{thm:main-lsi-explicit} and~\ref{thm:main-poincare-explicit}, we can simply use the proofs of Theorems~\ref{thm:main-lsi} and~\ref{thm:main-poincare} given in Appendices~\ref{app:highacc:slc-main} and~\ref{app:high-acc:poincare} respectively, substituting in the R\'enyi initialization bound of Lemma~\ref{lem:renyi-gaussian-init} and the bound on the size of the proximal sampler iterates in Lemma~\ref{lem:size-prox-iterates}.

    For the proof of Theorem~\ref{thm:main-logconcave-explicit}, we again follow the proof of Theorem~\ref{thm:main-logconcave} in Appendix~\ref{app:high-acc:logconcave} and use Lemmas~\ref{lem:renyi-gaussian-init} and~\ref{lem:size-prox-iterates}, but with the following additional ingredients.
    
    First, we bound the Wasserstein distance at initialization.
We note that for $\mu_0^X = \cN(x^*, (2\beta)^{-1}I_d)$,
\begin{align*}
    W_2(\mu_0^X, \pi^X)
    &\le W_2(\mu_0^X, \delta_{x^*}) + W_2(\delta_{x^*}, \pi^X)
    \le \sqrt{\frac{d}{2\beta}} + \sqrt{\E_{\pi^X}[\norm{\cdot - x^*}^2]}
    \lesssim \sqrt{\frac{d}{\beta}} + \mf m\,,
\end{align*}
where the last inequality makes use of a reverse H\"older inequality for log-concave measures; see~\cite[Proposition A.5]{alonsogutierrezbastero2015kls}.

    Second, although we could invoke Lemma~\ref{lem:size-prox-iterates} directly, this would incur a dependence on the Poincar\'e constant of $\pi$. Although all log-concave measures indeed satisfy a PI\@, it is not always straightforward to estimate the Poincar\'e constant. Therefore, we note that by invoking the reverse H\"older inequality in~\cite[Proposition A.5]{alonsogutierrezbastero2015kls} rather than the concentration under a PI (Lemma~\ref{lem:pi-subexp}), it is possible to eliminate the dependence on $C_{\msf{PI}}$ in Lemma~\ref{lem:size-prox-iterates} when $\pi$ is log-concave, which gives the final statement of the result.
\end{proof}

\begin{remark}[Interpretation of the bounds]
    We remark that the parameters $\mf m$ and $\Delta$ have already appeared in prior sampling analyses such as~\cite{chenetal2022proximalsampler, chewi2021analysis}.
    For the reader's convenience, we pause to discuss the interpretation of the dependencies in the final bounds.
    First, we note that the dependence on the failure probability $\delta$ is polynomial in $\log\log(1/\delta)$, and is therefore negligible.

    The dependence on the first moment bound $\mf m$ is logarithmic in all of the results and hence typically negligible, with the exception of Theorem~\ref{thm:main-logconcave} in which case it is usually the dominant term.

    Finally, in the isoperimetric settings, there is additionally a dependence on the objective gap $\Delta\deq f(x^*) - \min f$, which measures the quality of the stationary point.
    If $\Delta \lesssim d$, which is realistic for many applications (it reflects the situation in which the user has some reasonable prior knowledge about the mode), then this is never the dominant term in the bounds; and for Theorem~\ref{thm:main-lsi}, the dependence on $\Delta$ is only logarithmic.

    Although these parameters have to be controlled for any given application, for interpretability we give simplified statements of the bounds: we assume $\Delta \lesssim d$, that $\mf m$ is polynomially bounded, and we omit logarithmic factors for simplicity.

    \begin{center}
        \begin{tabular}{cc}
            \textbf{assumptions} & \textbf{complexity} \\
            log-concavity and log-smoothness & $\beta d^{1/2} \mf m^2/\eps^2$ \\
            LSI and log-smoothness & $\kappa d^{1/2} \;\polylog(1/\eps)$ \\
            PI and log-smoothness & $\kappa d^{3/2} \; \polylog(1/\eps)$
        \end{tabular}
    \end{center}
\end{remark}